\documentclass[10pt]{article}
\usepackage{amsmath,amsthm,amsfonts,amssymb,mathrsfs}
\usepackage{tikz}
\usepackage{graphicx}
\usepackage{hyperref}
\hypersetup{hidelinks}
\usepackage{color}
\usepackage[T1]{fontenc}
\usepackage{amsthm}
\usepackage[margin=1in]{geometry}
\allowdisplaybreaks

\newtheorem{thm}{Theorem}[section]
\newtheorem{lem}[thm]{Lemma}
\newtheorem{cor}[thm]{Corollary}
\newtheorem{prop}[thm]{Proposition}
\newtheorem{defin}[thm]{Definition}
\newtheorem{rem}[thm]{Remark}
\newtheorem{hyp}[thm]{Hypothesis}

\newtheorem{nota}[thm]{Notation}

\numberwithin{equation}{section}

\title{CLT for the trace functional of the IDS of magnetic random Schr\"odinger operators}
\author{
Dhriti Ranjan Dolai\thanks{
Department of Mathematics, Indian Institute of Technology Dharwad,
Dharwad 580011, Karnataka, India.
Email: \texttt{dhriti@iitdh.ac.in}
}
\quad and \quad
Naveen Kumar\thanks{
Department of Mathematics, Indian Institute of Technology Dharwad,
Dharwad 580011, Karnataka, India.
Email: \texttt{222071001@iitdh.ac.in}
}
}
\date{}
\begin{document}
\maketitle

\begin{abstract}
\noindent We consider the existence of the integrated density of states (IDS) of the magnetic Schr\"{o}dinger operator with a random potential on the Hilbert space \( L^2(\mathbb{R}^d) \), as an analogue of the law of large numbers (LLN) for trace functionals. In this work, we establish an analogue of the central limit theorem (CLT), which describes the fluctuations of the trace functionals of the IDS, for a class of test functions denoted by \( C^1_{d,0}(\mathbb{R}) \). This class consists of real-valued, continuously differentiable functions on \( \mathbb{R} \) that decay at the rate \( O(|x|^{-m}) \) as \( |x| \to \infty \), where \( m > d + 1 \).
\end{abstract}
\section{Introduction}
Random Schr\"odinger operators are fundamental objects in the study of disordered quantum systems, such as unordered alloys, amorphous solids, and liquids. Their physical relevance was first highlighted in the seminal work of Anderson~\cite{apw} in the late 1950s. A central spectral quantity associated with these random operators is the integrated density of states (IDS), which represents the average counting function of eigenvalues per unit volume and is defined via a thermodynamic limit as the volume tends to infinity. In this work, we investigate the fluctuations of finite-volume approximations of the IDS around their expectation and show that, under a suitable normalization, these fluctuations converge in distribution to a Gaussian random variable.\\~\\
Here, we consider the random Schr\"{o}dinger operator \( H^\omega \) on the Hilbert space \( L^2(\mathbb{R}^d) \) in the presence of a magnetic field, defined by
\begin{equation}
\label{model}
H^\omega := H_A + V^\omega, \quad \omega \in \Omega, \quad \text{where} \quad H_A := (i\nabla + A)^2.
\end{equation}
Here, \( A \in \big(L^2_{\mathrm{loc}}(\mathbb{R}^d)\big)^d \) denotes the vector potential of a constant magnetic field. The operator \( H_A\), known as the \emph{free magnetic Laplacian}, is a positive, unbounded, self-adjoint operator on \( L^2(\mathbb{R}^d) \), defined via an associated quadratic form. A detailed description of the operator \( H_A \) and its domain is provided in Remark~\ref{fl-d-lp}.
We consider the alloy-type random potential $V^\omega$ defined as the bounded multiplication operator
\begin{equation}
\label{pntal}
\big(V^\omega \varphi\big)(x) = \sum_{n \in \mathbb{Z}^d} \omega_n u(x - n)\varphi(x), \quad \forall \, \varphi \in L^2(\mathbb{R}^d),
\end{equation}
where the single site potential \( u \in L^\infty(\mathbb{R}^d) \) is real-valued and compactly supported. The coefficients \( \{\omega_n\}_{n \in \mathbb{Z}^d} \) are independent and identically distributed (i.i.d.), bounded, real-valued random variables with a common distribution \( \mu \), referred to as the \emph{single site distribution} (SSD).
The underlying probability space is the product measure space \( (\Omega, \mathcal{B}_\Omega, \mathbb{P}) := \big(\mathbb{R}^{\mathbb{Z}^d}, \mathcal{B}_{\mathbb{R}^{\mathbb{Z}^d}}, \mathbb{P} \big) \), where \( \mathbb{P} = \bigotimes_{n \in \mathbb{Z}^d} \mu \) is constructed via Kolmogorov’s extension theorem. We denote elements by \( \omega = (\omega_n)_{n \in \mathbb{Z}^d} \in \Omega \).
It follows from \cite{WF} and \cite{TLMS} that \( H^\omega \) is self-adjoint on \( L^2(\mathbb{R}^d) \) for almost every (a.e.) \( \omega \in \Omega \). Since \( V^\omega \) is a bounded operator, we have $D(H^\omega) = D(H_A),$ for a.e. $\omega \in \Omega$.
Moreover, the mapping \( \omega \mapsto H^\omega \) is measurable with respect to
the product $\sigma$-algebra \( \mathcal{B}_\Omega \), so that \( \{H^\omega\}_{\omega \in \Omega} \) forms a measurable family of random operators. For further details, we refer to \cite{Wernr} and \cite{CL}. 
 \\
Let $O\subset \mathbb{R}^d$ be a bounded open set. We denote
$H^\omega_{O, N}$ and $H^\omega_{O, D}$ as the finite-volume restriction of $H^\omega$ to the Hilbert space $L^2(O)$ with the Neumann and Dirichlet boundary conditions, respectively.   The finite-volume restriction of $H^\omega$ is defined by
\begin{equation}
\label{NDrst}
H^\omega_{O, X}=H^X_{A,O}+V^\omega_O~~\text{on}~~L^2(O),~\omega\in\Omega~~\text{and}~~X=D,N.
\end{equation}
In the above, $H^X_{A,O},~X=D, N$ denotes the finite-volume restriction of the free magnetic Laplacian $H_A$ to the Hilbert space $L^2(O)$ with Dirichlet and Neumann boundary conditions.  
The operators $H^D_{A,O}$ and $H^N_{A,O}$ are unbounded, positive, self-adjoint on $L^2(O)$, a detailed description is given in Lemma \ref{dir-neu-lp}. The bounded operator $V^\omega_O$ is the restriction of $V^\omega$ to the Hilbert space $L^2(O)$, defined by $\big(V^\omega_O\varphi\big)(x)=\big(V^\omega \varphi\big)(x)~\forall~x\in O$ where $\varphi\in L^2(O)$.   Also $\{H^\omega_{O,X}\}_{\omega\in\Omega}$ forms a measurable collection of random operators,  \big(see \cite{KMcre},  \cite{TLMS} and \cite{CL}\big).
Since $V^\omega_O$ is a bounded operator, therefore $H^\omega_{O,X}$ is an unbounded self-adjoint operator on $L^2(O)$ with $D\big( H^\omega_{O,X} \big)=D\big( H^X_{A,O} \big)$ for a.e. $\omega$,  here $X=D,N$.\\
Before describing our result, we formally make assumptions on the vector potential $A$, single site potential $u$,  and the single site distribution (SSD) $\mu$.
\begin{hyp}\indent
\label{hypo}
\begin{enumerate}
\item The single site distribution (SSD) $\mu$ is compactly supported.
\item The single site potential $u\in L^\infty(\mathbb{R}^d)$ is real-valued and compactly supported and denote $u_n(x)=u(x-n)$, $n\in\mathbb{Z}^d$.
\item  The vector potential $A(x)=\big(A_1(x),A_2(x),\ldots,A_d(x)\big):\mathbb{R}^d\to\mathbb{R}^d$ is Borel-measurable and it is given by 
$$A(x)=\frac{1}{2}Bx,~x=(x_1,x_2,\ldots,x_d)\in\mathbb{R}^d.$$
Here $B=\big[ B_{i,j} \big]_{1\leq i,j\leq d}$ is a real skew-symmetric $d\times d$ matrix
representing the constant magnetic field.
\end{enumerate}
\end{hyp}
\begin{rem}
\label{bnd-f-pnt}
\noindent Under Hypothesis~\ref{hypo}, there exists a deterministic constant $\|V\|_\infty$ such that $\|V^\omega_O\| \leq \|V\|_\infty~\text{for a.e.}~\omega$ and every open set $O \subseteq \mathbb{R}^d$.
\end{rem}
\begin{rem}
Under Hypothesis~\ref{hypo}, it is clear that each component $A_k$ of the vector potential $A$ can be written as 
$A_k(x_1,x_2,\ldots, x_d)=\frac{1}{2}\displaystyle\sum_{j=1}^dx_jB_{k,j}$, for $k=1,2,\ldots, d$.  Now $A_k\in L^2_{\mathrm{loc}}(\mathbb{R}^d)$ is immediate.
\end{rem}
\begin{rem}
The random operator $H^\omega$ can be self-adjoint even when the compact support assumptions on the single site distribution (SSD) $\mu$ and the single site potential $u$ are relaxed. In this case, appropriate moment conditions on $\mu$ and integrability conditions on $u$ are required. Under these assumptions, the existence of the integrated density of states (IDS) has also been studied; we refer to \cite{WF}, \cite{TLMS}, and \cite{TLMS1} for further details.
It is also possible to obtain our central limit theorem~\eqref{clt-cn} when the SSD $\mu$ and the single site potential $u$ do not have compact support. This can be achieved through a two-step approximation. In the first step, we assume that $u$ has compact support while $\mu$ does not. The CLT for this model can then be derived from the CLT corresponding to compactly supported $u$ and a compactly truncated $\mu$. Once this is established, the CLT for both non-compactly supported $\mu$ and $u$ can be obtained as a limit of CLTs for non-compactly supported $\mu$ and compactly truncated $u$. In both approximation steps, a modified upper bound for the concentration inequality \eqref{vr-est} of the limiting variance is required, involving the norm $\|\tilde{f}\|_\infty$, where $\tilde{f}(x) = (1 + x^2)^{\lfloor (1+d)/2 \rfloor}f'(x)$. In both cases, Theorem~\ref{app-con-thm} will be useful.
\end{rem}
\begin{rem}
The existence of the integrated density of states (IDS) has also been investigated in the case where the magnetic field $A^\omega(x)$ is random; details can be found in \cite{U} and \cite{MH}. In this setting, one can also prove a central limit theorem using same method, provided that the operator $H^\omega$ is ergodic and that the random fields $\{A^\omega(x)\}_{x\in\mathbb{R}^d}$ and $\{V^\omega(x)\}_{x\in\mathbb{R}^d}$ remain statistically independent.
\end{rem}
\noindent Denote $m=(m_1,m_2,\ldots,m_d)\in\mathbb{Z}^d$. We define a family of the unitary operators $\{U_m\}_{m\in\mathbb{Z}^d}$ on $L^2(\mathbb{R}^d)$ as
\begin{equation}
\label{unop}
\big(U_m\varphi\big)(x)=e^{i\Psi_m(x)}\varphi(x-m)~~\forall~~\varphi\in L^2(\mathbb{R}^d),~\text{here}~i=\sqrt{-1},
\end{equation}
here $\Psi_m(x)=\displaystyle\frac{1}{2}\sum_{j,k=1}^d(m_j-x_j)B_{k,j}m_k,~x=(x_1,x_2,\ldots,x_d)\in\mathbb{R}^d$.  We also define an ergodic family $\{T_m\}_{m\in\mathbb{Z}^d}$ of measure-preserving transformations on the product probability space $\big(\Omega, \mathcal{B}_\Omega, \mathbb{P} \big)$ as
\begin{equation}
\label{mpt}
\big(T_m\omega\big)_n=\omega_{n-m},~~\omega=(\omega_n)_{n\in\mathbb{Z}^d}\in\Omega=\mathbb{R}^{\mathbb{Z}^d}
\end{equation}
$\{H^\omega\}_{\omega\in\Omega}$ is an ergodic family of random self-adjoint operators (unbounded) on $L^2(\mathbb{R}^d)$,  in the sense that $U_mH^\omega U_m^*=H^{T_m\omega}~~\forall~~\omega\in\Omega$ and $m\in\mathbb{Z}^d$.  We refer to \cite{Wernr} for details about ergodic operators (self-adjoint); see also \cite{CL}.   Since $H^\omega$ is an ergodic operator (self-adjoint), its spectral components are non-random subsets of the real line for a.e. $\omega$, and the discrete spectrum of $H^\omega$ is empty for a.e. $\omega$, details can be found in \cite{Wernr},  \cite{CL}, and \cite{ U}.  The spectrum of the finite-volume restriction $H^\omega_{O, X}$ \big(of $H^\omega$\big) is always purely discrete for a.e. $\omega$, $X=D,N,$  see \cite{TLMS} and \cite{RS4} for details.\\
First, we set a few notations to define the integrated density of states (IDS) of the ergodic operator $H^\omega$.  Let $\Lambda_L\subseteq\mathbb{R}^d$ be the cube (open) of side length $L$ centered at the origin i.e. 
\begin{equation}
\label{cbL}
\Lambda_L=\bigg\{x=(x_1,x_2,\ldots, x_d)\in\mathbb{R}^d:|x_i|< \frac{L}{2}  \bigg\},~~L\in\mathbb{N}.
\end{equation}
Now the integrated density of states (IDS) of $H^\omega$ can be described by the limit
\begin{equation}
\label{ids}
\lim_{L\to\infty}\frac{1}{\big|\Lambda_L  \big|}\operatorname{Tr}\!\big( f\big(H^\omega_{\Lambda_L,X}  \big) \big)=
\mathbb{E}\left[\operatorname{Tr}\!\left(\chi_{\Lambda_1}f\big( H^\omega\big)\chi_{\Lambda_1} \right)\right]~~a.e.~\omega,~~\forall~f\in C_c(\mathbb{R}).
\end{equation}
In the above, $C_c(\mathbb{R})$ denotes the set of all continuous functions defined on $\mathbb{R}$ with compact support. The existence of the limit in $(\ref{ids})$ is independent of the choice of boundary condition $X=D,N$. The distribution function $\mathcal{N}(\cdot)$ defined by $\mathcal{N}(x)=\mathbb{E}\left[\operatorname{Tr}\!\left(\chi_{\Lambda_1} E_{ H^\omega}(-\infty,x]\chi_{\Lambda_1}\right) \right],~x\in\mathbb{R}$ is known as the integrated density of states (IDS) of $H^\omega$.  \\
The existence of the limit $(\ref{ids})$ is a well-studied topic in the literature.
A proof of the existence of IDS in the sense of vague convergence (a.e. $\omega$) of random measures can be found in \cite{TLMS1},  the generalisation of \cite{PF}, which deals with the case when the magnetic field is absent. Use of the functional-analytic argument is given in \cite{MH}; it is first presented in \cite{KM} without a magnetic field. A different approach using Feynman-Kac(-It\^{o}) functional-integral representation of Schr\"{o}dinger semigroups can be found in \cite{U, BHL} and it goes back to \cite{pasl,  Nas} for $A=0$. For the uniqueness of the limit $(\ref{ids})$,  \big(its independence of the boundary condition $X=D,N$\big), we refer to \cite{TLMS1,  DIM, Na, HS} for non-zero magnetic field and \cite{Wernr, CL} for zero magnetic field. More details about IDS and, in general, random Schr\"{o}dinger operator is well documented in \cite{KB, His, CFKS, ves}.\\
Let us regard the limit in \((\ref{ids})\) as the analogue of the law of large numbers (LLN) for the trace functional $\operatorname{Tr}\!\big(f(H^\omega_{\Lambda_L,X})\big)$.
Our main goal is to establish an analogue of the central limit theorem (CLT) for \((\ref{ids})\) and to demonstrate its independence from the choice of boundary conditions \(X = D, N\). More explicitly, we aim to study the fluctuations of the trace functional $\operatorname{Tr}\!\big(f(H^\omega_{\Lambda_L,X})\big)$ around its mean as \(L \to \infty\).
\begin{rem}
\label{lbs}
In view of Remark~\ref{bnd-f-pnt}, we have  
$\sigma(H^\omega_{O,X}) \subseteq \big[-\|V\|_\infty,\infty\big)$,  and
$\sigma(H^\omega) \subseteq \big[-\|V\|_\infty,\infty\big)$ a.e. $\omega$,
where $O \subset \mathbb{R}^d$ and $X \in \{D,N\}$.
\end{rem}
\noindent We now define the class of test functions $C^1_{d,0}[-\|V\|_\infty,\infty)$, for which we will establish an analogue of the central limit theorem for the limit (\ref{ids}).
\begin{defin}
\label{defcl}
We say $f\in C^1_{d, 0}\big[-\|V\|_\infty,\infty\big)$, if $f$ is real-valued, continuously differentiable function on $\big[-\|V\|_\infty,\infty\big)$ and $|f(x)|=O(x^{-{m_1}})$, $|f'(x)|=O(x^{-{m_2}})$ as $x\to\infty$ for some $m_1>d+1$,  $m_2>d+1$.
\end{defin}
\begin{rem}
The test function $f$ needs to have sufficient decay at infinity in order for the operator $f\big(H^\omega_{\Lambda_L,X}  \big)$ to be trace class.
\end{rem}
\noindent Now for any real-valued Borel measurable function $f$ on $\mathbb{R}$ we define the random variable $Y_{f,X,L}(\omega)$ as 
\begin{equation}
\label{rv}
Y_{f,X,L}(\omega)=\bigg(\operatorname{Tr}\!\big( f\big(H^\omega_{\Lambda_L,X}  \big) \big)-\mathbb{E}
\big[\operatorname{Tr}\!\big( f\big(H^\omega_{\Lambda_L,X}  \big) \big) \big]\bigg),~~X=D,N.
\end{equation}
\noindent
For $X=D,N$, we denote by $\sigma^2_{f,X}$ the limiting variance of 
$Y_{f,X,L}$ as $L\to\infty$, defined by
\begin{align}
\label{lim-vr}
\sigma^2_{f,X}&=\lim_{L\to\infty}\frac{1}{|\Lambda_L|} \mathrm{Var}\big(Y_{f,X,L}  \big)=\lim_{L\to\infty}\frac{1}{|\Lambda_L|}\mathbb{E}\bigg[\operatorname{Tr}\!\big( f\big(H^\omega_{\Lambda_L,X}  \big) \big)-\mathbb{E}\big[\operatorname{Tr}\!\big( f\big(H^\omega_{\Lambda_L,X}  \big) \big) \big]\bigg]^2.
\end{align}
\begin{rem}
We show that for \( f \in C^1_{d, 0}\big[-\|V\|_\infty,\infty\big) \), the limit~\eqref{lim-vr} exists, and the limiting variance \(\sigma^2_{f,X}\) is independent of the choice of boundary conditions \(X = D, N\).
\end{rem}
\noindent
To express $\sigma^2_{f,X}$ explicitly, we introduce certain subsets of
$\mathbb{Z}^d$ and their associated $\sigma$-algebras.
For $n = (n_1, n_2,\ldots, n_d) \in \mathbb{Z}^d$, define
\[
\begin{aligned}
A_0^1 &= \{\, n \in \mathbb{Z}^d : n_1 \le 0 \,\}, \quad\quad
A_1^1 = \{n\in\mathbb{Z}^d: n_1\le 1\},\\
A_{(1,0)}^2 &= A_0^1
\cup \{n\in \mathbb{Z}^d : n_1\le 1,\ n_2\le 0\}, \quad
A_{(1,1)}^2 = A_0^1 \,\cup\, \{\, n \in \mathbb{Z}^d : n_1 \le 1,\, n_2 \le 1 \,\},\\
A_{(1,1,0)}^3 &= A_{(1,0)}^2 \,\cup\, \{\, n \in \mathbb{Z}^d : n_1 \le 1,\, n_2 \le 1,\, n_3 \le 0 \,\}, \\
A_{(1,1,1)}^3 &= A_{(1,0)}^2 \,\cup\, \{\, n \in \mathbb{Z}^d : n_1 \le 1,\, n_2 \le 1,\, n_3 \le 1 \,\}.
\end{aligned}
\]
and inductively for $d\ge 3$,
\begin{equation*}
\label{sgmalg}
\begin{split}
& A^d_{\underbrace{\scriptstyle (1,1, \ldots,1,0)}_{\scriptstyle d}}=A^{d-1}_{\underbrace{\scriptstyle (1,1,\ldots,1,0)}_{\scriptstyle d-1}}
\cup \big\{n \in\mathbb{Z}^d: n_1\leq 1, n_2\leq 1\ldots,n_{d-1}\leq 1, n_d\leq 0\big\},\\
& A^d_{\underbrace{\scriptstyle (1,1,\ldots,1,1)}_{\scriptstyle d}}=A^{d-1}_{\underbrace{\scriptstyle (1,1,\ldots,1,0)}_{\scriptstyle d-1}}
\cup \big\{n \in\mathbb{Z}^d: n_1\leq 1, n_2\leq 1,\ldots,n_{d-1}\leq 1, n_d\leq 1\big\}.
\end{split}
\end{equation*}
Now we define the $\sigma$-algebras associated with the above two subsets of $\mathbb{Z}^d$ as 
\begin{align}
\label{smag}
\begin{split}
 \mathcal{F}^d_{\vec{1}_d}&=
\sigma\bigg(\omega_n:n\in A^d_{\underbrace{\scriptstyle (1,1,\ldots,1,1)}_{\scriptstyle d}}  \bigg),~~~~\vec{1}_d= (\underbrace { 1,1,\ldots,1,1,1}_{\scriptstyle d})\in\mathbb{Z}^d,\\
\mathcal{F}^d_{\vec{1}_{d-1,0}}&=
\sigma\bigg(\omega_n:n\in A^d_{\underbrace{\scriptstyle (1,1,\ldots,1,0)}_{\scriptstyle d}}  \bigg),~~\vec{1}_{d-1,0}= ( \underbrace { 1,1,\ldots,1,1}_{\scriptstyle d-1},0)\in\mathbb{Z}^d.
\end{split}
\end{align}
For each \( k \in \mathbb{Z}^d \) and \( t \in [0,1] \), we define the modified random operator
\[
H^\omega\big|_{(\omega_k \to t\omega_k)}
= H_A + t\omega_k\,u(\cdot - k)
+ \sum_{\substack{n \in \mathbb{Z}^d \\ n \ne k}} \omega_n\,u(\cdot - n).
\]
Then, using the spectral theorem for self-adjoint operators, we define, for any measurable function \( g : \mathbb{R} \to \mathbb{C} \),
\begin{equation}
\label{mdfmdl}
g(H^\omega)_{(\omega_k \to t\omega_k)} :=
g\!\left(H^\omega\big|_{(\omega_k \to t\omega_k)}\right).
\end{equation}
\noindent Now we are ready to state the main result of this work.
\begin{thm}
\label{mnthm}
Let $H^\omega$ be as in \eqref{model} and assume Hypothesis~\ref{hypo}. 
Then, for each $f \in C^1_{d,0}\big[-\|V\|_\infty,\infty\big)$, we have the following convergence of random variables:
\begin{equation}
\label{clt-cn}
\frac{1}{|\Lambda_L|^{\frac{1}{2}}}
\bigg(
  \operatorname{Tr}\!\big( f(H^\omega_{\Lambda_L,X}) \big)
  - \mathbb{E}\!\big[\operatorname{Tr}\!\big( f(H^\omega_{\Lambda_L,X}) \big)\big]
\bigg)
\xrightarrow[L\to\infty]{\;\;\text{distribution}\;\;}
\mathcal{N}\!\big(0,\sigma^2_{f,X}\big),
\end{equation}
where $X \in \{D,N\}$ and $\mathcal{N}\!\big(0,\sigma^2_{f,X}\big)$ denotes the normal distribution with mean $0$ and variance $\sigma^2_{f,X}$, as defined in \eqref{lim-vr}. Moreover, for any $f \in C^1_{d,0}\big[-\|V\|_\infty,\infty\big)$, the limiting variance $\sigma^2_{f,X}$ is finite and independent of the boundary condition $X \in \{D,N\}$. In particular,
$\sigma^2_f := \sigma^2_{f,N} = \sigma^2_{f,D} < \infty$,
and its exact expression is given by
\begin{align}
\label{extlmv}
\sigma^2_f
&= \mathbb{E}\!\Bigg[\omega_{\vec{1}_d}\,
  \mathbb{E}\!\bigg(
    \int_0^1 
    \operatorname{Tr}\!\Big( u_{\vec{1}_d}(x)\,
    f'(H^\omega)_{(\omega_{\vec{1}_d}\to t\omega_{\vec{1}_d})} \Big)\,dt
  \;\Big|\; \mathcal{F}^d_{\vec{1}_d}  \bigg) \nonumber\\
&\qquad\qquad
- \mathbb{E}\!\bigg(
    \int_0^1 \omega_{\vec{1}_d}\,
    \operatorname{Tr}\!\Big( u_{\vec{1}_d}(x)\,
    f'(H^\omega)_{(\omega_{\vec{1}_d}\to t\omega_{\vec{1}_d})} \Big)\,dt
    \;\Big|\; \mathcal{F}^d_{\vec{1}_{d-1,0}}
  \bigg)
\Bigg]^2.
\end{align}
Finally, if $f \in C^1_{d,0}\big[-\|V\|_\infty, \infty\big)$ is strictly monotone and $u \ge 0$ or $u \le 0$ with $\|u\|_2 \neq 0$, then $\sigma_f^2 > 0$.
\end{thm}
\noindent To the best of our knowledge, this is the first work establishing a central limit theorem (CLT) for the integrated density of states (IDS) of a magnetic Schr\"odinger operator acting on $L^2(\mathbb{R}^d)$ with an alloy-type random potential. Moreover, this is the first result of its kind for any continuum random Schr\"odinger operator in dimension $d \geq 2$.
\\~\\
\noindent Up to now, the only known result concerning the central limit theorem (CLT) for the integrated density of states (IDS) of a continuous random Schr\"{o}dinger operator on \( L^2(\mathbb{R}) \) is due to Re\v{z}nikova~\cite{R2} (see also~\cite{R}). In~\cite{R2}, the author studied the one-dimensional Schr\"{o}dinger operator  $H_L = -\frac{d^2}{dt^2} + q(t,\omega)$, 
defined on \( L^2(-L, L) \) with classical boundary conditions. The random potential is given by \( q(t,\omega) = F(X_t) \), where \( X_t \) is a Brownian motion on the \( v \)-dimensional torus \( S^v \), and \( F \) is a smooth function on \( S^v \) satisfying \(\displaystyle \min_{x \in S^v} F(x) = 0 \). Let \( \mathcal{N}_L(\lambda) \) denote the number of eigenvalues of \( H_L \) below a given energy \( \lambda \in \mathbb{R} \). Re\v{z}nikova proved that the centered and normalized eigenvalue counting function $\frac{\mathcal{N}_L(\lambda) - 2L \mathcal{N}(\lambda)}{\sqrt{2L}}$ converges in distribution to a continuous Gaussian process \( \mathcal{N}^*(\lambda) \), whose finite-dimensional distributions are non-degenerate for \( \lambda > 0 \). Moreover, the limiting process \( \mathcal{N}^*(\lambda) \) exhibits locally independent increments. This result implies that the convergence in~\eqref{clt-cn} holds for the one-dimensional Schr\"{o}dinger operator on \( L^2(\mathbb{R}) \) with a Markov-type random potential \( F(X_t) \), when the test function is the indicator function \( f(x) = \chi_{(-\infty, \lambda]}(x) \), for any \( \lambda \in \mathbb{R} \). The fluctuations of the integrated density of states for the one-dimensional continuum model with a decaying random potential were also studied by Nakano in \cite{nfjsp}.\\~\\
Earlier progress has also been made in the study of the CLT for the IDS in the discrete setting, namely the Anderson model on $\ell^2(\mathbb{Z}^d)$. In the one-dimensional case, CLT results were obtained by Re\v{z}nikova~\cite{Rez} for indicator test functions, by Kirsch--Pastur~\cite{KP} for the test function $f(x) = (x - E)^{-1}$, where $E$ satisfies $\operatorname{dist}(E, \sigma(H^\omega)) > 0$, and by Pastur--Shcherbina~\cite{PLSM} for test functions $f$ with sufficiently high regularity.
In higher dimensions, results were established by Grinshpon--White~\cite{GYWJ} for polynomial test functions and later generalized by Dolai~\cite{DCLT} to differentiable functions with polynomial growth. In the case of one-dimensional models with decaying randomness, related CLTs were proved by Breuer--Grinshpon--White~\cite{BGW} for polynomial test functions and later extended by Mashiko--Marui--Maruyama--Nakano~\cite{MMMN} to real-analytic test functions.
In the context of random matrix theory, this type of CLT, known as the fluctuation of linear eigenvalue statistics, has been studied in great detail. For a comprehensive overview, we refer to the review by Forrester~\cite{fos} and the references therein.\\~\\
Here, we study the fluctuations of the eigenvalue counting measure for finite-volume
restrictions of $H_\omega$ as the volume tends to infinity. These fluctuations can be
viewed as macroscopic-scale fluctuations of the eigenvalues. In contrast, substantial
attention has previously been devoted to microscopic (local) fluctuations of eigenvalues,
which are studied via the convergence of point processes associated with suitably rescaled
finite-volume spectra as the volume increases.
In the discrete setting, for microscopic spectral statistics of i.i.d.\ stationary
potentials, we refer to the works (and references therein) of Minami~\cite{Min96},
Klopp~\cite{Klo13}, Germinet--Klopp~\cite{GK14}, and Aizenman--Warzel~\cite{AMSW}.
See also the work of Cannizzaro--Labb\'e--Zuijlen~\cite{CLZ}, which considers correlated
Gaussian potentials. For one-dimensional models with decaying randomness, see
Kritchevski--Valk\'o--Vir\'ag~\cite{KVV} and Dolai--Mallick~\cite{MD2019}.
In the continuum setting, related microscopic results include the work of
Molchanov~\cite{Mol81} in one dimension; Dietlein--Elgart~\cite{DE21},
Hislop--Kirsch--Krishna~\cite{HKK20}, and Hislop--Krishna~\cite{HK15} for operators on
$L^2(\mathbb{R}^d)$; Nakano~\cite{Nak14} and Kotani--Nakano~\cite{KN} for one-dimensional
models with decaying randomness; and Dumaz--Labb\'e~\cite{DL20, DL24} for one-dimensional
models with white-noise potentials.
In~\cite{Min96, Klo13, AMSW, Mol81, DE21, HKK20, KN, CLZ, DL20, DL24}, the limiting point
process is Poisson. In~\cite{KVV}, the limiting process is the $\mathrm{sine}_\beta$
process. In~\cite{Nak14}, the limit is either a clock process or the circular-$\beta$
ensemble, depending on the decay rate, while in~\cite{MD2019} the limit is a clock
process. Finally, in~\cite{HK15}, the limiting process is a compound Poisson process.
\\~\\
To prove a central limit theorem for the integrated density of states in the discrete setting, one considers quantities of the form $\operatorname{Tr}\, P\!\left(H^\omega_\Lambda\right)$,
where $P$ is a polynomial.
In this case, $\mathrm{Tr}\, P(H^\omega_\Lambda)$ can be evaluated by counting lattice paths, where $H^\omega_\Lambda$ denotes the finite-volume restriction of $H^\omega$. 
In one-dimensional continuous models, the analysis can instead be carried out using Prüfer phase methods. However, in higher-dimensional continuum models, neither path-counting methods nor one-dimensional Prüfer techniques are applicable. 
It is also important to note that, in the discrete setting, the finite-volume restriction $H^\omega_\Lambda$ acts on a finite-dimensional Hilbert space, whereas in the continuum case it corresponds to an unbounded random operator on the infinite-dimensional space $L^2(\Lambda)$. 
Consequently, the techniques required to establish a central limit theorem (CLT) for trace functionals associated with the integrated density of states (IDS) in the continuum setting differ fundamentally from those used in the discrete case. 
In particular, discrete methods do not extend to the continuum setting; it demands new methods based on probabilistic techniques.
In this work, we introduce a novel framework to prove the CLT for continuum random Schrödinger operators with magnetic potentials acting on the Hilbert space $L^{2}(\mathbb{R}^{d})$. 
Our approach does not rely on any assumptions concerning localization or other spectral properties.
\\~\\
We outline the proof of Theorem~\ref{mnthm}, focusing on the main steps that require probabilistically demanding methods. Rather than treating general test functions directly, we first restrict attention to a special class of functions for which the relevant operator expressions can be analyzed precisely. Specifically, we consider test functions of the form \( P(x) = (x - E)^{-m} \sum_{k=0}^p a_k (x - E)^{-k},\) where \( E < -\|V\|_\infty \), \( m > d+1 \), \( p \in \mathbb{N} \cup \{0\} \), and \( a_k \in \mathbb{R} \). These functions are Laurent polynomials on the interval \( [-\|V\|_\infty, \infty) \) and, via functional calculus, correspond to powers of the resolvent. This representation allows us to express the relevant quantities in terms of resolvents and to control them despite the unboundedness of the operator.
We begin by studying the fluctuations of the resolvent-power trace
\( \displaystyle\operatorname{Tr}((H^\omega_{\Lambda_L, D} - E)^{-m}),\)
with \( m > d+1 \), around its mean. Instead of analyzing this trace directly on the large box \( \Lambda_L \), we decompose \( \Lambda_L \) into smaller annular regions \( \{\Lambda_{L,k}\}_k \), whose sizes depend explicitly on \( k \) and \( L \); see~\eqref{boxes}. We then relate the trace on \( \Lambda_L \) to the sum of the corresponding traces over these smaller annular regions,
\( \operatorname{Tr}((H^\omega_{\Lambda_{L,k}, D} - E)^{-m}) \), after normalization by \( |\Lambda_L|^{-1/2} \). This equivalence is established in Proposition~\ref{e-q-cn} and is not a straightforward calculation. In fact, Proposition~\ref{e-q-cn} shows that
\( \operatorname{Var}\big(
\operatorname{Tr}((H^\omega_{\Lambda_L, D} - E)^{-m})
- \sum_k \operatorname{Tr}((H^\omega_{\Lambda_{L,k}, D} - E)^{-m})
\big)\le \mathcal{O}(|\partial \Lambda_{L}|).\)
The proof requires an estimate of the form
\( \displaystyle
\mathbb{E}\!\left[
\operatorname{Tr}(\chi_F (H^\omega_{\Lambda_L, D} - E)^{-m})
- \operatorname{Tr}(\chi_F (H^\omega_{\Lambda_{L,k}, D} - E)^{-m})
\right]^2
\le C\, \mathrm{e}^{-\beta\, \operatorname{dist}(F, \partial \Lambda_{L,k})},
\)
for all \( F \subset \Lambda_{L,k}^\circ \subset \Lambda_L\). This estimate relies on a preliminary result stated in Proposition~\ref{ex-dy-tr}.
The annular regions \( \{\Lambda_{L,k}\}_k \) are divided into two categories: large annular regions \( \{\Lambda^B_{L,k}\}_k \) and small annular regions \( \{\Lambda^S_{L,k}\}_k \). They are arranged so that each small annular region lies between two large annular regions, and vice versa. Moreover, the volume of a small annular region \( \Lambda^S_{k,L} \) grows at a lower order than that of a large annular region. The precise construction of these regions is given in~\eqref{boxes} and~\eqref{vl-est-box}. We show that the cumulative traces over the small annular regions vanish in quadratic mean as \( L \to \infty \), after normalization by \( |\Lambda_L|^{-1/2} \); see Corollary~\ref{v-er-0}. Consequently, the asymptotic distribution of the trace on \( \Lambda_L \) coincides with that of
\( \sum_k \operatorname{Tr}((H^\omega_{\Lambda^B_{L,k}, D} - E)^{-m}),\) after normalization. Furthermore, the traces over different large annular regions, which are not identically distributed since the annular regions have different sizes, become asymptotically independent as \( L \to \infty \). This is proved in Proposition~\ref{wk-eq-cn}.
A key ingredient in the analysis is control of the growth of
\( \displaystyle\mathbb{E}[\operatorname{Tr}((H^\omega_{\Lambda, D} - E)^{-m})]\)
as the size of the box \( \Lambda \) increases. This estimate is established in Corollary~\ref{2-4-pol}. The argument relies on moment bounds for centered traces obtained via martingale techniques and the Burkholder inequality, proved in Proposition~\ref{bdd-2-4}. These results yield a central limit theorem (CLT) for Laurent-polynomial test functions by establishing a CLT for the sum of traces over large annular regions; see Proposition~\ref{wk-eq-cn}, which uses a classical CLT for triangular arrays.The extension to general test functions
\( f \in C^1_{d,0}([-\|V\|_\infty, \infty)) \)
is achieved by approximation. One constructs a sequence of Laurent polynomials \( \{P_n\}_n \) such that \( \displaystyle
\lim_{n \to \infty} \lim_{L \to \infty}
|\Lambda_L|^{-1} \operatorname{Var}(Y_{(f - P_n), D, L}) = 0.\)
Establishing this double limit requires uniform control of the variance functional and leads to a concentration estimate showing that the limiting variance $\sigma_f^2$ is bounded in terms of $\|\tilde f\|_\infty$, namely $\sigma_f^2 \le C \|\tilde f\|_\infty^2$, where $\tilde f(x) = (x - E)^{1 + \lfloor d/2 \rfloor} f'(x)$ with $E < -\|V\|_\infty$. This result is proved in Proposition~\ref{bdd-var-f}. The construction of the Laurent polynomials \( \{P_n\} \) is described in Remark~\ref{appr-fun-pol}, from which it follows that \( \|P_n - \tilde{f}\|_\infty \to 0 \).
The existence of the limiting variance \( \sigma_f^2 \) is proved separately for general functions \( f \) in Lemma~\ref{lmvrfrf} and for resolvent powers in Lemma~\ref{exst-vr-resl}. These proofs rely on martingale difference methods and the ergodicity of the random operator \( H^\omega \). They also require a derivative formula for traces of the form \( \operatorname{Tr}(f(T_\lambda)) \), where \( T_\lambda = T + \lambda K \). This part of the proof, which establishes both the existence of \( \sigma_f^2 \) and its explicit representation, is conceptually involved and requires careful analysis.
With these ingredients in place, we obtain the CLT for the normalized variables
\( \displaystyle\{|\Lambda_L|^{-1/2} Y_{f, D, L}\}_L\)
in Lemma~\ref{nr-appr}. The limiting variance admits an explicit representation in terms of the infinite-volume operator \( H^\omega \) and the derivative \( f' \); see~\eqref{lm-f-lv}. Using this representation, we show in Lemma~\ref{ps-vr} that \( \sigma_f^2 > 0 \) for strictly monotone functions \( f \). Finally, we show that the boundary conditions do not affect the limiting fluctuations by proving that \( \displaystyle\lim_{L \to \infty}
|\Lambda_L|^{-1} \operatorname{Var}(Y_{f, D, L} - Y_{f, N, L}) = 0.\)
To establish this limit for general \( f \), we first prove that for the function \( f(x) = (x - E)^{-m} \), \( \displaystyle
\operatorname{Var}(Y_{f, D, L} - Y_{f, N, L})
\le \mathcal{O}(|\partial \Lambda_L|).\)
This estimate is obtained using the bound
\( \displaystyle
\mathbb{E}\!\left[
\operatorname{Tr}(\chi_F (H^\omega_{\Lambda_L, D} - E)^{-m})
- \operatorname{Tr}(\chi_F (H^\omega_{\Lambda_L, N} - E)^{-m})
\right]^2
\le C\, \mathrm{e}^{-\beta\, \operatorname{dist}(F, \partial \Lambda_L)},
\)
valid for all \( F \subset \Lambda_L^{\circ} \). This result is proved in Lemma~\ref{trn=d} and relies on a preliminary estimate given in Proposition~\ref{df-n-d-in}.\\~\\
We divide the proof into three parts. In the first part, we present several preliminary results. In the second part, we prove the central limit theorem for the case where the test function is a Laurent polynomial. The third part addresses the case of a general test function \( f \in C^1_{d,0}([- \|V\|_\infty, \infty)) \).
In the appendix, we collect several results from probability theory in the form required for our arguments. We also present a number of results concerning the calculus of magnetic Schr\"odinger operators and their finite-volume restrictions. All results included in the appendix are used in the proof of our main theorem.
\section{Some preliminary results}
In this section, we estimate the moments of the random variable $\operatorname{Tr}\!\left( (H^\omega_{O,X} - E)^{-m} \right)$,  for \(m > \frac{d}{2}\) and a bounded open set \(O \subset \mathbb{R}^d\). 
We also provide several identities describing the difference between two resolvent operators, 
including formulas for the corresponding trace differences.\\~\\
Let the random potential \( V^\omega \) be defined as in~(\ref{pntal}). 
For any bounded open set \( O \subset \mathbb{R}^d \), define the truncated potential 
$V^\omega_O := \chi_O\, V^\omega$.  
For each \( n \in \mathbb{Z}^d \), set \( u_n(x) := u(x - n) \), where \( u \) is the single-site potential given in~(\ref{pntal}). 
Define the index set \( B_O \subset \mathbb{Z}^d \) by  
\begin{equation}
\label{ind-set}
B_O := \{ n \in \mathbb{Z}^d : \operatorname{supp}(u_n) \cap O \neq \emptyset \}.
\end{equation}
Since \( u \) is compactly supported in \( \mathbb{R}^d \), it follows that  $\#B_O = \mathcal{O}(|O|)$,  where \( |O| \) denotes the Lebesgue measure of \( O \subset \mathbb{R}^d \).
Assume that the single site distribution (SSD) \( \mu \) has bounded support and that \( u \) is bounded and compactly supported. 
Then there exists a deterministic constant \( \|V\|_\infty \) such that  
\begin{equation}
\label{sup-bd-ptn}
\|V^\omega_O\|_\infty \leq \|V\|_\infty~~
a.e. ~~ \omega ~ ~\forall~~O \subseteq \mathbb{R}^d .
\end{equation}
For any energy \( E < -\|V\|_\infty \) and any open set \( O \subset \mathbb{R}^d \), we define the centered random variable \( Y_{X,O} \) on the probability space \( \big(\Omega, \mathcal{B}_\Omega, \mathbb{P} \big) \) by  
\begin{equation}
\label{rv-O}
Y_{X,O}(\omega) := \operatorname{Tr}\big( (H^X_{A,O} + V^\omega_O - E)^{-m} \big) 
- \mathbb{E}\left[ \operatorname{Tr}\big( (H^X_{A,O} + V^\omega_O - E)^{-m} \big) \right],
\end{equation}  
where \( m > \frac{d}{2} \) and \( H^X_{A,O} \) denotes the magnetic Laplacian on \( O \) with boundary condition \( X \in \{D, N\} \) (Dirichlet or Neumann).\\
We now establish moment bounds for the centered random variable \( Y_{X,O} \) defined in~\eqref{rv-O}. These estimates are essential for proving concentration results and central limit theorems for trace functionals of the integrated density of states (IDS). The following proposition shows that the second and fourth moments of \( Y_{X,O} \) grow at most linearly and quadratically with the volume of the domain, respectively.
\begin{prop}
\label{bdd-2-4}
Let \( O \subset \mathbb{R}^d \) be a bounded open set, and let \( H^X_{A,O} \) be as in Lemma~\ref{dir-neu-lp}. 
Assume that the single site distribution (SSD) \( \mu \) has bounded support, and that \( u \in L^\infty(O) \) is a real-valued, compactly supported function. 
Then, for the centered random variable \( Y_{X,O} \) defined in~\eqref{rv-O}, the following bounds hold for the second and fourth moments:
\begin{equation}
\label{2nd-4th}
\mathbb{E}\!\left[ |Y_{X,O}|^2 \right] \le C |O|, \qquad
\mathbb{E}\!\left[ |Y_{X,O}|^4 \right] \le C |O|^2,
\end{equation}
where \( C \) is a positive constant independent of \( O \), and \( |O| \) denotes the Lebesgue measure of \( O \subset \mathbb{R}^d \).
\end{prop}
\begin{proof}
To estimate the moments of \( Y_{X,O} \), we apply a martingale approach. Observe that \( Y_{X,O} \) depends only on the collection of random variables \( \{\omega_n\}_{n \in B_O} \). Fix an enumeration \( \{\omega_{n_j}\}_{j=1}^{\# B_O} \) of the set \( \{\omega_n\}_{n \in B_O} \).
Define the filtration \( \{\mathcal{G}_k\}_{k=0}^{\# B_O} \) of $\sigma$-algebras by
$\mathcal{G}_k = \sigma\left( \omega_{n_j} : j \leq k \right),  \mathcal{G}_0 = \{\emptyset, \Omega\}$.
Then \( \left\{ \mathbb{E}\big( Y_{X,O} \mid \mathcal{G}_k \big) \right\}_{k=0}^{\# B_O} \) forms a Doob martingale. We express \( Y_{X,O} \) as the sum of its martingale differences:
\begin{equation}
\label{sm-mrt-dif}
Y_{X,O} = \sum_{k=1}^{\# B_O} \big( \mathbb{E}( Y_{X,O} \mid \mathcal{G}_k ) - \mathbb{E}( Y_{X,O} \mid \mathcal{G}_{k-1} ) \big).
\end{equation}
We now estimate each martingale difference in the above sum. Recall that
\begin{align}
\label{mrt-dif-est}
&\mathbb{E}( Y_{X,O} \mid \mathcal{G}_k ) - \mathbb{E}(Y_{X,O} \mid \mathcal{G}_{k-1} ) \nonumber\\
&= \mathbb{E}\left( \operatorname{Tr}\left( (H^X_{A,O} + V^\omega_O - E)^{-m} \right) \mid \mathcal{G}_k \right) - \mathbb{E}\left( \operatorname{Tr}\left( (H^X_{A,O} + V^\omega_O - E)^{-m} \right) \mid \mathcal{G}_{k-1} \right) \nonumber\\
&=\mathbb{E}\left( \operatorname{Tr}\left( (H^X_{A,O} + V^\omega_O - E)^{-m} \right) \mid \mathcal{G}_k \right)-\mathbb{E}\left( \operatorname{Tr}\left( (H^X_{A,O} + V^\omega_O - E)^{-m} \right)_{(\omega_{n_k}=0) } \mid \mathcal{G}_k \right)\nonumber\\
&+\mathbb{E}\left( \operatorname{Tr}\left( (H^X_{A,O} + V^\omega_O - E)^{-m} \right)_{(\omega_{n_k}=0) }\mid \mathcal{G}_{k-1} \right)- \mathbb{E}\left( \operatorname{Tr}\left( (H^X_{A,O} + V^\omega_O - E)^{-m} \right) \mid \mathcal{G}_{k-1} \right)
\nonumber\\
&= \mathbb{E}\left( \int_0^1 \frac{d}{dt} \operatorname{Tr}\left( \big(H^X_{A,O} + V^\omega_O - E\big)^{-m} \right)_{(\omega_{n_k} \to t \omega_{n_k})} dt \, \bigg| \, \mathcal{G}_k \right) \nonumber\\
&\qquad \qquad- \mathbb{E}\left( \int_0^1 \frac{d}{dt} \operatorname{Tr}\left( \big(H^X_{A,O} + V^\omega_O - E\big)^{-m} \right)_{(\omega_{n_k} \to t \omega_{n_k})} dt \, \bigg| \, \mathcal{G}_{k-1} \right).
\end{align}
Define \( u_{n_k}(x) := u(x - n_k) \), and note the shorthand
\begin{align}
\label{rpl}
&\left( H^X_{A,O} + V^\omega_O - E \right)^{-m}_{(\omega_{n_k} \to t\omega_{n_k})}:= \left( H^X_{A,O} + V^\omega_O - \omega_{n_k} \chi_O u_{n_k}(x) + t \omega_{n_k} \chi_O u_{n_k}(x) - E \right)^{-m}.
\end{align}
Applying the derivative formula for the trace of the resolvent, as given in Corollary~\ref{drv-tr-rsl}, we obtain:
\begin{align}
\label{est-by-drv}
&\mathbb{E}( Y_{X,O} \mid \mathcal{G}_k ) - \mathbb{E}( Y_{X,O} \mid \mathcal{G}_{k-1} ) \\
&= -m \, \mathbb{E}\left( \int_0^1 \omega_{n_k} \operatorname{Tr} \left( \chi_O u_{n_k}(x) \left( H^X_{A,O} + V^\omega_O - E \right)^{-m-1}_{(\omega_{n_k} \to t \omega_{n_k})} \right) dt \, \bigg| \, \mathcal{G}_k \right)\nonumber \\
&\quad \quad+ m \, \mathbb{E}\left( \int_0^1 \omega_{n_k} \operatorname{Tr} \left( \chi_O u_{n_k}(x) \left( H^X_{A,O} + V^\omega_O - E \right)^{-m-1}_{(\omega_{n_k} \to t \omega_{n_k})} \right) dt \, \bigg| \, \mathcal{G}_{k-1} \right)\nonumber.
\end{align}
Since the single site distribution (SSD) $\mu$ has compact support (i.e., $|\omega_n| \leq C$ almost surely) and the function $u$ is both bounded and compactly supported, we may apply Proposition~\ref{est-tr-chf} and Corollary~\ref{cr-tr-rst} to the identity~\eqref{est-by-drv} to obtain the following estimate:
\begin{equation}
\label{con-dif}
\left| \mathbb{E}( Y_{X,O} \mid \mathcal{G}_k ) - \mathbb{E}(Y_{X,O} \mid \mathcal{G}_{k-1} ) \right| \leq C 2m \|u\|_\infty |\text{supp}(u)|,
\end{equation}
where \( |\text{supp}(u)| \) denotes the Lebesgue measure of the support of \( u \).
Substituting this into \eqref{sm-mrt-dif}, we obtain:
\begin{align*}
\mathbb{E}\big[ |Y_{X,O}|^2 \big]
&= \sum_{k=1}^{\# B_O} \mathbb{E}\big[ \mathbb{E}(Y_{X,O} \mid \mathcal{G}_k ) - \mathbb{E}( Y_{X,O} \mid \mathcal{G}_{k-1} ) \big]^2 \leq {\# B_O} \big( C 2m \|u\|_\infty |\text{supp}(u)| \big)^2.
\end{align*}
Since \( \#B_O = \mathcal{O}(|O|) \), this yields the first part of \eqref{2nd-4th}.
To obtain the second part, we apply a Burkholder inequality as in (\ref{binq}), together with (\ref{con-dif}):
\begin{align*}
\mathbb{E}\big[ |Y_{X,O}|^4 \big]
&\leq C_4 \, \mathbb{E}\left[ \sum_{k=1}^{\# B_O} \left( \mathbb{E}( Y_{X,O} \mid \mathcal{G}_k ) - \mathbb{E}( Y_{X,O} \mid \mathcal{G}_{k-1} ) \right)^2 \right]^2 \\
&\leq C_4 (\# B_O)^2 \big( C 2m \|u\|_\infty |\text{supp}(u)| \big)^4.
\end{align*}
Again, since \(\#B_O = \mathcal{O}(|O|) \), we conclude the proof of \eqref{2nd-4th}.
\end{proof}
\noindent We next consider the effect of domain extension on the resolvents of magnetic Schr\"{o}dinger operators with either Dirichlet or Neumann boundary conditions. Specifically, we compare the resolvent powers of an operator defined on a bounded domain with those of its natural extension to a larger open set. The following proposition gives an explicit identity for the difference between these resolvent powers, expressed in terms of localized operators and cut-off functions. This identity is useful in various spectral and semiclassical analyses, especially when examining the impact of boundary conditions or domain truncations.
\begin{prop}
\label{idnty-ne-dir-rst}
Let \( O \) and \( \tilde{O} \) be open subsets of \( \mathbb{R}^d \) such that \( O \subset \tilde{O} \subseteq \mathbb{R}^d \), with \( O \) is bounded. Consider the operator \( H^X_{A,O} \) for \( X = D, N \), as defined in Lemma~\ref{dir-neu-lp}, and let \( \phi \in C^\infty(O) \) satisfy \( \operatorname{supp}(1 - \phi) \subset\subset O \). Denote \( \operatorname{supp}(\phi) = S \subseteq O \). Suppose further that \( V \in L^\infty(\tilde{O}) \) is real-valued function. Then, for every \( E < -\|V\|_\infty \), the following identity holds as an operator on \( L^2(O) \):
\begin{align}
\label{rst-dif}
&\chi_O \left( H^X_{A,\tilde{O}} + V - E \right)^{-m} \chi_O^*
- \left( H^X_{A,O} + V - E \right)^{-m} \nonumber\\
& =
\sum_{j=0}^{m-1} \Bigg[
\chi_O \left( H^X_{A,\tilde{O}} + V - E \right)^{j - m} \chi_S^*
\left( H^X_{A,O} + V - E \right) \phi \left( H^X_{A,O} + V - E \right)^{-j - 1} \nonumber\\
&\quad\quad~~~
- \chi_O \left( H^X_{A,\tilde{O}} + V - E \right)^{j - m + 1} \chi_S^*~
\phi \left( H^X_{A,O} + V - E \right)^{-j - 1}
\Bigg].
\end{align}
Here, \( \chi_O^* : L^2(O) \to L^2(\tilde{O}) \) and \( \chi_S^* : L^2(S) \to L^2(\tilde{O}) \) are the extension operators defined in Definition~\ref{cn-ext-rst}. In the special case \( \tilde{O} = \mathbb{R}^d \), we have \( H^X_{A,\tilde{O}} = H_A \), as explained in Remark~\ref{fl-d-lp}.
\end{prop}
\begin{proof}
Let \( f, g \in L^2(O) \). We begin by observing:
\begin{align}
& \left\langle f, \chi_O \left( H^X_{A,\tilde{O}} + V - E \right)^{-m} \chi_O^*g \right \rangle
- \left\langle f, \left( H^X_{A,O} + V - E \right)^{-m} g \right \rangle \nonumber\\
&\qquad= \left\langle \chi_O \left( H^X_{A,\tilde{O}} + V - E \right)^{-m} \chi_O^*f, g \right \rangle
- \left\langle f, \left( H^X_{A,O} + V - E \right)^{-m} g \right \rangle \nonumber\\
&\qquad= \sum_{j=0}^{m-1} \Big[
\left\langle \chi_O \left( H^X_{A,\tilde{O}} + V - E \right)^{j - m} \chi_O^*f, \left( H^X_{A,O} + V - E \right)^{-j} g \right\rangle \nonumber\\
&\qquad\quad - \left\langle \chi_O \left( H^X_{A,\tilde{O}} + V - E \right)^{j - m + 1} \chi_O^*f, \left( H^X_{A,O} + V - E \right)^{-j - 1} g \right\rangle
\Big].
\end{align}
Define \( f_j := \left( H^X_{A,\tilde{O}} + V - E \right)^{j - m} \chi_O^*f \) and \( g_j := \left( H^X_{A,O} + V - E \right)^{-j - 1} g \). Then the above becomes:
\begin{align}
\label{rl-sm-bg}
&=\sum_{j=0}^{m-1} \left[
\left\langle \chi_O f_j, \left( H^X_{A,O} + V - E \right) g_j \right\rangle
- \left\langle \chi_O \left( H^X_{A,\tilde{O}} + V - E \right) f_j, g_j \right\rangle
\right] \nonumber\\
&=\sum_{j=0}^{m-1} \bigg[
\left\langle \chi_O f_j, \left( H^X_{A,O} + V - E \right) \phi g_j \right\rangle
- \left\langle \chi_O \left( H^X_{A,\tilde{O}} + V - E \right) f_j, \phi g_j \right\rangle \nonumber\\
&\quad+\bigg\{
\left\langle \chi_O f_j, \left( H^X_{A,O} + V - E \right)(1 - \phi) g_j \right\rangle - \left\langle \chi_O \left( H^X_{A,\tilde{O}} + V - E \right) f_j, (1 - \phi) g_j \right\rangle\bigg\}
\bigg].
\end{align}
Since \( g_j \in D(H^X_{A,O}) \) and \( 1 - \phi \in C^\infty(O) \) with \( \operatorname{supp}(1 - \phi) \subset\subset O \), Lemma~\ref{mult-dom} and~Lemma~\ref{mult-dom-neu} imply that \( (1 - \phi) g_j \in D(H^X_{A,O}) \). Moreover, since \( \operatorname{supp}((1 - \phi)g_j) \subset\subset O \), Lemma~\ref{ext-dir} and Lemma~\ref{ext-neu} yield:
\begin{equation}
\label{u-ext}
\chi_O^*(1 - \phi)g_j \in D(H^X_{A,\tilde{O}}), \quad
H^X_{A,\tilde{O}}\big(\chi_O^*(1 - \phi)g_j\big) = \chi_O^* H^X_{A,O} \big((1 - \phi)g_j\big).
\end{equation}
Substituting this into~\eqref{rl-sm-bg}, the second bracketed term cancels, yielding:
\begin{align}
&= \sum_{j=0}^{m-1} \Big[
\left\langle \chi_O f_j, \left( H^X_{A,O} + V - E \right) \phi g_j \right\rangle
- \left\langle \chi_O \left( H^X_{A,\tilde{O}} + V - E \right) f_j, \phi g_j \right\rangle
\Big] \nonumber\\
&= \sum_{j=0}^{m-1} \left[
\left\langle f, \chi_O \left( H^X_{A,\tilde{O}} + V - E \right)^{j - m} \chi_O^* \left( H^X_{A,O} + V - E \right) \phi \left( H^X_{A,O} + V - E \right)^{-j - 1} g \right\rangle \right. \nonumber\\
&\quad \qquad\qquad- \left.
\left\langle f, \chi_O \left( H^X_{A,\tilde{O}} + V - E \right)^{j - m + 1} \chi_O^* \phi \left( H^X_{A,O} + V - E \right)^{-j - 1} g \right\rangle\nonumber
\right].
\end{align}
Finally, since \( \operatorname{supp}(\phi) = S \), and by Lemmas~\ref{mult-dom} and~\ref{mult-dom-neu}, we have
\[
\operatorname{supp} \left( \left( H^X_{A,O} + V - E \right) \phi \left( H^X_{A,O} + V - E \right)^{-j - 1} g \right) \subseteq S,
\]
so we may replace \( \chi_O^* \) with \( \chi_S^* \) in the above expression, which completes the proof.
\end{proof}
\noindent We now compare the resolvents of magnetic Schr\"{o}dinger operators with Dirichlet and Neumann boundary conditions, defined on the same open domain. The difference between the \( m \)-th powers of these resolvents plays a central role in various spectral and functional-analytic estimates. The following proposition provides an explicit identity for this difference, expressed in terms of localized operators and cut-off functions. Such identities are useful, for instance, in studying the asymptotic behavior of traces or spectral shift functions under boundary perturbations.
\begin{prop}
For an open set \( O \subset \mathbb{R}^d \), let \( H_{A,O}^X \), for \( X = D, N \), be as defined in Lemma \ref{dir-neu-lp}. Assume that \( \phi \in C^\infty(O) \) satisfies \( \operatorname{supp}(1 - \phi) \subset\subset O \).  Denote $\operatorname{supp}(\phi)=S\subseteq O$ and let \( V \in L^\infty(O) \) be a real-valued function. Then, for each \( E < -\|V\|_\infty \), we have the following operator identity:
\begin{align}
\label{dif-nu-dr-lp}
&\left( H^N_{A,O} + V - E \right)^{-m} - \left( H^D_{A,O} + V - E \right)^{-m} \nonumber\\
&\quad = \sum_{j=0}^{m-1} \Bigg[
 \left( H^N_{A,O} + V - E \right)^{j - m} \chi_S^*
\left( H^D_{A,O} + V - E \right) \phi \left( H^D_{A,O} + V - E \right)^{-j - 1} \nonumber\\
&\qquad\qquad\quad
-  \left( H^N_{A,O} + V - E \right)^{j - m + 1} \chi_S^*
\phi \left( H^D_{A,O} + V - E \right)^{-j - 1}
\Bigg].
\end{align}
\end{prop}
\begin{proof}
For \( f, g \in L^2(O) \), we compute
\begin{align*}
& \left\langle f,  \left( H^N_{A,O} + V - E \right)^{-m} g \right\rangle
- \left\langle f, \left( H^D_{A,O} + V - E \right)^{-m} g \right\rangle \\
&\quad\qquad= \left\langle \left( H^N_{A,O} + V - E \right)^{-m} f, g \right\rangle
- \left\langle f, \left( H^D_{A,O} + V - E \right)^{-m} g \right\rangle \\
&\quad\qquad= \sum_{j=0}^{m-1} \Big[
\left\langle \left( H^N_{A,O} + V - E \right)^{j - m} f, \left( H^D_{A,O} + V - E \right)^{-j} g \right\rangle \\
&\qquad\qquad \qquad- \left\langle \left( H^N_{A,O} + V - E \right)^{j - m + 1} f, 
\left( H^D_{A,O} + V - E \right)^{-j - 1} g \right\rangle
\Big].
\end{align*}
Set \( f_j := \left( H^N_{A,O} + V - E \right)^{j - m} f \) and \( g_j := \left( H^D_{A,O} + V - E \right)^{-j - 1} g \). Then, by applying Proposition~\ref{dr-nu-agr} to handle the terms involving \( (1 - \phi) \), and following the same steps as in the proof of Proposition~\ref{idnty-ne-dir-rst}, we obtain the identity~\eqref{dif-nu-dr-lp}. Therefore, we omit the detailed calculations.
\end{proof}
\begin{rem}
\label{ind-rf}
Let $O_k \subset \mathbb{R}^d$, $k=1,2$.  
Since the function $u$ has compact support, there exists $R>0$ such that 
the random fields $V^\omega_{O_1}$ and $V^\omega_{O_2}$ are independent 
whenever $\operatorname{dist}(O_1, O_2) > R$,  
where $V^\omega_{O_k}(x) := \chi_{O_k}(x)\, V^\omega(x)$.
\end{rem}
\noindent We decompose $\Lambda_L = \left(-\tfrac{L}{2}, \tfrac{L}{2}\right)^d$ into disjoint smaller open sets (box annuli).
Set  $M_L = \frac{1}{2} \big[ L^\epsilon \big], ~r_{_L} = \big[ L^\delta \big]$,  where
$\epsilon > \delta > 0,   ~\epsilon + \delta = 1$.
For $k=1,2,\ldots, r_{_L}-1$,  define 
\begin{align}
\begin{split}
\label{boxes}
\Lambda^B_{L,k} &:= \left\{ x \in \mathbb{R}^d : (k-1)M_L + k(3R) < \|x\|_\infty < kM_L + (k-1)3R \right\},\\
\Lambda^S_{L,k} &:= \left\{ x \in \mathbb{R}^d : kM_L + (k-1)3R < \|x\|_\infty < kM_L + (k+1)3R \right\},\\
\Lambda_L^{r_{_L}} &:= \left\{ x \in \mathbb{R}^d : \frac{L}{2} - \left( M_L - r_{_L} 3R \right) < \|x\|_\infty < \frac{L}{2} \right\},\\
\Lambda_{3R} &:= \left\{ x \in \mathbb{R}^d : \|x\|_\infty < 3R \right\},~\text{where}~\|x\|_\infty = \max_{1 \le i \le d} |x_i|.
\end{split}
\end{align}
It is clear that 
\begin{equation}
\label{int-cl}
\left(-\frac{L}{2}, \frac{L}{2}\right)^d  
= \left( \overline{ \Lambda_{3R}  
\ \sqcup\ \bigsqcup_{k=1}^{r_{_L}-1} \Lambda^B_{L,k}  
\ \sqcup\ \bigsqcup_{k=1}^{r_{_L}-1} \Lambda^S_{L,k}  
\ \sqcup\ \Lambda_L^{r_{_L}} } \right)^\circ .
\end{equation}
We can also obtain the asymptotic volume estimates of all the above boxes, given by
\begin{align} \begin{split} 
\label{vl-est-box} \big|\Lambda_{L,k}^B \big|=\mathcal{O}\big( (2M_L)^dk^{d-1} \big)\qquad&\&\qquad \big|\Lambda_{L,k}^S \big|=\mathcal{O}\big( (2M_L)^{d-1}k^{d-1} \big)\\
 \big|\Lambda_{L}^{r_L} \big|=\mathcal{O}\big( (2M_L)L^{d-1} \big)\qquad&\&\qquad \big|\Lambda_{3R} \big|=(6R)^d. \end{split} \end{align}
Let $\ell_L=M_L^\gamma$,  $0<\gamma<1$ and we define the $\ell_L$-interior of the above boxes as
\begin{align}
\label{boxes-int}
\big(\Lambda^B_{L,k}\big)^{\circ}_{\ell_L} &:= \left\{ x \in \mathbb{R}^d : (k-1)M_L + k(3R)+\ell_L < \|x\|_\infty < kM_L + (k-1)3R-\ell_L \right\},\nonumber\\
\big(\Lambda^S_{L,k}\big)^\circ_{\ell_L} &:= \left\{ x \in \mathbb{R}^d : kM_L + (k-1)3R+\ell_L < \|x\|_\infty < kM_L + (k+1)3R-\ell_L \right\},\nonumber\\
\big(\Lambda_L^{r_{_L}} \big)^{\circ}_{\ell_L}&:= \left\{ x \in \mathbb{R}^d : \frac{L}{2} - \left( M_L - r_{_L} 3R \right)+\ell_L < \|x\|_\infty < \frac{L}{2}-\ell_L \right\}.
\end{align}
For \( m > d+1 \), we define the random variables associated with the above boxes as follows:  
\begin{align}
\label{rv-box}
\begin{split}
Y_{\Lambda^B_{L,k}}(\omega) &:= \operatorname{Tr} \big( H^\omega_{\Lambda^B_{L,k}, D} - E \big)^{-m} - \mathbb{E}\bigg[ \operatorname{Tr} \big( H^\omega_{\Lambda^B_{L,k}, D} - E \big)^{-m} \bigg], \\
Y_{\Lambda^S_{L,k}}(\omega) &:= \operatorname{Tr} \big( H^\omega_{\Lambda^S_{L,k}, D} - E \big)^{-m} - \mathbb{E}\bigg[ \operatorname{Tr} \big( H^\omega_{\Lambda^S_{L,k}, D} - E \big)^{-m} \bigg], \\
Y_{\Lambda^{r_{_L}}_L}(\omega) &:= \operatorname{Tr} \big( H^\omega_{\Lambda^{r_{_L}}_{L}, D} - E \big)^{-m} - \mathbb{E}\bigg[ \operatorname{Tr} \big( H^\omega_{\Lambda^{r_{_L}}_{L}, D} - E \big)^{-m} \bigg], \\
Y_{\Lambda_{3R}}(\omega) &:= \operatorname{Tr} \big( H^\omega_{\Lambda_{3R}, D} - E \big)^{-m} - \mathbb{E}\bigg[ \operatorname{Tr} \big( H^\omega_{\Lambda_{3R}, D} - E \big)^{-m} \bigg], \\
Y_{D,L}(\omega) &:=Y_{\Lambda_L}(\omega)= \operatorname{Tr} \big( H^\omega_{\Lambda_L, D} - E \big)^{-m} - \mathbb{E}\bigg[ \operatorname{Tr} \big( H^\omega_{\Lambda_L, D} - E \big)^{-m} \bigg].
\end{split}
\end{align}
Here, \( H^\omega_{O,D} \) denotes the finite-volume Dirichlet restriction of \( H^\omega \) to the open set \( O \), as defined in~\eqref{NDrst}.
\begin{rem}
\label{ind-resl}
In view of Remark~\ref{ind-rf}, from the construction~\eqref{boxes} of the box $\Lambda_{L,k}^B$, it is clear that $\big\{Y_{\Lambda_{L,k}^B}\big\}_{k=1}^{r_{_L}-1}$ forms an independent collection of random variables, and the same holds for $\big\{Y_{\Lambda_{L,k}^S}\big\}_{k=1}^{r_{_L}-1}$.
\end{rem}
\begin{rem}
\label{smth-fun}
We can always construct a function $\varphi_{L,k}^B \in C^\infty(\Lambda_{L,k}^B)$
such that $\operatorname{supp}(\varphi_{L,k}^B) 
\subseteq \Lambda_{L,k}^B \setminus \big(\Lambda_{L,k}^B\big)^{\circ}_{\tfrac{\ell_L}{2}}, ~0 \leq \varphi_{L,k}^B \leq 1$,
and $\operatorname{supp}(1 - \varphi_{L,k}^B) \subset\subset \Lambda_{L,k}^B $.
Since $\ell_L \to \infty$ as $L \to \infty$, it is also possible to construct 
$\varphi_{L,k}^B$ in such a way that $\displaystyle\sup_{L,k} \|\nabla \varphi_{L,k}^B\|_\infty < \infty~ \text{and}~\displaystyle\sup_{L,k} \|\Delta \varphi_{L,k}^B\|_\infty < \infty $.\\
The function $\varphi_{L,k}^B$ is equal to $1$ in a neighbourhood of the inner boundary of $\Lambda_{L,k}^B$, and equal to $0$ in the $\frac{\ell_L}{2}$-interior
$\big(\Lambda_{L,k}^B\big)^{\circ}_{\tfrac{\ell_L}{2}}$.  Between these two regions, the function is smoothly interpolated, taking values between $0$ and $1$. As $L$ and $k$ vary, the volume of the box $\Lambda_{L,k}^B$ changes, and accordingly we adjust the regions where the function is identically $1$ or $0$. However, the transition layer connecting these two regions is kept of the same structure, independent of $L$ and $k$. Moreover, we can ensure that the Lebesgue measure of the transition region remains bounded below by a fixed positive constant as $L$ and $k$ vary.\\~\\
Similarly, with respect to the boxes $\Lambda_{L,k}^S$ and $\Lambda_{L}^{r_{_L}}$, we can construct functions $\varphi_{L,k}^S \in C^\infty(\Lambda_{L,k}^S)$ and $\varphi_{L}^{r_{_L}} \in C^\infty(\Lambda_{L}^{r_{_L}})$ that have the same properties as $\varphi_{L,k}^B$.
\end{rem}
\noindent The following result shows that if we consider a region \(F\) well inside the interior of the boxes, then the local traces of the resolvent on the whole domain \(\Lambda_L\) and on the smaller subdomains are nearly the same; their difference decays exponentially fast with the distance of \(F\) from the boundary of \(\Lambda_L\).
\begin{prop}\label{ex-dy-tr}
Consider the open sets $\Lambda_L$, $\Lambda_{L,k}^B$, $\Lambda_{L,k}^S$, and $\Lambda_L^{r_{_L}}$ together with their $\ell_L$-interiors 
$\big(\Lambda_{L,k}^B\big)^{\circ}_{\ell_L}$, $\big(\Lambda_{L,k}^S\big)^{\circ}_{\ell_L}$, and $\big(\Lambda_{L}^{r_{_L}}\big)^{\circ}_{\ell_L}$, as defined in \eqref{boxes} and \eqref{boxes-int}.  
Assume that $V\in L^\infty(\mathbb{R}^d)$ is real-valued, $m>d+1$, and $E<-\|V\|_\infty$.  Then, for any $F\subset \big(\Lambda^B_{L,k}\big)^{\circ}_{\ell_L}$ or
$F\subset \big(\Lambda^S_{L,k}\big)^{\circ}_{\ell_L}$ or $F\subset \big(\Lambda^{r_{_L}}_{L}\big)^{\circ}_{\ell_L}$,
the following exponential decay of the trace difference holds:
\begin{equation}
\label{tr-df-dec}
\bigg| \operatorname{Tr}\!\left(\chi_{F}\,(H^D_{A,\Lambda_{L}}+V-E)^{-m}\chi_{F}^*\right)
- \operatorname{Tr}\!\left(\chi_{F}\,(H^D_{A,\Lambda}+V-E)^{-m}\chi_{F}^*\right) \bigg|
\;\leq\; C e^{-\beta \ell_L/2}.
\end{equation}
Here $\Lambda \in \{\Lambda_{L,k}^B,\, \Lambda_{L,k}^S,\, \Lambda^{r_{_L}}_L\}$, and the constants $C,\beta>0$ are independent of $V$ and of the boxes $\Lambda_{L,k}^B$, $\Lambda_{L,k}^S$, and $\Lambda^{r_{_L}}_L$. The operators $\chi_F$ and $\chi_F^*$ are defined as in Definition \ref{cn-ext-rst}  
\end{prop}

\begin{proof}
We treat the case $\Lambda=\Lambda_{L,k}^B$; the other cases are entirely analogous.
Since $F\subset \big(\Lambda^B_{L,k}\big)^{\circ}_{\ell_L}\subset\Lambda^B_{L,k}$, in the first trace term on the left-hand side of \eqref{sam-tr-df} (below), we consider $\chi_F : L^2(\Lambda_L) \longrightarrow L^2(F),$ and $\chi_F^* : L^2(F) \longrightarrow L^2(\Lambda_L)$ denotes its adjoint, and in the second trace term, $\chi_F : L^2(\Lambda_{L,k}^B) \longrightarrow L^2(F)$,
and $\chi_F^* : L^2(F) \longrightarrow L^2(\Lambda_{L,k}^B)$ is again the corresponding adjoint. Then we may write
\begin{align}
\label{sam-tr-df}
&\bigg| \operatorname{Tr}\!\left(\chi_{F}\,(H^D_{A,\Lambda_{L}}+V-E)^{-m}\chi_{F}^*\right)
- \operatorname{Tr}\!\left(\chi_{F}\,(H^D_{A,\Lambda_{L,k}^B}+V-E)^{-m}\chi_{F}^*\right) \bigg|\nonumber\\
&~=\bigg| \operatorname{Tr}\!\left(\chi_{F}\chi_{\Lambda_{L,k}^B}\,(H^D_{A,\Lambda_{L}}+V-E)^{-m}\chi_{\Lambda_{L,k}^B}^*\chi_F^*\right)- \operatorname{Tr}\!\left(\chi_{F}\,(H^D_{A,\Lambda_{L,k}^B}+V-E)^{-m}\chi_F^*\right) \bigg|.
\end{align}
Applying Proposition \ref{idnty-ne-dir-rst} with $\tilde{O}=\Lambda_L$, $O=\Lambda_{L,k}^B$, $\phi=\varphi_{L,k}^B$, and $S=\operatorname{supp}(\varphi_{L,k}^B)\subseteq \Lambda_{L,k}^B \setminus \big(\Lambda_{L,k}^B\big)^{\circ}_{\tfrac{\ell_L}{2}}$, where $\varphi_{L,k}^B$ is as in Remark \ref{smth-fun}, and noting that $\chi_F\chi_{\Lambda_{L,k}^B}=\chi_F$, we obtain
\begin{align}
\label{rst-dif-lk}
&\chi_F\chi_{\Lambda_{L,k}^B} \left( H^D_{A,\Lambda_L} + V - E \right)^{-m} \chi_{\Lambda_{L,k}^B}^*\chi_{F}^* - \chi_F\left( H^D_{A,\Lambda_{L,k}^B} + V - E \right)^{-m} \chi_{F}^*\nonumber\\
& =\sum_{j=0}^{m-1} \Bigg[
\chi_{F} \left( H^D_{A,\Lambda_L} + V - E \right)^{j - m} \chi_S^*
\left( H^D_{A,\Lambda_{L,k}^B} + V - E \right) \varphi_{L,k}^B \left( H^D_{A,\Lambda_{L,k}^B} + V - E \right)^{-j - 1}\chi_{F}^* \nonumber\\
& \quad - \chi_{F} \left( H^D_{A,\Lambda_L} + V - E \right)^{j - m + 1} \chi_S^*
\varphi_{L,k}^B \left( H^D_{A,\Lambda_{L,k}^B} + V - E \right)^{-j - 1}\chi_{F}^*
\Bigg].
\end{align}
Since $\max\{m-j-1, j\}\ge \tfrac{m-1}{2}>\tfrac{d}{2}$, either $\left( H^D_{A,\Lambda_L} + V - E \right)^{j - m + 1}$ or $\left( H^D_{A,\Lambda_{L,k}^B} + V - E \right)^{-j }$ belongs to the trace class.  
From the construction of $\varphi_{L,k}^B$ in Remark \ref{smth-fun} and Corollary \ref{inv-mupt-op}, it follows that $\left( H^D_{A,\Lambda_{L,k}^B} + V - E \right) \varphi_{L,k}^B \left( H^D_{A,\Lambda_{L,k}^B} + V - E \right)^{-1}$
is a bounded operator, with operator norm uniformly bounded independently of $\Lambda_{L,k}^B$.   Also we have $\| \chi_F^*\|\leq 1$.\\
Case 1. Suppose $\left( H^D_{A,\Lambda_L} + V - E \right)^{j - m + 1}$ is trace class. Then
\begin{align}
\label{rst-dif-lk-1}
&\bigg\|\chi_F\chi_{\Lambda_{L,k}^B} \left( H^D_{A,\Lambda_L} + V - E \right)^{-m} \chi_{\Lambda_{L,k}^B}^*\chi_{F}^*- \chi_F\left( H^D_{A,\Lambda_{L,k}^B} + V - E \right)^{-m}\chi_{F}^*\bigg\|_1 \nonumber\\
& \leq\sum_{j=0}^{m-1}\Bigg[ \Bigg\|
\chi_{F} \left( H^D_{A,\Lambda_L} + V - E \right)^{j - m} \chi_S^*\Bigg\|_1\nonumber\\
&\qquad\times
\Bigg\|\left( H^D_{A,\Lambda_{L,k}^B} + V - E \right) \varphi_{L,k}^B \left( H^D_{A,\Lambda_{L,k}^B} + V - E \right)^{-j - 1}\Bigg\|_\infty \nonumber\\
& \quad +
\Bigg\|
\chi_{F} \left( H^D_{A,\Lambda_L} + V - E \right)^{j - m+1} \chi_S^*\Bigg\|_1
\Bigg\|\varphi_{L,k}^B \left( H^D_{A,\Lambda_{L,k}^B} + V - E \right)^{-j - 1}\Bigg\|_\infty\Bigg]\nonumber\\
&\leq m C \Bigg(
\Bigg\|
\chi_{F} \left( H^D_{A,\Lambda_L} + V - E \right)^{j - m} \chi_S^*\Bigg\|_1+
\Bigg\|\chi_{F} \left( H^D_{A,\Lambda_L} + V - E \right)^{j - m+1} \chi_S^*\Bigg\|_1\Bigg)\nonumber\\
& \leq 2m C_1 e^{-\beta\ell_L/2}.
\end{align}
In the last inequality, we use that $\operatorname{dist}(F,S)>\tfrac{\ell_L}{2}$ together with the Combes--Thomas estimate \eqref{exp-decy}.\\
Case 2.  Suppose instead that $\left( H^D_{A,\Lambda_{L,k}^B} + V - E \right)^{-j }$ is trace class. Then
\begin{align}
\label{rst-dif-lk-2}
&\bigg\|\chi_F\chi_{\Lambda_{L,k}^B} \left( H^D_{A,\Lambda_L} + V - E \right)^{-m} \chi_{\Lambda_{L,k}^B}^*- \chi_F\left( H^D_{A,\Lambda_{L,k}^B} + V - E \right)^{-m}\bigg\|_1 \nonumber\\
& \leq\sum_{j=0}^{m-1}\Bigg[ \Bigg\|
\chi_{F} \left( H^D_{A,\Lambda_L} + V - E \right)^{j - m} \chi_S^*\Bigg\|_\infty\nonumber\\
&\qquad\times
\Bigg\|\left( H^D_{A,\Lambda_{L,k}^B} + V - E \right) \varphi_{L,k}^B \left( H^D_{A,\Lambda_{L,k}^B} + V - E \right)^{-1} \left( H^D_{A,\Lambda_{L,k}^B} + V - E \right)^{-j }\Bigg\|_1\nonumber\\
& \quad+
\Bigg\|
\chi_{F} \left( H^D_{A,\Lambda_L} + V - E \right)^{j - m+1} \chi_S^*\Bigg\|_\infty
\Bigg\|\varphi_{L,k}^B \left( H^D_{A,\Lambda_{L,k}^B} + V - E \right)^{-j - 1}\Bigg\|_1\Bigg]\nonumber\\
&\leq m C \Bigg(
\Bigg\|
\chi_{F} \left( H^D_{A,\Lambda_L} + V - E \right)^{j - m} \chi_S^*\Bigg\|_\infty+
\Bigg\|\chi_{F} \left( H^D_{A,\Lambda_L} + V - E \right)^{j - m+1} \chi_S^*\Bigg\|_\infty\Bigg)\nonumber\\
& \leq 2m C_2 e^{-\beta\ell_L/2}.
\end{align}
Here again we used the Combes--Thomas estimate \eqref{exp-decy}, together with the fact that $\operatorname{dist}(F,S)>\tfrac{\ell_L}{2}$.
Combining \eqref{sam-tr-df}, \eqref{rst-dif-lk}, \eqref{rst-dif-lk-1}, and \eqref{rst-dif-lk-2}, the claimed estimate (\ref{tr-df-dec}) follows.
\end{proof}
\noindent Define the $\tilde{\ell}_L$-interior of $\Lambda_L=\Big(-\tfrac{L}{2}, \tfrac{L}{2} \Big)^d$ as
\begin{equation}
\label{int-bl}
\big(\Lambda_L\big)_{\tilde{\ell}_L}^{\circ}=\Big(-\tfrac{L}{2}+ \tilde{\ell}_L,\, \tfrac{L}{2}-\tilde{\ell}_L \Big)^d,  ~~\tilde{\ell}_L=L^\alpha, ~~ 0<\alpha<1.
\end{equation}
\begin{rem}
\label{fn-bl}
We can construct a function $\varphi_L \in C^\infty(\Lambda_L)$ such that  
$\operatorname{supp}(\varphi_L)\subseteq \Lambda_L \setminus \big(\Lambda_L\big)_{\frac{\tilde{\ell}_L}{2}}^{\circ}, ~0\leq \varphi_L \leq 1$,
and $\operatorname{supp}(1-\varphi_L)\subset\subset \Lambda_L$.  
Moreover, we may choose $\varphi_L$ so that  
$\sup_{L}\|\nabla \varphi_L\|_\infty < \infty, ~ 
\sup_{L}\|\Delta \varphi_L\|_\infty < \infty$. 
The construction is analogous to that one described in Remark~\ref{smth-fun}.
\end{rem}
\noindent By applying the same method with $\varphi_L=\varphi_{L,k}^B$, as in Proposition~\ref{ex-dy-tr}, we obtain the following result in the case where the smaller cube is $\Lambda_L$ and the larger domain is the entire space $\mathbb{R}^d$.
\begin{prop}
\label{tr-dec-bg}
Let \( V \in L^\infty(\mathbb{R}^d) \) be real-valued, \( E < -\|V\|_\infty \), and \( m > d + 1 \). 
Then, for any set \( F \subset (\Lambda_L)_{\tilde{\ell}_L}^{\circ} \), we have
\begin{equation*}
\bigg| \operatorname{Tr}\!\left(\chi_{F}\,(H_{A}+V-E)^{-m}\chi_{F}^*\right)
- \operatorname{Tr}\!\left(\chi_{F}\,(H^D_{A,\Lambda_L}+V-E)^{-m}\chi_{F}^*\right) \bigg|
\;\leq\; C e^{-\beta \tilde{\ell}_L/2}.
\end{equation*}
Here $H_A: = H^D_{A,\mathbb{R}^d}$.  The constants $C,~\beta>0$ are independent of $\Lambda_L$ and $V$.
\end{prop}
\begin{proof}
The proof follows the same method as in Proposition~\ref{ex-dy-tr}, making use of the function $\varphi_L$ described in Remark~\ref{fn-bl}. We therefore omit the details.
\end{proof}
\noindent The following result provides an exponential estimate for the difference 
between the Neumann and Dirichlet resolvent traces associated with the box \( \Lambda_L \), 
when restricted to a set \( F \subset \Lambda_L \) that is located sufficiently far from the boundary of \( \Lambda_L \).
\begin{prop}
\label{df-n-d-in}
Let \( V \in L^\infty(\mathbb{R}^d) \) be real-valued, \( E < -\|V\|_\infty \), and \( m > d + 1 \). 
Then, for any set \( F \subset (\Lambda_L)_{\tilde{\ell}_L}^{\circ} \), we have
\begin{equation}
\label{tr-df-nd}
\bigg| 
\operatorname{Tr}\!\left(\chi_{F}\,(H^N_{A,\Lambda_{L}}+V-E)^{-m}\chi_{F}^*\right)
- \operatorname{Tr}\!\left(\chi_{F}\,(H^D_{A,\Lambda_{L}}+V-E)^{-m}\chi_{F}^*\right) 
\bigg|
\;\leq\; C\, e^{-\beta \tilde{\ell}_L/2},
\end{equation}
where the constants \( C, \beta > 0 \) are independent of \( V \) and \( \Lambda_L \).
\begin{proof}
The proof proceeds similarly to Proposition~\ref{ex-dy-tr}, by applying the identity~(\ref{dif-nu-dr-lp}) with \( \phi = \varphi_L \), 
as in Remark~\ref{fn-bl}, together with the Combes--Thomas estimate~\eqref{exp-decy} 
for the Neumann and Dirichlet restrictions.
\end{proof}
\end{prop}

\section{CLT for Laurent Polynomials}
In this section, we derive the CLT~\eqref{clt-cn} for test functions that are Laurent polynomials, that is,
$f(x) = (x - E)^{-m} \sum_{k=0}^p a_k (x - E)^{-k}$,
defined on the interval $\big[-\|V\|_\infty, \infty\big)$, where $E < -\|V\|_\infty$, 
$m > d+1,~p \in \mathbb{N} \cup \{0\}$, and $a_k \in \mathbb{R}$.\\~\\
Before stating the next result, we note that, under suitable conditions,  global quantities can be approximated by sums of their localized counterparts.  
The following proposition makes this precise by showing that the variance of the discrepancy 
between \(Y_{D,L}\) and the sum of its localized contributions vanishes in the large-volume limit.
\begin{prop}
\label{e-q-cn}
Let the random variables $Y_{\Lambda}$, 
$\Lambda \in \{\Lambda_L, \Lambda^B_{L,k}, \Lambda^S_{L,k}, 
\Lambda_L^{r_{_L}}, \Lambda_{3R}\}$, be defined as in~\eqref{rv-box}.
Assume that the single site distribution (SSD) $\mu$ has compact support, and that $u \in L^\infty(\mathbb{R}^d)$ is a real-valued, compactly supported function.
Furthermore, let $A \in \big(L^2_{\mathrm{loc}}(\mathbb{R}^d)\big)^d$.
Then, the following convergence holds:
\begin{equation}
\label{eqv-cn}
\frac{1}{|\Lambda_L|}\mathbb{E}\bigg[
Y_{D,L}-\sum_{k=1}^{r_{_L}-1} Y_{\Lambda^B_{L,k}}
-\sum_{k=1}^{r_{_L}-1} Y_{\Lambda^S_{L,k}}
-Y_{\Lambda_L^{r_{_L}}}
-Y_{\Lambda_{3R}}
\bigg]^2 \xrightarrow[L \to \infty]{} 0 .
\end{equation}
\end{prop}
\begin{proof}
First, we define the index set $B_L$ by
\begin{equation}
\label{indx-set}
B_L := \Big\{ n \in \mathbb{Z}^d : \operatorname{supp}(u_n) \cap \Lambda_L \neq \emptyset \Big\} \subset \mathbb{Z}^d, ~\text{where}~
u_n(x) := u(x-n).
\end{equation}
Then the random variable
\begin{equation}
\label{dfn-G}
G_L(\omega) = Y_{D,L}(\omega) - \sum_{k=1}^{r_{_L}-1} Y_{\Lambda^B_{L,k}}(\omega) - \sum_{k=1}^{r_{_L}-1} Y_{\Lambda^S_{L,k}}(\omega) - Y_{\Lambda_L^{r_{_L}}}(\omega) - Y_{\Lambda_{3R}}(\omega)
\end{equation}
depends only on the collection $\{\omega_n\}_{n \in B_L}$.  
Since $u$ has compact support, we also have  
$\displaystyle\lim_{L \to \infty} \frac{\#B_L}{|\Lambda_L|} = 1$. 
Now define a subset $B_{L,\ell_L} \subset B_L$ by
\begin{align}
\label{sp-int}
B_{L,\ell_L} &= \Big\{ n \in B_L : \operatorname{supp}(u_n) \subset (\Lambda^B_{L,k})^{\circ}_{\ell_L}~\text{or}~\operatorname{supp}(u_n) \subset (\Lambda^S_{L,k})^{\circ}_{\ell_L}\nonumber\\
&\qquad \qquad\text{or}~\operatorname{supp}(u_n) \subset (\Lambda_L^{r_{_L}})^{\circ}_{\ell_L}~\text{for some}~k=1,2,\ldots,r_{_L}-1\Big\}.
\end{align}
Choose an enumeration $\{\omega_{n_j}\}_{j=1}^{\#B_L}$ of $\{\omega_n\}_{n \in B_L}$ such that the first $\#B_{L,\ell_L}$ terms,  
$\{\omega_{n_j}\}_{j=1}^{\#B_{L,\ell_L}}$, correspond exactly to $\{\omega_n\}_{n \in B_{L,\ell_L}}$.  That is,  $\{\omega_{n_j}\}_{j=1}^{\#B_{L,\ell_L}} = \{\omega_n\}_{n \in B_{L,\ell_L}}$.
Now define a filtration $\{\mathcal{F}_s\}_{s=0}^{\#B_L}$ of $\sigma$-algebras by  
$\mathcal{F}_s = \sigma(\omega_{n_j} : 1 \leq j \leq s)$,  $\mathcal{F}_0 = \{\emptyset, \Omega\}$.  Since $\mathbb{E}[G_L(\omega)]=0$, the variance of $G_L$ can be expressed as the sum of squares of a martingale difference sequence:
\begin{align}
\label{est-th-mrt}
\begin{split}
&\mathbb{E}\big[G_L^2\big]
= \sum_{s=1}^{\#B_L} \mathbb{E}\Big[\mathbb{E}(G_L \mid \mathcal{F}_s) - \mathbb{E}(G_L \mid \mathcal{F}_{s-1}) \Big]^2 \\
&= \sum_{s=1}^{\#B_{L,\ell_L}} \mathbb{E}\Big[ \mathbb{E}(G_L \mid \mathcal{F}_s) - \mathbb{E}(G_L \mid \mathcal{F}_{s-1}) \Big]^2+ \sum_{s=1+\#B_{L,\ell_L}}^{\#B_L} \mathbb{E}\Big[ \mathbb{E}(G_L \mid \mathcal{F}_s) - \mathbb{E}(G_L \mid \mathcal{F}_{s-1}) \Big]^2.
\end{split}
\end{align}
We now estimate each martingale difference term.  
For $s = 1,2,\dots,\#B_{L,\ell_L}$, we write
\begin{align}
\label{tr-dif-drv}
&\mathbb{E}(G_L \mid \mathcal{F}_s) - \mathbb{E}(G_L \mid \mathcal{F}_{s-1})= \mathbb{E}(G_L(\omega) \mid \mathcal{F}_s) - \mathbb{E}(G_L(\omega : \omega_{n_s}=0) \mid \mathcal{F}_s) \nonumber\\
&\quad \qquad\qquad\qquad+ \mathbb{E}(G_L(\omega : \omega_{n_s}=0) \mid \mathcal{F}_{s-1}) - \mathbb{E}(G_L(\omega) \mid \mathcal{F}_{s-1}) \nonumber\\
&= \mathbb{E}\Bigg(\int_0^1 \frac{d}{dt} G_L(\omega : \omega_{n_s} \to t \omega_{n_s}) \,dt\Big|\, \mathcal{F}_s \Bigg) - \mathbb{E}\Bigg( \int_0^1\frac{d}{dt} G_L(\omega : \omega_{n_s} \to t \omega_{n_s}) \,dt\Big|\, \mathcal{F}_{s-1} \Bigg).
\end{align}
For such $s$, at most one of the random variables in the collection
$\{ Y_{\Lambda^B_{L,k}},\, Y_{\Lambda^S_{L,k}},\, Y_{\Lambda^{r_{_L}}_L} : k=1,2,\dots,r_{_L}-1 \}$
depends on $\omega_{n_s}$.  
Let $n_s \in \Lambda$, with $\Lambda \in \{\Lambda^B_{L,k}, \Lambda^S_{L,k}, \Lambda^{r_{_L}}_L\}$.  Using the definition \eqref{dfn-G} of $G_L(\omega)$, we obtain
\begin{align} \label{tr-ag} &\mathbb{E}\big(G_L|\mathcal{F}_s \big)-\mathbb{E}\big(G_L|\mathcal{F}_{s-1} \big)\\ &=\mathbb{E}\bigg(\int_0^1\frac{d}{dt}Y_{D,L}\big(\omega:\omega_{n_s}\to t\omega_{n_s} \big)dt\big|\mathcal{F}_s \bigg) -\mathbb{E}\bigg(\int_0^1\frac{d}{dt}Y_{\Lambda}\big(\omega:\omega_{n_s}\to t\omega_{n_s} \big)dt\big|\mathcal{F}_s \bigg)\nonumber\\ 
&\quad-\mathbb{E}\bigg(\int_0^1\frac{d}{dt}Y_{D,L}\big(\omega:\omega_{n_s}\to t\omega_{n_s} \big)dt\big|\mathcal{F}_{s-1} \bigg) +\mathbb{E}\bigg(\int_0^1\frac{d}{dt}Y_{\Lambda}\big(\omega:\omega_{n_s}\to t\omega_{n_s} \big)dt\big|\mathcal{F}_{s-1} \bigg)\nonumber\\ &=-m \, \mathbb{E}\left( \int_0^1 \omega_{n_s} \operatorname{Tr} \left( \chi_{\Lambda_L} u_{n_s}(x) \left( H^D_{A,\Lambda_L} + V^\omega_{\Lambda_L} - E \right)^{-m-1}_{(\omega_{n_s} \to t \omega_{n_s})} \right) dt \, \bigg| \, \mathcal{F}_s \right)\nonumber \\ & \qquad+m \, \mathbb{E}\left( \int_0^1 \omega_{n_s} \operatorname{Tr} \left( \chi_{\Lambda} u_{n_s}(x) \left( H^D_{A,\Lambda} + V^\omega_{\Lambda} - E \right)^{-m-1}_{(\omega_{n_s} \to t \omega_{n_s})} \right) dt \, \bigg| \, \mathcal{F}_s \right)\nonumber\\ &\qquad +m \, \mathbb{E}\left( \int_0^1 \omega_{n_s} \operatorname{Tr} \left( \chi_{\Lambda_L} u_{n_s}(x) \left( H^D_{A,\Lambda_L} + V^\omega_{\Lambda_L} - E \right)^{-m-1}_{(\omega_{n_s} \to t \omega_{n_s})} \right) dt \, \bigg| \, \mathcal{F}_{s-1} \right)\nonumber \\ & \qquad-m \, \mathbb{E}\left( \int_0^1 \omega_{n_s} \operatorname{Tr} \left( \chi_{\Lambda} u_{n_s}(x) \left( H^D_{A,\Lambda} + V^\omega_{\Lambda} - E \right)^{-m-1}_{(\omega_{n_s} \to t \omega_{n_s})} \right) dt \, \bigg| \, \mathcal{F}_{s-1} \right).\nonumber 
\end{align}
Now consider the case $\Lambda = \Lambda_{L,k}^B$; the other cases are analogous.  
Since $F = \operatorname{supp}(u_{n_s}) \subset (\Lambda_{L,k}^B)^{\circ}_{\ell_L} \subset \Lambda_{L,k}^B \subset \Lambda_L$,
it follows from Corollary~\ref{cr-tr-rst} and the estimate~\eqref{tr-df-dec} that
\begin{align} 
\label{tr-ag-ag}
&\bigg| \operatorname{Tr} \left( \chi_{\Lambda_L} u_{n_s}(x) \left( H^D_{A,\Lambda_L} + V^\omega_{\Lambda_L} - E \right)^{-m-1}_{(\omega_{n_s} \to t \omega_{n_s})} \right)\nonumber\\ & \qquad\qquad \qquad-\operatorname{Tr} \left( \chi_{\Lambda_{L,k}^B} u_{n_s}(x) \left( H^D_{A,\Lambda_{L,k}^B} + V^\omega_{\Lambda_{L,k}^B} - E \right)^{-m-1}_{(\omega_{n_s} \to t \omega_{n_s})} \right) \bigg|\nonumber\\ &\qquad\qquad = \bigg| \operatorname{Tr} \left( u_{n_s}\chi_{F} \left( H^D_{A,\Lambda_L} + V^\omega_{\Lambda_L} - E \right)^{-m-1}_{(\omega_{n_s} \to t \omega_{n_s})}\chi_F^* \right)\nonumber\\ & \qquad\qquad \qquad\qquad-\operatorname{Tr} \left(u_{n_s} \chi_{F} \left( H^D_{A,\Lambda_{L,k}^B} + V^\omega_{\Lambda_{L,k}^B} - E \right)^{-m-1}_{(\omega_{n_s} \to t \omega_{n_s})} \chi_F^* \right) \bigg|\nonumber\\ &\qquad\qquad \leq \|u\|_\infty\bigg| \operatorname{Tr} \left( \chi_{F} \left( H^D_{A,\Lambda_L} + V^\omega_{\Lambda_L} - E \right)^{-m-1}_{(\omega_{n_s} \to t \omega_{n_s})}\chi_F^* \right)\nonumber\\ & \qquad\qquad \qquad\qquad-\operatorname{Tr} \left( \chi_{F} \left( H^D_{A,\Lambda_{L,k}^B} + V^\omega_{\Lambda_{L,k}^B} - E \right)^{-m-1}_{(\omega_{n_s} \to t \omega_{n_s})} \chi_F^* \right) \bigg|\nonumber\\ &\qquad\qquad \leq \| u\|_\infty Ce^{-\beta \ell_L/2}. 
\end{align}
Substituting \eqref{tr-ag-ag} into \eqref{tr-ag}, we obtain, for $s=1,\dots,\#B_{L,\ell_L}$,
\begin{align}
\label{tr-aggg}
\Big| \mathbb{E}(G_L \mid \mathcal{F}_s) - \mathbb{E}(G_L \mid \mathcal{F}_{s-1}) \Big|
\leq \|\omega_{n_s}\|_\infty \|u\|_\infty C e^{-\beta \ell_L/2}.
\end{align}
For $s = 1+\#B_{L,\ell_L}, 2+\#B_{L,\ell_L},\dots,\#B_L$, at most finitely many, say $M$ (independent of $s$), random variables from the collection
$\{ Y_{\Lambda^B_{L,k}},\, Y_{\Lambda^S_{L,k}},\, Y_{\Lambda^{r_{_L}}_L},\, Y_{\Lambda_{3R}} : k=1,2,\dots,r_{_L}-1 \}$ depend on $\omega_{n_s}$.  
Proceeding as in \eqref{tr-dif-drv} and using \eqref{dfn-G}, we obtain
\begin{align} \label{tr-ag-ut-bd} &\mathbb{E}\big(G_L|\mathcal{F}_s \big)-\mathbb{E}\big(G_L|\mathcal{F}_{s-1} \big)\\ &=\mathbb{E}\bigg(\int_0^1\frac{d}{dt}Y_{D,L}\big(\omega:\omega_{n_s}\to t\omega_{n_s} \big)dt\big|\mathcal{F}_s \bigg) -\sum_{j=1}^M\mathbb{E}\bigg(\int_0^1\frac{d}{dt}Y_{\Lambda_j}\big(\omega:\omega_{n_s}\to t\omega_{n_s} \big)dt\big|\mathcal{F}_s \bigg)\nonumber\\ 
&-\mathbb{E}\bigg(\int_0^1\frac{d}{dt}Y_{D,L}\big(\omega:\omega_{n_s}\to t\omega_{n_s} \big)dt\big|\mathcal{F}_{s-1} \bigg) +\sum_{j=1}^M\mathbb{E}\bigg(\int_0^1\frac{d}{dt}Y_{\Lambda_j}\big(\omega:\omega_{n_s}\to t\omega_{n_s} \big)dt\big|\mathcal{F}_{s-1} \bigg)\nonumber\\ &=-m \, \mathbb{E}\left( \int_0^1 \omega_{n_s} \operatorname{Tr} \left( \chi_{\Lambda_L} u_{n_s}(x) \left( H^D_{A,\Lambda_L} + V^\omega_{\Lambda_L} - E \right)^{-m-1}_{(\omega_{n_s} \to t \omega_{n_s})} \right) dt \, \bigg| \, \mathcal{F}_s \right) \nonumber\\ 
& +m\sum_{j=1}^M \, \mathbb{E}\left( \int_0^1 \omega_{n_s} \operatorname{Tr} \left( \chi_{\Lambda_j} u_{n_s}(x) \left( H^D_{A,\Lambda_j} + V^\omega_{\Lambda_j} - E \right)^{-m-1}_{(\omega_{n_s} \to t \omega_{n_s})} \right) dt \, \bigg| \, \mathcal{F}_s \right)\nonumber\\ 
& +m \, \mathbb{E}\left( \int_0^1 \omega_{n_s} \operatorname{Tr} \left( \chi_{\Lambda_L} u_{n_s}(x) \left( H^D_{A,\Lambda_L} + V^\omega_{\Lambda_L} - E \right)^{-m-1}_{(\omega_{n_s} \to t \omega_{n_s})} \right) dt \, \bigg| \, \mathcal{F}_{s-1} \right) \nonumber\\ 
& -m\sum_{j=1}^M \, \mathbb{E}\left( \int_0^1 \omega_{n_s} \operatorname{Tr} \left( \chi_{\Lambda_j} u_{n_s}(x) \left( H^D_{A,\Lambda_j} + V^\omega_{\Lambda_j} - E \right)^{-m-1}_{(\omega_{n_s} \to t \omega_{n_s})} \right) dt \, \bigg| \, \mathcal{F}_{s-1} \right).\nonumber \end{align}
In the above $\Lambda_j\in \big\{Y_{\Lambda_{L,k}^B}, Y_{\Lambda_{L,k}^B}, Y_{\Lambda_L^{r_{_L}}}, Y_{\Lambda_{3R}}: k=1,2,\ldots, r_{_L}-1 \big\}$ for $j=1,2,\ldots, M$.
Using the trace estimate \eqref{tr-est}, we get
\begin{align}
\label{tr-aggg-ut}
\Big| \mathbb{E}(G_L \mid \mathcal{F}_s) - \mathbb{E}(G_L \mid \mathcal{F}_{s-1}) \Big|
\leq 2m(1+M)\|\omega_{n_s}\|_\infty | \operatorname{supp}(u)|.
\end{align}
The above is true for $s=1+\#B_{L,\ell_L}, 2+\#B_{L,\ell_L},\ldots, \#B_L$.
Note that $|\cdot|$ denotes Lebesgue measure on $\mathbb{R}^d$, and since $u_{n_s}(x) = u(x-n_s)$, we have $|\operatorname{supp}(u_{n_s})| = |\operatorname{supp}(u)|$.
Finally, combining \eqref{tr-aggg} and \eqref{tr-aggg-ut} with \eqref{est-th-mrt}, we obtain
\begin{align}
\frac{1}{|\Lambda_L|}\mathbb{E}\big[G_L^2\big] \leq 
C\left( \frac{\#B_L e^{-\beta \ell_L}}{|\Lambda_L|} + \frac{\#B_L - \#B_{L,\ell_L}}{|\Lambda_L|} \right).
\end{align}
Since  
$\frac{\#B_L - \#B_{L,\ell_L}}{|\Lambda_L|} \to 0 \quad \text{as } L \to \infty$,
we obtain \eqref{eqv-cn} from the definition \eqref{dfn-G} of $G_L$.  
This completes the proof.
\end{proof}
\noindent We now prove that the limiting variance of the normalized random variable $|\Lambda_L|^{-1/2}Y_{D,L}$ exists, and we compute it explicitly.  This result will later be used to establish the existence and positivity of the limiting variance of the random variable $|\Lambda_L|^{-1/2} Y_{f,D,L}$, as defined in (\ref{rv}).
\begin{lem}
\label{exst-vr-resl}
Let $Y_{D,L}$ be the random variable defined in (\ref{rv-box}).  
Then, under Hypothesis \ref{hypo}, the limiting variance of $\displaystyle\frac{1}{|\Lambda_L|^{1/2}} Y_{D,L}$ exists and is given by
\begin{align}
\label{lmv-resol}
\lim_{L \to \infty} \frac{1}{|\Lambda_L|} \,\mathrm{Var}\!\big( Y_{D,L} \big)
&=\mathbb{E}\bigg[\mathbb{E}\!\bigg(\int_0^1 \omega_{\vec{1}_d} \,\operatorname{Tr}\!\bigg(u_{\vec{1}_d}(x)
\big(f'_E(H^\omega)\big) _{(\omega_{\vec{1}_d}\to t\omega_{\vec{1}_d})}\bigg)\,dt \,\bigg|\, \mathcal{F}^d_{\vec{1}_d} \bigg)\nonumber\\
&\quad-\mathbb{E}\!\bigg(\int_0^1 \omega_{\vec{1}_d} \,\operatorname{Tr}\!\bigg(u_{\vec{1}_d}(x)\big(f'_E(H^\omega)\big) _{(\omega_{\vec{1}_d}\to t\omega_{\vec{1}_d})}\bigg)\,dt \,\bigg|\, \mathcal{F}^d_{\vec{1}_{d-1,0}} \bigg)\bigg]^2.\nonumber
\end{align}
Here $\vec{1}_d$, $\vec{1}_{d-1,0}$, $\mathcal{F}^d_{\vec{1}_d}$, and $\mathcal{F}^d_{\vec{1}_{d-1,0}}$ are as defined in (\ref{smag}), and $f_E(x) = (x - E)^{-m}$ on $\big[-\|V\|_\infty,\infty\big)$ with $E < -\|V\|_\infty$, where $\|V^\omega\|\leq \|V\|_\infty$ a.e. $\omega$.
\end{lem}
\begin{proof}
It is clear from (\ref{indx-set}) that the random variable $Y_{D,L}(\omega)$ depends only on the 
$\{\omega_n\}_{n\in B_L}$.  Since $u_n=u(x-n)$ and $u$ has compact support, there exist a $\delta>0$ (independent of $L$), such that 
\begin{equation}
\label{lr-bx}
B_L \subseteq \Lambda_{L_\delta} 
= \Bigl\{ n = (n_1, n_2, \ldots, n_d) \in \mathbb{Z}^d : |n_i| \leq L_\delta\Bigr\},  ~
L_\delta=\left\lfloor \frac{L+\delta}{2} \right\rfloor.
\end{equation}
Define the interior set $\big(B_{L}\big)_{\tilde{\ell}_{L}}^\circ\subset \mathbb{Z}^d $ as
\begin{equation}
\label{indx-int}
\big(B_{L}\big)_{\tilde{\ell}_L}^\circ=\big\{n\in\mathbb{Z}^d: S_{n}\subset \big(\Lambda_{L}  \big)_{\tilde{\ell}_{L}}^\circ\big\}, ~~S_{n}=\operatorname{supp}(u_{n}).
\end{equation}
Here $\big(\Lambda_{L}  \big)_{\tilde{\ell}_{L}}^\circ$ is defined as in \eqref{int-bl}. To understand the existence of the limit
$\displaystyle\lim_{L \to \infty} \frac{1}{|\Lambda_L|} \,\mathrm{Var}\!\big( Y_{D,L} \big)$,
and its explicit form, we first prove it for \(L^2(\mathbb{R}^d)\) with \(d = 1,2\).
Once this is established, the case for general \(d\) follows by an induction-type argument.\\~\\
\underline{Case \(d=1\):} Denote $A^1_{k_1}=\{n\in\mathbb{Z}:n\leq k_1 \}$,  $k_1\in \mathbb{Z}$, and define the $\sigma$-algebra $\mathcal{F}^1_{k_1}=\sigma\{\omega_n : n\in A^1_{k_1} \}$. Then we have $Y_{D,L}=\mathbb{E}\!\left(Y_{D,L} \,\middle|\, \mathcal{F}^1_{L_\delta} \right)$, $\mathbb{E}[Y_{D,L}]=\mathbb{E}\!\left(Y_{D,L} \,\middle|\, \mathcal{F}^1_{-L_\delta-1} \right)$ and
$\mathcal{F}^1_{k_1}\setminus\mathcal{F}^1_{k_1-1}=\sigma\{\omega_{k_1}\}$.
Using the martingale difference technique, we can write the variance of $Y_{D,L}$ as
\begin{equation}
\label{v-d1}
\mathrm{Var}(Y_{D,L})
= \sum_{k_1=-L_\delta}^{L_\delta} 
\mathbb{E}\!\Bigg[ 
\mathbb{E}\!\big(Y_{D,L} \,\big|\, \mathcal{F}^1_{k_1} \big)
- \mathbb{E}\!\big(Y_{D,L} \,\big|\, \mathcal{F}^1_{k_1-1} \big)
\Bigg]^2=\sum_{k_1=-L_\delta}^{L_\delta} 
\mathbb{E}\!\Bigg[ Y^1_{L,k_1}\Bigg]^2.
\end{equation}
Now we estimate each martingale difference as
\begin{align}
\label{mrt-d1}
Y^1_{L,k_1}&=\mathbb{E}\!\big(Y_{D,L} \,\big|\, \mathcal{F}^1_{k_1} \big)
- \mathbb{E}\!\big(Y_{D,L} \,\big|\, \mathcal{F}^1_{k_1-1} \big)\\
&=\mathbb{E}\!\big(Y_{D,L}(\omega) \,\big|\, \mathcal{F}^1_{k_1} \big)-\mathbb{E}\!\big(Y_{D,L}(\omega: \omega_{k_1}=0)\big| \,\mathcal{F}^1_{k_1} \big)\nonumber\\
&\qquad\qquad +\mathbb{E}\!\big(Y_{D,L}(\omega:\omega_{k_1}=0) \,\big|\, \mathcal{F}^1_{k_1-1} \big)- \mathbb{E}\!\big(Y_{D,L}(\omega) \,\big|\, \mathcal{F}^1_{k_1-1} \big)\nonumber\\
&=\mathbb{E}\!\bigg(\int_0^1\frac{d}{dt}Y_{D,L}(\omega: \omega_{k_1}\to t\omega_{k_1}) dt \,\big|\, \mathcal{F}^1_{k_1} \bigg)\nonumber\\
&\qquad\qquad
-\mathbb{E}\!\bigg(\int_0^1\frac{d}{dt}Y_{D,L}(\omega: \omega_{k_1}\to t\omega_{k_1}) dt\,\big|\, \mathcal{F}^1_{k_1-1} \bigg)\nonumber\\
& =-m\mathbb{E}\!\bigg(\int_0^1\omega_{k_1} \operatorname{Tr}\bigg(u_{k_1}(x)
\big(H^\omega_{\Lambda_L, D} -E\big)^{-m-1} _{(\omega_{k_1}\to t\omega_{k_1})}\bigg) dt \,\bigg|\, \mathcal{F}^1_{k_1} \bigg)\nonumber\nonumber\\
&\qquad\qquad
+m\mathbb{E}\!\bigg(\int_0^1\omega_{k_1} \operatorname{Tr}\bigg(u_{k_1}(x)
\big(H^\omega_{\Lambda_L, D} -E\big)^{-m-1} _{(\omega_{k_1}\to t\omega_{k_1})}\bigg) dt \,\bigg|\, \mathcal{F}^1_{k_1-1} \bigg)\nonumber.
\end{align}
Let $S_{k_1}=\operatorname{supp}(u_{k_1})$,  $u_{k_1}=u(x-k_1)$. Using Corollary \ref{cr-tr-rst}, we obtain
\begin{equation}
\label{uscan}
\operatorname{Tr}\bigg(u_{k_1}(x)
\big(H^\omega_{\Lambda,D} -E\big)^{-m-1} _{(\omega_{k_1}\to t\omega_{k_1})}  \bigg)
=\operatorname{Tr}\bigg(u_{k_1}(x)\chi_{S_{k_1}}
\big(H^\omega_{\Lambda, D} -E\big)^{-m-1} _{(\omega_{k_1}\to t\omega_{k_1})}\chi_{S_{k_1}}^* \bigg)
\end{equation}
Let denote $S=\operatorname{supp}(u)$, now the Proposition \ref{est-tr-chf} gives
\begin{equation}
\label{bbd-sq}
\bigg| \omega_{k_1}\operatorname{Tr}\bigg(u_{k_1}(x)
\big(H^\omega_{\Lambda,D} -E\big)^{-m-1} _{(\omega_{k_1}\to t\omega_{k_1})}\bigg)  \bigg|
\leq C\|\omega_0\|_\infty \|u \|_\infty|S|.
\end{equation}
In the above (\ref{uscan}) and (\ref{bbd-sq}), $\Lambda=\Lambda_L$ or $\mathbb{R}^d$ and $H^\omega_{\mathbb{R}^d, D}:=H^\omega$.\\
Now for each $k_1\in \big(B_L\big)_{\tilde{\ell}_L}^\circ$ together with (\ref{uscan}), we have from the Proposition \ref{tr-dec-bg}
\begin{align}
\label{apr-fl}
&\bigg|\mathbb{E}\!\bigg(\int_0^1\omega_{k_1} \operatorname{Tr}\bigg(u_{k_1}(x)
\big(H^\omega_{\Lambda_L, D} -E\big)^{-m-1} _{(\omega_{k_1}\to t\omega_{k_1})}\bigg) dt \,\bigg|\, \mathcal{F}^1_{k_1} \bigg)\nonumber\\
&\qquad-\mathbb{E}\!\bigg(\int_0^1\omega_{k_1} \operatorname{Tr}\bigg(u_{k_1}(x)
\big(H^\omega-E\big)^{-m-1} _{(\omega_{k_1}\to t\omega_{k_1})}\bigg) dt \,\bigg|\, \mathcal{F}^1_{k_1} \bigg)\bigg|\leq Ce^{-\beta \tilde{\ell}_L/2}.
\end{align}
Now it follows from (\ref{apr-fl}) and (\ref{bbd-sq}) that
\begin{align}
\label{decay-sq-ex}
&\sum_{k_1\in\big(B_L\big)_{\tilde{\ell}_L}^\circ}\bigg|\bigg\{\mathbb{E}\!\bigg[\mathbb{E}\!\bigg(\int_0^1\omega_{k_1} \operatorname{Tr}\bigg(u_{k_1}(x)
\big(H^\omega_{\Lambda_L, D} -E\big)^{-m-1} _{(\omega_{k_1}\to t\omega_{k_1})}\bigg) dt\,\bigg|\, \mathcal{F}^1_{k_1} \bigg)\nonumber\\
&\qquad\qquad
-\mathbb{E}\!\bigg(\int_0^1\omega_{k_1} \operatorname{Tr}\bigg(u_{k_1}(x)
\big(H^\omega_{\Lambda_L, D} -E\big)^{-m-1} _{(\omega_{k_1}\to t\omega_{k_1})}\bigg) dt \,\bigg|\, \mathcal{F}^1_{k_1-1} \bigg)\bigg]^2\bigg\}\nonumber\\
&\qquad\qquad\qquad
-\bigg\{\mathbb{E}\!\bigg[\mathbb{E}\!\bigg(\int_0^1\omega_{k_1} \operatorname{Tr}\bigg(u_{k_1}(x)
\big(H^\omega -E\big)^{-m-1} _{(\omega_{k_1}\to t\omega_{k_1})}\bigg) dt\,\bigg|\, \mathcal{F}^1_{k_1} \bigg)\nonumber\\
&\qquad\qquad\qquad\qquad
-\mathbb{E}\!\bigg(\int_0^1\omega_{k_1} \operatorname{Tr}\bigg(u_{k_1}(x)
\big(H^\omega -E\big)^{-m-1} _{(\omega_{k_1}\to t\omega_{k_1})}\bigg) dt\,\bigg|\, \mathcal{F}^1_{k_1-1} \bigg)\bigg]^2\bigg\}\bigg|\nonumber\\
&\qquad\qquad\qquad\qquad\qquad \qquad\qquad\qquad\qquad
\leq C \big| \big(B_L\big)_{\tilde{\ell}_L}^\circ \big|e^{-\beta \tilde{\ell}_L/2}.
\end{align}
Substituting (\ref{mrt-d1}) in (\ref{v-d1}), we obtain
\begin{align}
\label{vumrt}
&\mathrm{Var}(Y_{D,L})\nonumber\\
&= \sum_{k_1=-L_\delta}^{L_\delta}\mathbb{E}\!\bigg[
m\mathbb{E}\!\bigg(\int_0^1\omega_{k_1} \operatorname{Tr}\bigg(u_{k_1}(x)
\big(H^\omega_{\Lambda_L, D} -E\big)^{-m-1} _{(\omega_{k_1}\to t\omega_{k_1})}\bigg) dt\,\bigg|\, \mathcal{F}^1_{k_1} \bigg)\nonumber\\
&\qquad\qquad
-m\mathbb{E}\!\bigg(\int_0^1\omega_{k_1} \operatorname{Tr}\bigg(u_{k_1}(x)
\big(H^\omega_{\Lambda_L, D} -E\big)^{-m-1} _{(\omega_{k_1}\to t\omega_{k_1})}\bigg) dt\,\bigg|\, \mathcal{F}^1_{k_1-1} \bigg)\bigg]^2\nonumber\\
&=\sum_{k_1\in\big(B_L\big)_{\tilde{\ell}_L}^\circ}\mathbb{E}\!\bigg[
m\mathbb{E}\!\bigg(\int_0^1\omega_{k_1} \operatorname{Tr}\bigg(u_{k_1}(x)
\big(H^\omega_{\Lambda_L, D} -E\big)^{-m-1} _{(\omega_{k_1}\to t\omega_{k_1})}\bigg) dt\,\bigg|\, \mathcal{F}^1_{k_1} \bigg)\nonumber\nonumber\\
&\qquad\qquad
-m\mathbb{E}\!\bigg(\int_0^1\omega_{k_1} \operatorname{Tr}\bigg(u_{k_1}(x)
\big(H^\omega_{\Lambda_L, D} -E\big)^{-m-1} _{(\omega_{k_1}\to t\omega_{k_1})}\bigg) dt\,\bigg|\, \mathcal{F}^1_{k_1-1} \bigg)\bigg]^2\nonumber\\
&\qquad+\sum_{k_1\in\Lambda_{L_\delta}\setminus\big(B_L\big)_{\tilde{\ell}_L}^\circ}\mathbb{E}\!\bigg[
m\mathbb{E}\!\bigg(\int_0^1\omega_{k_1} \operatorname{Tr}\bigg(u_{k_1}(x)
\big(H^\omega_{\Lambda_L, D} -E\big)^{-m-1} _{(\omega_{k_1}\to t\omega_{k_1})}\bigg) dt\,\bigg|\, \mathcal{F}^1_{k_1} \bigg)\nonumber\nonumber\\
&\qquad\qquad\qquad
-m\mathbb{E}\!\bigg(\int_0^1\omega_{k_1} \operatorname{Tr}\bigg(u_{k_1}(x)
\big(H^\omega_{\Lambda_L, D} -E\big)^{-m-1} _{(\omega_{k_1}\to t\omega_{k_1})}\bigg) dt\,\bigg|\, \mathcal{F}^1_{k_1-1} \bigg)\bigg]^2\nonumber\\
&:=\Gamma^1_L+\mathcal{E}^1_L.
\end{align}
Owing to estimate~(\ref{decay-sq-ex}), we can replace the operator 
$H^\omega_{\Lambda_L, D}$ with $H^\omega$ in the expression 
$\operatorname{Tr}\!\left(u_{k_1}(x)\,(H^\omega_{\Lambda_L, D} - E)^{-m-1}_{(\omega_{k_1}\to t\omega_{k_1})}\right)$
whenever $k_1$ belongs to the interior $\big(B_L\big)_{\tilde{\ell}_L}^\circ$ as $L \to \infty$.
It is immediate from (\ref{decay-sq-ex}) and the above that
\begin{align}
\label{est-1}
&\lim_{L\to\infty}\frac{1}{|\Lambda_L|}\Gamma^1_L\nonumber\\
&=\lim_{L\to\infty}
\frac{1}{|\Lambda_L|}\sum_{k_1\in\big(B_L\big)_{\tilde{\ell}_L}^\circ}\mathbb{E}\!\bigg[
m\mathbb{E}\!\bigg(\int_0^1\omega_{k_1} \operatorname{Tr}\bigg(u_{k_1}(x)
\big(H^\omega -E\big)^{-m-1} _{(\omega_{k_1}\to t\omega_{k_1})}\bigg) dt \,\bigg|\, \mathcal{F}^1_{k_1} \bigg)\nonumber\nonumber\\
&\qquad\qquad\qquad
-m\mathbb{E}\!\bigg(\int_0^1\omega_{k_1} \operatorname{Tr}\bigg(u_{k_1}(x)
\big(H^\omega -E\big)^{-m-1} _{(\omega_{k_1}\to t\omega_{k_1})}\bigg) dt \,\bigg|\, \mathcal{F}^1_{k_1-1} \bigg)\bigg]^2\nonumber\\
&=\mathbb{E}\!\bigg[
m\mathbb{E}\!\bigg(\int_0^1\omega_{1} \operatorname{Tr}\bigg(u_{1}(x)
\big(H^\omega -E\big)^{-m-1} _{(\omega_{1}\to t\omega_{1})}\bigg)dt \,\bigg|\, \mathcal{F}^1_{1} \bigg)\nonumber\\
&\qquad\qquad\qquad
-m\mathbb{E}\!\bigg(\int_0^1\omega_{1} \operatorname{Tr}\bigg(u_{1}(x)
\big(H^\omega -E\big)^{-m-1} _{(\omega_{1}\to t\omega_{1})}\bigg) dt\,\bigg|\, \mathcal{F}^1_{0} \bigg)\bigg]^2.
\end{align}
In the above we have used the fact that $\displaystyle\lim_{L\to\infty}\frac{\#\big(B_L\big)_{\tilde{\ell}_L}^\circ}{|\Lambda_L|}=1$ and the Corollary \ref{smdscn}.
Now from (\ref{bbd-sq}), we have
\begin{align}
\label{est-2}
\big|\mathcal{E}^1_L\big|\leq \big| \Lambda_{L_\delta} \setminus\big(B_L\big)_{\tilde{\ell}_L}^\circ\big|\big(2C\|\omega_0\|_\infty \|u \|_\infty|S|\big)^2~~\text{and}~~
\lim_{L\to\infty}\frac{\mathcal{E}^1_L}{|\Lambda_L|}=0.
\end{align}
Use of (\ref{est-1}) and (\ref{est-2}) in (\ref{vumrt}) gives
\begin{align}
\label{vlmd1}
&\lim_{L\to\infty}\frac{1}{|\Lambda_L|}\mathrm{Var}(Y_{D,L})\nonumber\\
&\qquad=
\mathbb{E}\!\bigg[
m\mathbb{E}\!\bigg(\int_0^1\omega_{1} \operatorname{Tr}\bigg(u_{1}(x)
\big(H^\omega -E\big)^{-m-1} _{(\omega_{1}\to t\omega_{1})}\bigg) dt\,\bigg|\, \mathcal{F}^1_{1} \bigg)\nonumber\\
&\qquad\qquad\qquad
-m\mathbb{E}\!\bigg(\int_0^1\omega_{1} \operatorname{Tr}\bigg(u_{1}(x)
\big(H^\omega -E\big)^{-m-1} _{(\omega_{1}\to t\omega_{1})}\bigg) dt\,\bigg|\, \mathcal{F}^1_{0} \bigg)\bigg]^2.
\end{align}
\underline{Case \(d=2\):} Denote $n=(n_1,n_2)\in\mathbb{Z}^2$. Define $A^1_{k_1}=\{n \in \mathbb{Z}^2: n_1\leq k_1  \}$, $k_1\in\mathbb{Z}$ and $A^2_{(k_1,k_2)}=A^1_{k_1-1}\cup\{n\in \mathbb{Z}^2: n_1\leq k_1, n_2\leq k_2 \}$.  It is easy to observe that $A^2_{(k_1,k_2)}\setminus A^2_{(k_1,k_2-1)}=\{(k_1,k_2)\}$.  Define the $\sigma$-algebras $\mathcal{F}^1_{k_1}=\sigma\{\omega_n:n\in A^1_{k_1}  \}$ and
$\mathcal{F}^2_{(k_1,k_2)}=\sigma\{\omega_n: n\in A^2_{(k_1,k_2)} \}$.   It is observed that
$Y_{D,L}=\mathbb{E}\big(Y_{D,L}\big|\mathcal{F}^1_{L_\delta}  \big)$ and 
$\mathbb{E}\big(Y_{D,L}\big)=\mathbb{E}\big(Y_{D,L}\big|\mathcal{F}^1_{-L_\delta-1}  \big)$.
We can now write the variance of $Y_{D,L}$ as
\begin{equation}
\label{vd11-0}
\mathrm{Var}(Y_{D,L})
= \sum_{k_1=-L_\delta}^{L_\delta} 
\mathbb{E}\!\Bigg[ \mathbb{E}\!\big(Y_{D,L} \,\big|\, \mathcal{F}^1_{k_1} \big)
- \mathbb{E}\!\big(Y_{D,L} \,\big|\, \mathcal{F}^1_{k_1-1} \big)\Bigg]^2=\sum_{k_1=-L_\delta}^{L_\delta} \mathbb{E}\!\Big[Y^1_{L,k_1} \Big]^2  .
\end{equation}
In the above we denote $Y^1_{L,k_1}=\mathbb{E}\!\big(Y_{D,L} \,\big|\, \mathcal{F}^1_{k_1} \big)
- \mathbb{E}\!\big(Y_{D,L} \,\big|\, \mathcal{F}^1_{k_1-1} \big)$.  It is also clear that 
$Y^1_{L,k_1}=\mathbb{E}\big( Y^1_{L,k_1}\big|\mathcal{F}^2_{(k_1,L_\delta)} \big)$ and
$\mathbb{E}\big( Y^1_{L,k_1} \big)=\mathbb{E}\big( Y^1_{L,k_1}\big|\mathcal{F}^2_{(k_1,-L_\delta-1)} \big)=0$.  Now we can write each term of the above sum as 
\begin{align}
\label{vd12}
\mathbb{E}\!\Big[Y^1_{L,k_1} \Big]^2 &=\sum_{k_2=-L_\delta}^{L_\delta} 
\mathbb{E}\!\Bigg[ \mathbb{E}\!\big(Y^1_{L,k_1} \,\big|\, \mathcal{F}^2_{(k_1, k_2)} \big)
- \mathbb{E}\!\big(Y^1_{L,k_1}\,\big|\, \mathcal{F}^2_{(k_1, k_2-1)} \big)\Bigg]^2\nonumber\\
 &=\sum_{k_2=-L_\delta}^{L_\delta} 
\mathbb{E}\!\Bigg[ \mathbb{E}\!\big(Y_{D,L} \,\big|\, \mathcal{F}^2_{(k_1, k_2)} \big)
- \mathbb{E}\!\big(Y_{D,L}\,\big|\, \mathcal{F}^2_{(k_1, k_2-1)} \big)\Bigg]^2\nonumber\\
&=\sum_{k_2=-L_\delta}^{L_\delta} 
\mathbb{E}\!\Bigg[Y^2_{L,k_1,k_2}\Bigg]^2,
\end{align}
here $Y^2_{L,k_1,k_2}=\mathbb{E}\!\big(Y_{D,L} \,\big|\, \mathcal{F}^2_{(k_1, k_2)} \big)
- \mathbb{E}\!\big(Y_{D,L}\,\big|\, \mathcal{F}^2_{(k_1, k_2-1)} \big)$.
In the second line above, we make use of the following inclusions of $\sigma$-algebras:
$\mathcal{F}^2_{(k_1, k_2)} \subset \mathcal{F}^1_{k_1}, 
\mathcal{F}^1_{k_1-1} \subset \mathcal{F}^2_{(k_1, k_2)}, 
\mathcal{F}^2_{(k_1, k_2-1)} \subset \mathcal{F}^1_{k_1}, 
\mathcal{F}^1_{k_1-1} \subset \mathcal{F}^2_{(k_1, k_2-1)}$.
We also use the fact that, for any random variable $X$ and $\sigma$-algebras 
$\mathcal{G}_1\subseteq\mathcal{G}_2$, we have
$\mathbb{E}\!\left( \,\mathbb{E}(X \mid \mathcal{G}_1)\,\middle|\, \mathcal{G}_2 \right) 
= \mathbb{E}\!\left( X \,\middle|\, \mathcal{G}_1  \right)=\mathbb{E}\!\left( \,\mathbb{E}(X \mid \mathcal{G}_2)\,\middle|\, \mathcal{G}_1 \right) $.
Now the expression (\ref{vd11-0}) becomes
\begin{align}
\label{vd11}
\mathrm{Var}(Y_{D,L})=\sum_{k_1,k_2=-L_\delta}^{L_\delta} 
\mathbb{E}\!\Bigg[Y^2_{L,k_1,k_2}\Bigg]^2.
\end{align}
Also, we have 
$\mathcal{F}^2_{(k_1, k_2)}\setminus\mathcal{F}^2_{(k_1, k_2-1)}=\sigma\{\omega_{(k_1,k_2)}\}$.  Applying the same method used in (\ref{mrt-d1}) gives
\begin{align}
\label{mrt-d2}
&Y^2_{L,k_1,k_2}\\
&=\mathbb{E}\!\big(Y_{D,L} \,\big|\, \mathcal{F}^2_{(k_1, k_2)} \big)
- \mathbb{E}\!\big(Y_{D,L}\,\big|\, \mathcal{F}^2_{(k_1, k_2-1)} \big)\nonumber\\
&=-m\mathbb{E}\!\bigg(\int_0^1\omega_{(k_1, k_2)} \operatorname{Tr}\bigg(u_{(k_1, k_2)}(x)\big(H^\omega_{\Lambda_L, D} -E\big)^{-m-1} _{(\omega_{(k_1, k_2)}\to t\omega_{(k_1, k_2)})}\bigg)dt \,\bigg|\, \mathcal{F}^2_{(k_1,k_2)} \bigg)\nonumber\\
&\qquad+m\mathbb{E}\!\bigg(\int_0^1\omega_{(k_1, k_2)} \operatorname{Tr}\bigg(u_{(k_1, k_2)}(x)\big(H^\omega_{\Lambda_L, D} -E\big)^{-m-1} _{(\omega_{(k_1, k_2)}\to t\omega_{(k_1, k_2)})}\bigg)dt \,\bigg|\, \mathcal{F}^2_{(k_1,k_2-1)} \bigg)\nonumber.
\end{align}
Substituting (\ref{mrt-d2}) and (\ref{vd12}) in (\ref{vd11}), we obtain
\begin{align}
\label{vr-d2-fr}
&\mathrm{Var}\big(Y_{D,L} \big)\nonumber\\
&=\sum_{k_1,k_2=-L_\delta}^{L_\delta}
\mathbb{E}\bigg[  -m\mathbb{E}\!\bigg(\int_0^1\omega_{(k_1, k_2)} \operatorname{Tr}\bigg(u_{(k_1, k_2)}(x)\big(H^\omega_{\Lambda_L, D} -E\big)^{-m-1} _{(\omega_{(k_1, k_2)}\to t\omega_{(k_1, k_2)})}\bigg) dt \,\bigg|\, \mathcal{F}^2_{(k_1,k_2)} \bigg)\nonumber\\
&\qquad+m\mathbb{E}\!\bigg(\int_0^1\omega_{(k_1, k_2)} \operatorname{Tr}\bigg(u_{(k_1, k_2)}(x)\big(H^\omega_{\Lambda_L, D} -E\big)^{-m-1} _{(\omega_{(k_1, k_2)}\to t\omega_{(k_1, k_2)})}\bigg) dt\,\bigg|\, \mathcal{F}^2_{(k_1,k_2-1)} \bigg)\bigg]^2\nonumber\\
&=\sum_{\big(k_1,k_2\big)\in \big(B_L\big)_{\tilde{\ell}_L}^\circ}
\mathbb{E}\bigg[  -m\mathbb{E}\!\bigg(\int_0^1\omega_{(k_1, k_2)} \operatorname{Tr}\bigg(u_{(k_1, k_2)}(x)\big(H^\omega_{\Lambda_L, D} -E\big)^{-m-1} _{(\omega_{(k_1, k_2)}\to t\omega_{(k_1, k_2)})}\bigg) dt \,\bigg|\, \mathcal{F}^2_{(k_1,k_2)} \bigg)\nonumber\\
&\qquad+m\mathbb{E}\!\bigg(\int_0^1\omega_{(k_1, k_2)} \operatorname{Tr}\bigg(u_{(k_1, k_2)}(x)\big(H^\omega_{\Lambda_L, D} -E\big)^{-m-1} _{(\omega_{(k_1, k_2)}\to t\omega_{(k_1, k_2)})}\bigg) dt\,\bigg|\, \mathcal{F}^2_{(k_1,k_2-1)} \bigg)\bigg]^2
  \nonumber\\
&+\sum_{\big(k_1,k_2\big)\in\Lambda_{L_\delta}\setminus \big(B_L\big)_{\tilde{\ell}_L}^\circ}
\mathbb{E}\bigg[  -m\mathbb{E}\!\bigg(\int_0^1\omega_{(k_1, k_2)} \operatorname{Tr}\bigg(u_{(k_1, k_2)}(x)\big(H^\omega_{\Lambda_L, D} -E\big)^{-m-1} _{(\omega_{(k_1, k_2)}\to t\omega_{(k_1, k_2)})}\bigg) dt\,\bigg|\, \mathcal{F}^2_{(k_1,k_2)} \bigg)\nonumber\\
&\qquad+m\mathbb{E}\!\bigg(\int_0^1\omega_{(k_1, k_2)} \operatorname{Tr}\bigg(u_{(k_1, k_2)}(x)\big(H^\omega_{\Lambda_L, D} -E\big)^{-m-1} _{(\omega_{(k_1, k_2)}\to t\omega_{(k_1, k_2)})}\bigg) dt\,\bigg|\, \mathcal{F}^2_{(k_1,k_2-1)} \bigg)\bigg]^2 \nonumber\\ 
&=\Gamma_L^2+\mathcal{E}^2_L.
\end{align}
With the help of Corollary \ref{smdscn}, the exact same argument as in (\ref{est-1}) gives\begin{align}
\label{estd2-1}
& \lim_{L\to\infty}\frac{1}{|\Lambda_L|}\Gamma^2_L\\
&=
\mathbb{E}\bigg[  -m\mathbb{E}\!\bigg(\int_0^1\omega_{(1, 1)} \operatorname{Tr}\bigg(u_{(1, 1)}(x)\big(H^\omega -E\big)^{-m-1} _{(\omega_{(1, 1)}\to t\omega_{(1, 1)})}\bigg)dt \,\bigg|\, \mathcal{F}^2_{(1,1)} \bigg)\nonumber\\
&~~+m\mathbb{E}\!\bigg(\int_0^1\omega_{(1, 1)} \operatorname{Tr}\bigg(u_{(1, 1)}(x)\big(H^\omega -E\big)^{-m-1} _{(\omega_{(1, 1)}\to t\omega_{(1, 1)})}\bigg) dt \,\bigg|\, \mathcal{F}^2_{(1,0)} \bigg)\bigg]^2\nonumber.
\end{align}
Arguing as in (\ref{est-2}), 
we obtain
\begin{equation}
\label{estd2-2}
\lim_{L\to\infty}\frac{1}{|\Lambda_L|}\mathcal{E}^2_L=0.
\end{equation}
Using (\ref{estd2-2}) and (\ref{estd2-1}) in  (\ref{vr-d2-fr}), we deduce
\begin{align}
\label{vrlmd2}
& \lim_{L\to\infty}\frac{1}{|\Lambda_L|}\mathrm{Var}\big(Y_{D,L}\big)\nonumber\\
&=
\mathbb{E}\bigg[  -m\mathbb{E}\!\bigg(\int_0^1\omega_{(1, 1)} \operatorname{Tr}\bigg(u_{(1, 1)}(x)\big(H^\omega -E\big)^{-m-1} _{(\omega_{(1, 1)}\to t\omega_{(1, 1)})}\bigg) dt\,\bigg|\, \mathcal{F}^2_{(1,1)} \bigg)\nonumber\\
&~~+m\mathbb{E}\!\bigg(\int_0^1\omega_{(1, 1)} \operatorname{Tr}\bigg(u_{(1, 1)}(x)\big(H^\omega -E\big)^{-m-1} _{(\omega_{(1, 1)}\to t\omega_{(1, 1)})}\bigg) dt \,\bigg|\, \mathcal{F}^2_{(1,0)} \bigg)\bigg]^2.
\end{align}
\underline{Case \(d\ge 3\):} Denote $n=(n_1,n_2,n_3,\ldots,n_d)\in\mathbb{Z}^d$ and define the following subsets of $\mathbb{Z}^d$;
\begin{equation*}
\begin{split}
A^1_{k_1}
&=\{n\in\mathbb{Z}^d:n_1\le k_1\},
\qquad
A^2_{(k_1,k_2)}
= A^1_{k_1-1}
\cup \{n\in\mathbb{Z}^d:n_1\le k_1,\; n_2\le k_2\},\\
A^3_{(k_1,k_2,k_3)}
&= A^2_{(k_1,k_2-1)}
\cup \{n\in\mathbb{Z}^d:n_1\le k_1,\; n_2\le k_2,\; n_3\le k_3\}.
\end{split}
\end{equation*}
Now, for any $4\leq r\leq d-1$, we can define successively
\begin{align}
\label{atmstg}
&A^r_{(k_1,k_2,k_3,\ldots,k_{r-1},k_r)}\\
&=A^{r-1}_{(k_1,k_2,k_3,\ldots,k_{r-2},k_{r-1}-1)}
\cup\{n\in\mathbb{Z}^d:n_1\leq k_1, n_2\leq k_2, \ldots, n_{r-1}\leq k_{r-1},n_r\leq k_r  \}.\nonumber
\end{align}
In particular, we write for $r=d$
\begin{align}
\label{dsbsts}
&A^d_{(k_1,k_2,k_3,\ldots,k_{d-1},k_d)}\\
&=A^{d-1}_{(k_1,k_2,k_3,\ldots,k_{d-2},k_{d-1}-1)}
\cup\{n\in\mathbb{Z}^d:n_1\leq k_1, n_2\leq k_2, \ldots, n_{d-1}\leq k_{d-1},n_d\leq k_d  \}\nonumber.
\end{align}
We note that $A^d_{(k_1,k_2,k_3,\ldots,k_{d-1},k_d)}\setminus A^d_{(k_1,k_2,k_3,\ldots,k_{d-1},k_d-1)}=\big\{\big(k_1,k_2,\ldots,k_{d-1},k_d \big)\big\}$.
Now, define the collection of $\sigma$-algebras associated with the subset 
$A^r_{(k_1,k_2,k_3,\ldots,k_{r-1},k_r)} \subset \mathbb{Z}^d$ (as defined above):
\begin{equation}
\label{msalgb}
\mathcal{F}^r_{(k_1,k_2,\ldots,k_{r-1},k_r)}
= \sigma\big(\omega_n : n \in A^r_{(k_1,k_2,\ldots,k_{r-1},k_r)}\big),
~~~1 \leq r \leq d-1.
\end{equation}
For $r=d$ we have
\begin{equation}
\label{dsgal}
\mathcal{F}^d_{(k_1,k_2,k_3,\ldots,k_{d-1},k_d)}=\sigma\big(\omega_n:n\in A^d_{(k_1,k_2,k_3,\ldots,k_{d-1},k_d)}\big).
\end{equation}
It is also clear from the above construction that $\mathcal{F}^r_{(k_1,k_2,\ldots,k_{r-1},k_r)}\subset \mathcal{F}^{r-1}_{(k_1,k_2,\ldots,k_{r-1})}$,  $ \mathcal{F}^{r-1}_{(k_1,k_2,\ldots,k_{r-1}-1)}\subset \mathcal{F}^r_{(k_1,k_2,\ldots,k_{r-1}, k_r)}$,  $\mathcal{F}^r_{(k_1,k_2,\ldots,k_{r-1},k_r-1)}\subset \mathcal{F}^{r-1}_{(k_1,k_2,\ldots,k_{r-1})}$ and \\ 
$ \mathcal{F}^{r-1}_{(k_1,k_2,\ldots,k_{r-1}-1)}\subset \mathcal{F}^r_{(k_1,k_2,\ldots,k_{r-1}, k_r-1)}$. We now define a family of random variables recursively via conditional expectations
\begin{equation}
\label{dfn}
\begin{split}
Y^1_{L,k_1}&=\mathbb{E}\big( Y_{D,L}\big|\mathcal{F}^1_{k_1} \big)-\mathbb{E}\big( Y_{D,L}\big|\mathcal{F}^1_{k_1-1} \big)\\
Y^2_{L,k_1,k_2}&=\mathbb{E}\big( Y^1_{L,k_1}\big|\mathcal{F}^2_{(k_1,k_2)} \big)-\mathbb{E}\big( Y^1_{L,k_1}\big|\mathcal{F}^2_{(k_1,k_2-1)} \big)\\
&=\mathbb{E}\big( Y_{D,L}\big|\mathcal{F}^2_{(k_1,k_2)} \big)-\mathbb{E}\big( Y_{D,L}\big|\mathcal{F}^2_{(k_1,k_2-1)}\big).
\end{split}
\end{equation}
We define, for $3 \leq r \leq d$,
\begin{align}
\label{mthcnex}
Y^r_{L,k_1,k_2,\ldots,k_{r-1},k_r}&=\mathbb{E}\big(Y^{r-1}_{L,k_1,k_2,\ldots,k_{r-1}} \big|\mathcal{F}^r_{(k_1,k_2,\ldots,k_{r-1},k_r)}\big)\nonumber\\
&\qquad\qquad-
\mathbb{E}\big(Y^{r-1}_{L,k_1,k_2,\ldots,k_{r-1}}   \big|\mathcal{F}^r_{(k_1,k_2,\ldots,k_{r-1},k_r-1)}\big)\nonumber\\
&=\mathbb{E}\big(Y_{D,L} \big|\mathcal{F}^r_{(k_1,k_2,\ldots,k_{r-1},k_r)}\big)-
\mathbb{E}\big(Y_{D,L} \big|\mathcal{F}^r_{(k_1,k_2,\ldots,k_{r-1},k_r-1)}\big).
\end{align}
We repeat the method used to obtain (\ref{vd11}) a total of $d$ times, and thereby obtain
\begin{equation}
\label{lstvrex}
\mathrm{Var}(Y_{D,L})=\sum_{k_1,k_2,\ldots,k_{d-1},k_d=-L_\delta}^{L_\delta}\mathbb{E}
\bigg[Y^d_{L,k_1,k_2,\ldots,k_{d-1},k_d}  \bigg]^2,
\end{equation}
here $Y^d_{L,k_1,k_2,\ldots,k_{d-1},k_d}=\mathbb{E}\big(Y_{D,L} \big|\mathcal{F}^d_{(k_1,k_2,\ldots,k_{d-1},k_d)}\big)-
\mathbb{E}\big(Y_{D,L} \big|\mathcal{F}^d_{(k_1,k_2,\ldots,k_{d-1},k_d-1)}\big)$. Also we have
$\mathcal{F}^d_{(k_1,k_2,\ldots,k_{d-1},k_d)}\setminus\mathcal{F}^d_{(k_1,k_2,\ldots,k_{d-1},k_d-1)}=\sigma\{\omega_{(k_1,k_2,\ldots,k_{d-1},k_d)}\}$.\\ 
Now using the same method as in (\ref{mrt-d2}), we get
\begin{align}
\label{mrtdcs}
&Y^d_{L,k_1,k_2,\ldots,k_{d-1},k_d}\\
&=\mathbb{E}\big(Y_{D,L} \big|\mathcal{F}^d_{(k_1,k_2,\ldots,k_{d-1},k_d)}\big)-
\mathbb{E}\big(Y_{D,L} \big|\mathcal{F}^d_{(k_1,k_2,\ldots,k_{d-1},k_d-1)}\big)\nonumber\\
&=-m\mathbb{E}\!\bigg(\int_0^1\omega_{(k_1, k_2,\ldots,k_{d-1},k_d)} \operatorname{Tr}\bigg(u_{(k_1, k_2,\ldots,k_{d-1},k_d)}(x)\nonumber\\
&\qquad\qquad\big(H^\omega_{\Lambda_L, D} -E\big)^{-m-1} _{(\omega_{(k_1, k_2,\ldots,k_{d-1},k_d)}\to t\omega_{(k_1, k_2,\ldots,k_{d-1},k_d)})}\bigg)dt \,\bigg|\, \mathcal{F}^d_{(k_1,k_2,\ldots,k_{d-1},k_d)} \bigg)\nonumber\\
&+m\mathbb{E}\!\bigg(\int_0^1\omega_{(k_1, k_2,\ldots,k_{d-1},k_d)} \operatorname{Tr}\bigg(u_{(k_1, k_2,\ldots,k_{d-1},k_d)}(x)\nonumber\\
&\qquad\qquad\big(H^\omega_{\Lambda_L, D} -E\big)^{-m-1} _{(\omega_{(k_1, k_2,\ldots,k_{d-1},k_d)}\to t\omega_{(k_1, k_2,\ldots,k_{d-1},k_d)})}\bigg)dt \,\bigg|\, \mathcal{F}^d_{(k_1,k_2,\ldots,k_{d-1},k_d-1)} \bigg)\nonumber.
\end{align}
Now, with the help of Corollary~\ref{smdscn}, the same method as in \eqref{vr-d2-fr}, \eqref{estd2-1}, \eqref{estd2-2}, and \eqref{vrlmd2} yields
\begin{align}
\label{lmvext}
&\lim_{L\to\infty}\frac{1}{|\Lambda_L|}\mathrm{Var}(Y_{D,L})\\
&=\mathbb{E}\bigg[-m\mathbb{E}\!\bigg(\int_0^1\omega_{\vec{1}_d} \operatorname{Tr}\bigg(u_{\vec{1}_d}(x)
\big(H^\omega -E\big)^{-m-1} _{(\omega_{\vec{1}_d}\to t\omega_{\vec{1}_d})}\bigg)dt \,\bigg|\, \mathcal{F}^d_{\vec{1}_d} \bigg)\nonumber\\
&\qquad\qquad+m\mathbb{E}\!\bigg(\int_0^1\omega_{\vec{1}_d} \operatorname{Tr}\bigg(u_{\vec{1}_d}(x)
\big(H^\omega -E\big)^{-m-1} _{(\omega_{\vec{1}_d}\to t\omega_{\vec{1}_d})}\bigg)dt \,\bigg|\, \mathcal{F}^d_{\vec{1}_{d-1,0}} \bigg)\bigg]^2.\nonumber
\end{align}
In the above $\vec{1}_d= (\underbrace { 1,1,\ldots,1,1,1}_{\scriptstyle d})\in\mathbb{Z}^d$ and
$\vec{1}_{d-1,0}= ( \underbrace { 1,1,\ldots,1,1}_{\scriptstyle d-1},0)\in\mathbb{Z}^d$.
Hence the lemma.
\end{proof}
\begin{rem}
\label{fnlmres}
An upper bound on $\displaystyle\limsup_{L\to\infty}\frac{1}{|\Lambda_L|}\mathrm{Var}(Y_{D,L})$ is established in Proposition~\ref{bdd-var-f} (below).
\end{rem}
\noindent We now introduce a family of centered random variables obtained by evaluating a Laurent polynomial in the spectral parameter on the finite-volume operators.  
These will serve as localized building blocks in our variance estimates.\\
For $m > d+1$, consider the Laurent polynomial  
\begin{equation}
\label{rec-pol}
P(x) = \frac{1}{(x-E)^m} \sum_{j=0}^p \frac{a_j}{(x-E)^j},~p\in\mathbb{N}\cup\{0\},
\end{equation}
where $a_j \in \mathbb{R}$, $x \in [-\|V\|_\infty, \infty)$, and $E < -\|V\|_\infty$.  
We define, for this $P(x)$ the following centered random variables:
\begin{align}
\label{rv-pol}
\begin{split}
Y_{P,\Lambda^B_{L,k}}(\omega) &:= \operatorname{Tr} P\big(H^\omega_{\Lambda^B_{L,k}, D}\big) 
- \mathbb{E}\!\left[ \operatorname{Tr} P\big(H^\omega_{\Lambda^B_{L,k},D}\big) \right], \\[0.3em]
Y_{P,\Lambda^S_{L,k}}(\omega) &:= \operatorname{Tr} P\big(H^\omega_{\Lambda^S_{L,k},D}\big) 
- \mathbb{E}\!\left[ \operatorname{Tr} P\big(H^\omega_{\Lambda^S_{L,k},D}\big) \right], \\[0.3em]
Y_{P,\Lambda^{r_{_L}}_L}(\omega) &:= \operatorname{Tr} P\big(H^\omega_{\Lambda^{r_{_L}}_{L}, D}\big) 
- \mathbb{E}\!\left[ \operatorname{Tr} P\big(H^\omega_{\Lambda^{r_{_L}}_L,D}\big) \right], \\[0.3em]
Y_{P,\Lambda_{3R}}(\omega) &:= \operatorname{Tr} P\big(H^\omega_{\Lambda_{3R}, D}\big) 
- \mathbb{E}\!\left[ \operatorname{Tr} P\big(H^\omega_{\Lambda_{3R}, D}\big) \right], \\[0.3em]
Y_{P,D,L}(\omega) &:=Y_{P,\Lambda_L}(\omega)= \operatorname{Tr} P\big(H^\omega_{\Lambda_L, D}\big) 
- \mathbb{E}\!\left[ \operatorname{Tr} P\big(H^\omega_{\Lambda_L, D}\big) \right].
\end{split}
\end{align}
Where $\Lambda^B_{L,k}$,  $\Lambda^S_{L,k}$,  $\Lambda^{r_{_L}}_{L}$,  $\Lambda_{3R}$ are defined in (\ref{boxes}).
\begin{rem}
\label{ind-lrpl}
In view of Remark~\ref{ind-resl}, $\big\{Y_{P,\Lambda_{L,k}^B}\big\}_{k=1}^{r_{_L}-1}$ forms an independent collection of random variables, and the same is true for $\big\{Y_{P,\Lambda_{L,k}^S}\big\}_{k=1}^{r_{_L}-1}$.
\end{rem}
\noindent We next establish the existence of the limiting variance of the normalized random variable $|\Lambda_L|^{-\frac{1}{2}}Y_{P,D,L}$ as $L\to\infty$.
\begin{cor}
\label{lm-vr-lor}
Consider the random variable \(Y_{P,D,L}\) as defined in (\ref{rv-pol}). 
Then the limiting variance of the normalized random variable 
\(\tfrac{1}{|\Lambda_L|^{1/2}} Y_{P,D,L}\) exists, and it is given by \begin{align}
\label{lmv-lor}
\sigma^2_{P,D}:=
&\lim_{L \to \infty} \frac{1}{|\Lambda_L|} \,\mathrm{Var}\!\big( Y_{P,D,L} \big)\\
&=\mathbb{E}\bigg[\mathbb{E}\!\bigg(\int_0^1 \omega_{\vec{1}_d} \,\operatorname{Tr}\!\bigg(u_{\vec{1}_d}(x)
\big(P'(H^\omega)\big) _{(\omega_{\vec{1}_d}\to t\omega_{\vec{1}_d})}\bigg)\,dt \,\bigg|\, \mathcal{F}^d_{\vec{1}_d} \bigg)\nonumber\\
&\qquad\qquad-\mathbb{E}\!\bigg(\int_0^1 \omega_{\vec{1}_d} \,\operatorname{Tr}\!\bigg(u_{\vec{1}_d}(x)\big(P'(H^\omega)\big) _{(\omega_{\vec{1}_d}\to t\omega_{\vec{1}_d})}\bigg)\,dt \,\bigg|\, \mathcal{F}^d_{\vec{1}_{d-1,0}} \bigg)\bigg]^2.\nonumber
\end{align}
\end{cor}
\begin{proof}
Since $P(x)=\frac{1}{(x-E)^m}\sum_{j=0}^p\frac{a_j}{(x-E)^j}$ the exact same method used to prove Lemma~\ref{exst-vr-resl} yields the existence of the limit together with its explicit form.
\end{proof}
\begin{rem}
An upper bound on $\displaystyle\limsup_{L\to\infty}\frac{1}{|\Lambda_L|}\mathrm{Var}(Y_{P,D,L})$, follows as a special case of Proposition~\ref{bdd-var-f} (below).
\end{rem}
\noindent As a direct consequence of Proposition~\ref{bdd-2-4}; we obtain the following bounds for the second and fourth moments of the random variables in~\eqref{rv-pol}.
\begin{cor}
\label{2-4-pol}
Under the same assumptions as in Proposition~\ref{bdd-2-4},  
the second and fourth moments of the random variables defined in~\eqref{rv-pol} satisfy the bounds
\begin{equation*}
\mathbb{E}\big[ |Y_{P,\Lambda}|^2 \big] \leq C|\Lambda|, ~~
\mathbb{E}\big[ |Y_{P,\Lambda}|^4 \big] \leq C|\Lambda|^2, 
~~ \Lambda \in \big\{ \Lambda^B_{L,k},\, \Lambda^S_{L,k},\, \Lambda^{r_{_L}}_L,\, \Lambda_{3R}, \Lambda_L\big\},
\end{equation*}
here $C>0$ is independent of $\Lambda$.
\end{cor}
\begin{proof}
By Minkowski's inequality and Proposition \ref{bdd-2-4}, the result follows.
\end{proof}
\begin{cor}
\label{v-er-0}
Consider the random variables $Y_{P,\Lambda}$,  $\Lambda \in \big\{ \Lambda^B_{L,k},\, \Lambda^S_{L,k},\, \Lambda^{r_{_L}}_L,\, \Lambda_{3R}, \Lambda_L\big\}$ as defined in \eqref{rv-pol}. Then we have  
\begin{equation}\label{eq-limit-Y}
\lim_{L\to\infty} \frac{1}{|\Lambda_L|}\, 
\mathbb{E} \Bigg[
\sum_{k=1}^{r_{_L}-1} Y_{P,\Lambda^S_{L,k}}
+ Y_{P,\Lambda_L^{r_{_L}}}
+ Y_{P,\Lambda_{3R}}
\Bigg]^2 = 0,
\end{equation}
and moreover, the set of limit points (in the sense of distribution) of the sequence of random variables $\displaystyle\left\{ \frac{1}{|\Lambda_L|^{\frac{1}{2}}}\sum_{k=1}^{r_{_L}-1} Y_{P,\Lambda^B_{L,k}} \right\}_L$
coincides with that of the sequence $\left\{
\displaystyle\frac{1}{|\Lambda_L|^{\frac{1}{2}}} \bigg(\sum_{k=1}^{r_{_L}-1} Y_{P,\Lambda^B_{L,k}}
+  \displaystyle\sum_{k=1}^{r_{_L}-1} Y_{P,\Lambda^S_{L,k}}
+ Y_{P,\Lambda_L^{r_{_L}}}
+ Y_{P,\Lambda_{3R}}\bigg)
\right\}_L$.
\end{cor}
\begin{proof}
Since $\{Y_{P,\Lambda^S_{L,k}}\}_{k=1}^{r_{_L}-1}$ is an independent collection of random variables, we obtain the inequality  
\begin{align}
\label{wk-tnd-0}
&\frac{1}{|\Lambda_L|}\, 
\mathbb{E} \Bigg[
\sum_{k=1}^{r_{_L}-1} Y_{P,\Lambda^S_{L,k}}
+ Y_{P,\Lambda_L^{r_{_L}}}
+ Y_{P,\Lambda_{3R}}
\Bigg]^2\nonumber\\
&\qquad\qquad \leq \frac{9}{|\Lambda_L|}\bigg(\sum_{k=1}^{r_{_L}-1}
\mathbb{E}\big[ Y_{P,\Lambda^S_{L,k} } \big]^2+\mathbb{E}\big[Y_{P,\Lambda_L^{r_{_L}}}\big]^2 + \mathbb{E}\big[Y_{P,\Lambda_{3R}}\big]^2\bigg)\nonumber\\
&\qquad\qquad \leq \frac{9C}{|\Lambda_L|}\bigg(\sum_{k=1}^{r_{_L}-1}\big|  \Lambda^S_{L,k}\big|+\big|\Lambda_L^{r_{_L}}\big|+\big|\Lambda_{3R}\big|\bigg)\nonumber\\
&\qquad\qquad = \mathcal{O}\big(M_L^{-1}+L^{-(1-\epsilon)}+L^{-d}\big).
\end{align}
In the last line above, we have used the volume estimate given in (\ref{vl-est-box}).  
From (\ref{wk-tnd-0}),  the limit (\ref{eq-limit-Y}) follows immediately.  
It then follows from (\ref{wk-tnd-0}) that  
\begin{align}
\label{wkt0}
\frac{1}{|\Lambda_L|^{\frac{1}{2}}} \bigg( \sum_{k=1}^{r_{_L}-1} Y_{P,\Lambda^S_{L,k}}
+ Y_{P,\Lambda_L^{r_{_L}}}
+ Y_{P,\Lambda_{3R}}\bigg) \xrightarrow[L\to\infty]{\text{in distribution}} 0.
\end{align}
Finally, by Slutsky's theorem, we conclude the second part of the corollary.
\end{proof}
\noindent As a direct consequence of Proposition~\ref{e-q-cn}, we have:
\begin{cor}
\label{e-v-cn-pol}
Let $Y_{P,\Lambda}$,  $\Lambda \in \big\{ \Lambda^B_{L,k},\, \Lambda^S_{L,k},\, \Lambda^{r_{_L}}_L,\, \Lambda_{3R}, \Lambda_L\big\}$, be the random variables defined in~\eqref{rv-pol}.  
Under the same assumptions as in Proposition~\ref{e-q-cn}, the following convergence of second moments holds:
\begin{equation}
\label{eqlcn}
\frac{1}{|\Lambda_L|} \, \mathbb{E} \bigg[
Y_{P,D,L} 
- \sum_{k=1}^{r_{_L}-1} Y_{P,\Lambda^B_{L,k}}
- \sum_{k=1}^{r_{_L}-1} Y_{P,\Lambda^S_{L,k}}
- Y_{P,\Lambda_L^{r_{_L}}}
- Y_{P,\Lambda_{3R}}
\bigg]^2 \xrightarrow[L \to \infty]{} 0 .
\end{equation}
\end{cor}
\begin{proof} Since $P(x)=\frac{1}{(x-E)^m}\sum_{j=0}^p\frac{a_j}{(x-E)^j}$, 
the proof follows immediately from Minkowski's inequality and Proposition~\ref{e-q-cn}.
\end{proof}
\begin{cor}
\label{vr-eq}
Let $Y_{P,\Lambda}$,  $\Lambda \in \big\{ \Lambda^B_{L,k},\, \Lambda^S_{L,k},\, \Lambda^{r_{_L}}_L,\, \Lambda_{3R}, \Lambda_L\big\}$, be the random variables defined in~\eqref{rv-pol}.  
Under the same assumptions as in Proposition~\ref{e-q-cn}, we have
\begin{align}
\label{smvral}
\sigma^2_{P,D} &= \lim_{L\to\infty} \frac{1}{|\Lambda_L|} \mathrm{Var}\big( Y_{P,D,L} \big)\\
&= \lim_{L\to\infty} \frac{1}{|\Lambda_L|} \mathbb{E} \bigg[
\sum_{k=1}^{r_{_L}-1} Y_{P,\Lambda^B_{L,k}}
+ \sum_{k=1}^{r_{_L}-1} Y_{P,\Lambda^S_{L,k}}
+ Y_{P,\Lambda_L^{r_{_L}}}
+ Y_{P,\Lambda_{3R}} \bigg]^2 \nonumber\\
&= \lim_{L\to\infty} \frac{1}{|\Lambda_L|} \mathbb{E} \bigg[
\sum_{k=1}^{r_{_L}-1} Y_{P,\Lambda^B_{L,k}} \bigg]^2 \nonumber
\end{align}
\end{cor}
\begin{proof}
The first equality is given in Corollary~\ref{lm-vr-lor}.
The second equality follows from (\ref{eqlcn}) and Proposition~\ref{cn-vr-rv}.
The third equality follows directly from Proposition~\ref{vr-sm-lm} together with (\ref{eq-limit-Y}).
\end{proof}
\noindent Because of (\ref{wkt0}) and (\ref{eqlcn}), to establish the limiting distribution of $|\Lambda_L|^{-\tfrac{1}{2}}Y_{P,D,L}$ it is sufficient to determine the limiting distribution of $\displaystyle|\Lambda_L|^{-\frac{1}{2}}\sum_{k}Y_{P,\Lambda^B_{L,k}}$.  For this, we make use of the following version of the Central Limit Theorem for triangular arrays.
\begin{thm}
\label{clt-trar}
Suppose that for each $n$, the sequence of independent real-valued random variables $\{ Z_{n,k} \}_{k=1}^{r_n}$ has zero mean and satisfies
\begin{equation}
\label{cond}
\lim_{n\to\infty}\frac{1}{\sigma_n^{2+\delta}}\sum_{k=1}^{r_n}\mathbb{E}\big[|Z_{n,k}|^{2+\delta}  \big]=0,~\text{for some}~\delta>0.
\end{equation}
Then we have,
$ \frac{S_n}{\sigma_n}\xrightarrow[n\to\infty]{\substack{in\\ distribution}}\mathcal{N}(0,1)$.
Where $S_n$ and $\sigma_n$ are given by $S_n=\sum_{k=1}^{r_n}Z_{n,k}~~\text{and}~~\\
\sigma^2_n=\sum_{k=1}^{r_n}\mathbb{E}\big[ Z_{n,k}^2\big]$
\end{thm}
\begin{proof}
The proof of the theorem is given in \cite[Theorem 27.3]{PB}.
\end{proof}
\noindent Next, we state a result that identifies the exact limiting distribution of
$|\Lambda_L|^{-1/2}\sum_k Y_{P,\Lambda^B_{L,k}}$.
\begin{prop}
\label{wk-eq-cn}
Let $Y_{P,\Lambda}$,  $\Lambda \in \big\{ \Lambda^B_{L,k},\, \Lambda^S_{L,k},\, \Lambda^{r_{_L}}_L,\, \Lambda_{3R}, \Lambda_L\big\}$, be the random variables defined in~\eqref{rv-pol}, and let \(\sigma^2_{P,D} \ge 0\) be as in~\eqref{smvral}.  
Under the same assumptions as in Proposition~\ref{e-q-cn}, the following convergence in distribution holds:
\begin{align}
\label{clt1}
\lim_{L\to\infty} \frac{1}{|\Lambda_L|^{1/2}} Y_{P,D,L}
= \lim_{L\to\infty} \frac{1}{|\Lambda_L|^{1/2}} \sum_{k=1}^{r_{_L}-1} Y_{P,\Lambda^B_{L,k}}
= \mathcal{N}\big(0, \sigma^2_{P,D} \big),
\end{align}
with the convention that \(\mathcal{N}(0,0) = 0\).
\end{prop}
\begin{proof}
To begin with, we first establish (\ref{clt1}) when $\sigma^2_{P,D}=0$. 
If $\sigma^2_{P,D}=0$, then from (\ref{smvral}) we have
\begin{equation*}
\lim_{L\to\infty} \frac{1}{|\Lambda_L|} \mathrm{Var}\big( Y_{P,D,L} \big)=
\lim_{L\to\infty}\frac{1}{|\Lambda_L|} \mathbb{E}\bigg[ \sum_{k=1}^{r_{_L}-1} Y_{P,\Lambda^B_{L,k}} \bigg]^2
=0,
\end{equation*}
which implies
\begin{equation*}
\lim_{L\to\infty}\frac{1}{|\Lambda_L|^{1/2}}\sum_{k=1}^{r_{_L}-1} Y_{P,\Lambda^B_{L,k}}=\lim_{L\to\infty}\frac{1}{|\Lambda_L|^{1/2}}Y_{P,D,L}=0,~\text{in distribution}.
\end{equation*} 
On the other hand, when $\sigma^2_{P,D}>0$ we proceed as follows. By Remark~\ref{ind-lrpl}, we define
\begin{equation}
\label{vrtr}
\sigma^2_L:=\mathbb{E}\bigg[ \sum_{k=1}^{r_{_L}-1} Y_{P,\Lambda^B_{L,k}} \bigg]^2
=\sum_{k=1}^{r_{_L}-1}\mathbb{E}\big[ Y_{P,\Lambda^B_{L,k}}^2 \big].
\end{equation}
From (\ref{smvral}) it follows that $\sigma^2_L/|\Lambda_L|\to\sigma^2_{P,D}>0$ as $L\to\infty$.
Next, we verify assumption (\ref{cond}) for $\delta=2$. We compute
\begin{align*}
\frac{1}{\sigma^4_L}\sum_{k=1}^{r_{_L}-1}\mathbb{E}\big[ Y_{P,\Lambda^B_{L,k}}^4 \big]
&\leq \frac{1}{\sigma^4_L}\sum_{k=1}^{r_{_L}-1}\big| \Lambda^B_{L,k} \big|^2 \leq \frac{1}{\sigma^4_L}\sum_{k=1}^{r_{_L}-1}\mathcal{O}\!\big((2M_L)^{2d} k^{2d-2} \big) \\
&\leq \frac{|\Lambda_L|^2}{\sigma^4_L}\cdot \frac{1}{|\Lambda_L|^2}\,\mathcal{O}\!\big(M_L^{2d}r_{_L}^{2d-1}\big) \to 0 \quad \text{as } L\to\infty.
\end{align*}
In the first and second inequalities of the above, we used Corollary~\ref{2-4-pol} and estimate (\ref{vl-est-box}), respectively.  
By Theorem~\ref{clt-trar}, we conclude that
\[
\frac{1}{\sigma_L} \sum_{k=1}^{r_{_L}-1} Y_{P,\Lambda^B_{L,k}}
\xrightarrow[L\to\infty]{\text{in distribution}}\mathcal{N}(0, 1).
\]
Since $\sigma_L/|\Lambda_L|^{1/2}\to \sigma_{P,D}$ as $L\to\infty$, we equivalently obtain\begin{align}
\label{nrtrar}
\frac{1}{|\Lambda_L|^{1/2}}\sum_{k=1}^{r_{_L}-1} Y_{P,\Lambda^B_{L,k}}
=\frac{\sigma_L}{|\Lambda_L|^{1/2}}  \frac{1}{\sigma_L} \sum_{k=1}^{r_{_L}-1} Y_{P,\Lambda^B_{L,k}}
\xrightarrow[L\to\infty]{\text{in distribution}}\mathcal{N}\!\big(0,\sigma^2_{P,D}\big).
\end{align}
Finally, using (\ref{eqlcn}) and Proposition~\ref{cn-vr-rv}, we write, in distribution,
\begin{align}
\label{cltfrresl}
\lim_{L\to\infty }\frac{1}{|\Lambda_L|^{1/2}}Y_{P,D,L}
&= \lim_{L\to\infty }\frac{1}{|\Lambda_L|^{1/2}}
\Bigg( \sum_{k=1}^{r_{_L}-1} Y_{P,\Lambda^B_{L,k}}
+ \sum_{k=1}^{r_{_L}-1} Y_{P,\Lambda^S_{L,k}}
+ Y_{P,\Lambda_L^{r_{_L}}}
+ Y_{P,\Lambda_{3R}} \Bigg)\nonumber \\
&= \lim_{L\to\infty }\frac{1}{|\Lambda_L|^{1/2}}\sum_{k=1}^{r_{_L}-1} Y_{P,\Lambda^B_{L,k}} \nonumber \\
&= \mathcal{N}\!\big(0,\sigma^2_{P,D}\big). 
\end{align}
In the second and third lines of the above we used Corollary~\ref{v-er-0} and (\ref{nrtrar}), respectively.  
Thus, (\ref{clt1}) follows from (\ref{nrtrar}) and (\ref{cltfrresl}).
\end{proof}
\noindent Thus, the above result establishes the convergence (\ref{clt-cn}) when the function $f$ is a Laurent polynomial of the form (\ref{rec-pol}).

\section{CLT for test functions in 
\texorpdfstring{$C^1_{d,0}\big[-\|V\|_\infty,\infty\big)$}
{C1 d,0 [-V infinity, infinity)}}
Here we present the proof of our main result, Theorem~\ref{mnthm}.
To prove the CLT~\eqref{clt-cn} for more general test functions in 
$C^1_{d,0}[-\|V\|_\infty, \infty)$ (as in Definition \ref{defcl}), we make use of the following known result.
\begin{thm}
\label{app-con-thm}
Let $\{Z_L\}_{L=1}^\infty$ and $\{Y_{n,L}\}_{n,L=1}^\infty$ be real-valued random variables. Assume that
\begin{center}
\begin{enumerate}
\item [(a)] $Y_{n,L}\xrightarrow[L\to\infty]{\text{in distribution}}Y_n$,~ for each fixed $n$.\\
\item [(b)] $Y_n\xrightarrow[n\to\infty]{\text{in distribution}}Y$.\\
\item [(c)] For each $\delta>0$,  $\displaystyle \lim_{n\to\infty}\limsup_{L\to\infty}\mathbb{P}\bigg( \big|Y_{n,L}-Z_L  \big|\ge \delta \bigg)=0$.
\end{enumerate}
\end{center}
Then $Z_L\xrightarrow[L\to\infty]{\text{in distribution}}Y$. 
\end{thm}
\begin{proof}
A proof of this result can be found in \cite[Theorem 25.5]{PB}.
\end{proof}
\noindent To apply Theorem \ref{app-con-thm}, we first need some preliminary results. 
In particular, we require an estimate of how small the limiting variance 
$\sigma^2_{f,D}$ (see~(\ref{lim-vr})) becomes as $\|f\|_\infty$ decreases. 
To verify assumption~(c) of the theorem, we show that every function in 
$C_{d,0}^1[-\|V\|_\infty, \infty)$ can be approximated in the supremum norm 
by Laurent polynomials of the form~(\ref{rec-pol}).
\begin{defin}
\label{sbcl}
Fix a real number $E < -\|V\|_\infty$.  
A real-valued Laurent polynomial $P$ on $[-\|V\|_\infty, \infty)$ is said to belong to the class $\mathcal{A}_{E, d}\big[-\|V\|_\infty, \infty\big)$ if it can be written in the form
\[
P(x) = \frac{1}{(x-E)^m} \sum_{k=0}^p \frac{a_k}{(x-E)^k},
~~m \in \mathbb{N} \cap (d+1, \infty),  a_k \in \mathbb{R}~ \&~ p \in \mathbb{N} \cup \{0\}.
\]
\end{defin}
\noindent It follows from the Stone--Weierstrass Theorem that the above collection of
functions is dense, with respect to the supremum norm, in the space of all
continuous functions vanishing at infinity.
\begin{prop}
\label{st-wes}
Let $C_0\big[-\|V\|_\infty, \infty\big)$ denote the set of all continuous real-valued functions on $\big[-\|V\|_\infty, \infty\big)$ that vanish at infinity.  
Then the algebra $\mathcal{A}_{E,d}\big[-\|V\|_\infty, \infty\big)$ is dense in $C_0\big[-\|V\|_\infty, \infty\big)$ with respect to the supremum norm.
\end{prop}
\begin{proof}
The proof follows from the Stone–Weierstrass Theorem, see \cite[8.3 Corollary]{conwy}.
\end{proof}
\noindent We proceed to estimate the upper bound for the limsup of
the variance of the random variable $Y_{f,D,L}$ in terms of the supremum norm of the function $f$.
\begin{prop}
\label{bdd-var-f}
Let $f \in C^1_{d,0}[-\,\|V\|_\infty,\infty)$ (see Definition~\ref{defcl}), and consider the random variable $Y_{f,D,L}$ as in~\eqref{rv}.  
Under the same assumptions as in Proposition~\ref{e-q-cn}, we have
\begin{equation}
\label{vr-est}
\bar{\sigma}^2_{f,D} := \limsup_{L\to\infty}\frac{1}{|\Lambda_L|}\,\mathrm{Var}\!\big( Y_{f,D,L} \big) \,\leq\, C \, \| \tilde{f} \|^2_\infty,
\end{equation}
where $\tilde{f}(x)=(x-E)^{\,1 + \lfloor d/2 \rfloor}\, f'(x)$,
and $C>0$ is independent of $L$.
\end{prop}
\begin{proof}
It follows from (\ref{rv}) and (\ref{indx-set}) that the random variable $Y_{f,D,L}$ depends only on $\{\omega_n\}_{n\in B_L}$.  
Let $\{\omega_{n_j}\}_{j=1}^{\#B_L}$ be an enumeration of $\{\omega_n\}_{n\in B_L}$.  
Define $\mathcal{F}_k=\sigma\big(\omega_{n_j}: 1\leq j\leq k \big)$,
the $\sigma$-algebra generated by the random variables $\{\omega_{n_j}: 1\leq j\leq k\}$, and let $\mathcal{F}_0=\{\emptyset,  \Omega \}$.  
Using the martingale difference method, we can write the variance of $Y_{f,D,L}$ as
\begin{align}
\label{up-bbd-vr}
\mathrm{Var} \big( Y_{f,D,L} \big)
&=\sum_{k=1}^{\#B_L}
\mathbb{E}\Big[\mathbb{E}\big(Y_{f,D,L}\mid\mathcal{F}_k\big)- \mathbb{E}\big(Y_{f,D,L}\mid\mathcal{F}_{k-1} \big) \Big]^2\\
&=\sum_{k=1}^{\#B_L}
\mathbb{E}\Big[\mathbb{E}\big(\operatorname{Tr}\big(f(H^\omega_{\Lambda_L,D}) \big)\mid\mathcal{F}_k\big)
- \mathbb{E}\big(\operatorname{Tr}\big(f(H^\omega_{\Lambda_L,D}) \big)\mid\mathcal{F}_{k-1} \big) \Big]^2.\nonumber
\end{align}
To estimate the martingale difference inside the above sum, we write
\begin{align}
\label{up-bd-mntf}
&\mathbb{E}\big(\operatorname{Tr}(f(H^\omega_{\Lambda_L,D}))\mid\mathcal{F}_k\big)
- \mathbb{E}\big(\operatorname{Tr}(f(H^\omega_{\Lambda_L,D}))\mid\mathcal{F}_{k-1}\big)\nonumber\\
&\quad= \mathbb{E}\big(\operatorname{Tr}(f(H^\omega_{\Lambda_L,D}))\mid\mathcal{F}_k\big)
- \mathbb{E}\big(\operatorname{Tr}(f(H^\omega_{\Lambda_L,D}))_{(\omega_{n_k}=0)}\mid\mathcal{F}_k\big)\nonumber\\
&\qquad+\mathbb{E}\big(\operatorname{Tr}(f(H^\omega_{\Lambda_L,D}))_{(\omega_{n_k}=0)}\mid\mathcal{F}_{k-1}\big)
- \mathbb{E}\big(\operatorname{Tr}(f(H^\omega_{\Lambda_L,D}))\mid\mathcal{F}_{k-1}\big)\nonumber\\
&\quad=\mathbb{E}\bigg(\int_0^1 \frac{d}{dt}\,\operatorname{Tr}\!\big(f(H^\omega_{\Lambda_L,D})\big)_{(\omega_{n_k}\to t\omega_{n_k})}\,dt \,\Big|\,\mathcal{F}_{k}\bigg)\nonumber\\
&\qquad-\mathbb{E}\bigg(\int_0^1 \frac{d}{dt}\,\operatorname{Tr}\!\big(f(H^\omega_{\Lambda_L,D})\big)_{(\omega_{n_k}\to t\omega_{n_k})}\,dt \,\Big|\,\mathcal{F}_{k-1}\bigg)\nonumber\\
&\quad=\mathbb{E}\bigg(\omega_{n_k}\int_0^1 \operatorname{Tr}\!\big( u_{n_k}(x)f'(H^\omega_{\Lambda_L,D}) \big)_{(\omega_{n_k}\to t\omega_{n_k})}\,dt \,\Big|\,\mathcal{F}_{k}\bigg)\nonumber\\
&\qquad-\mathbb{E}\bigg(\omega_{n_k}\int_0^1 \operatorname{Tr}\!\big( u_{n_k}(x)f'(H^\omega_{\Lambda_L,D}) \big)_{(\omega_{n_k}\to t\omega_{n_k})}\,dt \,\Big|\,\mathcal{F}_{k-1}\bigg).
\end{align}
Let $S_{n_k}=\operatorname{supp}(u_{n_k})$. Then, by Corollary \ref{tr-est-fn}, we obtain
\begin{align}
\label{use-cor}
&\Big| \mathbb{E}\big(\operatorname{Tr}(f(H^\omega_{\Lambda_L,D}))\mid\mathcal{F}_k\big)
- \mathbb{E}\big(\operatorname{Tr}(f(H^\omega_{\Lambda_L,D}))\mid\mathcal{F}_{k-1}\big)\Big|
\nonumber\\
&\qquad\qquad \leq 2C \|\omega_{n_k} \|_\infty\|u_{n_k}\|_\infty \|\tilde{f}\|_\infty \,|S_{n_k}|=2C\|\omega_0 \|_\infty \|u\|_\infty \|\tilde{f}\|_\infty \,|S|.
\end{align}
Here $\tilde{f}(x)=(x-E)^{\,1 + \lfloor d/2 \rfloor}\, f'(x)$, with $E<-\|V\|_\infty$.  
In the last equality, we used the facts that $\{\omega_n\}_{n \in \mathbb{Z}^d}$ are i.i.d., $u_{n_k} = u(x-n_k)$, and $S = \operatorname{supp}(u)$.
Substituting (\ref{use-cor}) and (\ref{up-bd-mntf}) into (\ref{up-bbd-vr}), we obtain
\begin{equation}
\label{fnl-est}
\mathrm{Var} \big( Y_{f,D,L} \big)\leq 4 \big(\#B_L\big) C^2\|\tilde{f}\|_\infty^2,
\quad C>0~\text{independent of}~L.
\end{equation}
Since $\frac{\#B_L}{|\Lambda_L|}\to 1$ as $L\to\infty$, we obtain (\ref{vr-est}) from (\ref{fnl-est}).
\end{proof}
\noindent Consider the random variable $Y_{f,N,L}$ associated with the Neumann restriction $H^\omega_{\Lambda_L, N}$ as defined in (\ref{rv}). Then, the variance of $|\Lambda_L|^{-1/2} Y_{f,N,L}$ admits a bound analogous to that of the Dirichlet case, as given in (\ref{vr-est}).
\begin{cor}
Let $f \in C^1_{d,0}[-\,\|V\|_\infty,\infty)$ (see Definition~\ref{defcl}), and consider the random variable $Y_{f,N,L}$ as in~\eqref{rv}.  
Under the same assumptions as in Proposition~\ref{bdd-var-f}, we have
\begin{equation}
\label{vr-estn}
\bar{\sigma}^2_{f,N} := \limsup_{L\to\infty}\frac{1}{|\Lambda_L|}\,\mathrm{Var}\!\big( Y_{f,N,L} \big) \,\leq\, C \, \| \tilde{f} \|^2_\infty,
\end{equation}
where $\tilde{f}(x)=(x-E)^{\,1 + \lfloor d/2 \rfloor}\, f'(x)$,
and $C>0$ is independent of $L$.
\end{cor}
\begin{proof}
The proof proceeds exactly as in the Dirichlet case, as presented in Proposition~\ref{bdd-var-f}.
\end{proof}
\begin{rem}
\label{appr-fun-pol}
Let $f \in C^1_{d,0}\big[-\|V\|_\infty,\infty\big)$.  
Define $\tilde{f}(x) = (x-E)^{\,1 + \lfloor d/2 \rfloor}\, f'(x)$. Then $\tilde{f}\in C_0\big[-\|V\|_\infty,\infty \big)$, and by Proposition~\ref{st-wes}, there exists a sequence of Laurent polynomials $\{P_n\}_{n=1}^\infty$ in $\mathcal{A}_{E,d}\big[-\|V\|_\infty, \infty\big)$ such that  
$\| P_n - \tilde{f} \|_\infty \to 0$ as $n \to \infty$.  
Moreover, we can choose $Q_n \in \mathcal{A}_{E,d}\big[-\|V\|_\infty, \infty\big)$ satisfying  
$Q_n'(x) = (x-E)^{-\left( 1 + \lfloor d/2 \rfloor \right)} P_n(x)$.  
Since $|x-E| > -\|V\|_\infty - E$ for all $x \in \big[-\|V\|_\infty,\infty\big)$ and $E < -\|V\|_\infty$, it follows that $\| Q_n' - f' \|_\infty \to 0$ as $n \to \infty$.
\end{rem}
\noindent The above choice of $Q_n$ will later allow us to show that the limit~(\ref{vr-est}) actually exists, not merely its $\limsup$.
\begin{lem}
\label{lmvrfrf}
Let \(Q_n \in \mathcal{A}_{E,d}\big[-\|V\|_\infty, \infty\big)\) be as defined in Remark~\ref{appr-fun-pol},  let \(Y_{Q_n,D,L}\) be defined as in~\eqref{rv-pol},  and let \(\sigma^2_{Q_n, D}\) be as given in Corollary~\ref{lm-vr-lor}.  Further, let \(f \in C^1_{d,0}\big[-\|V\|_\infty, \infty\big)\) and consider \(Y_{f,D,L}\) as in~\eqref{rv}.  Then
\begin{align}
\label{lm-f-lv}
\sigma^2_{f,D} :&= \lim_{L\to\infty} \frac{1}{|\Lambda_L|} \mathrm{Var}(Y_{f,D,L}) = \lim_{n\to\infty} \sigma^2_{Q_n,D} \nonumber\\
&= \mathbb{E}\!\Bigg[\omega_{\vec{1}_d}\,
  \mathbb{E}\!\bigg(
    \int_0^1 
    \operatorname{Tr}\!\Big( u_{\vec{1}_d}(x)\,
    f'(H^\omega)_{(\omega_{\vec{1}_d}\to t\omega_{\vec{1}_d})} \Big)\,dt
  \;\Big|\; \mathcal{F}^1_{\vec{1}_d}  \bigg) \nonumber\\
&\qquad
- \mathbb{E}\!\bigg(
    \int_0^1 \omega_{\vec{1}_d}\,
    \operatorname{Tr}\!\Big( u_{\vec{1}_d}(x)\,
    f'(H^\omega)_{(\omega_{\vec{1}_d}\to t\omega_{\vec{1}_d})} \Big)\,dt
    \;\Big|\; \mathcal{F}^1_{\vec{1}_{d-1,0}}
  \bigg)
\Bigg]^2.
\end{align}
\end{lem}
\begin{proof}
First, we observe that 
\[
Q'_n(x)-f'(x)=(x-E)^{-\left( 1 + \lfloor d/2 \rfloor \right)}\,(P_n(x)-\tilde{f}(x)), 
\quad E<-\|V\|_\infty.
\]
Therefore,  
\begin{align*}
\big(Q'_n(H^\omega)-f'(H^\omega)\big)_{(\omega_{\vec{1}_d}\to t\omega_{\vec{1}_d})}
&=\Big(\,(H^\omega-E)^{-\left( 1 + \lfloor d/2 \rfloor \right)}\,(P_n(H^\omega)-\tilde{f}(H^\omega))\Big)_{(\omega_{\vec{1}_d}\to t\omega_{\vec{1}_d})}.
\end{align*}
Since $V^\omega\ge -\|V\|_\infty$ and $E<-\|V\|_\infty$, we have from the Corollary \ref{tr-est-fn} that
\begin{align*}
 &\Big|\operatorname{Tr}\!\Big( u_{\vec{1}_d}(x)\,
    Q_n'(H^\omega)_{(\omega_{\vec{1}_d}\to t\omega_{\vec{1}_d})} \Big)\,
    - \operatorname{Tr}\!\Big( u_{\vec{1}_d}(x)\,
    f'(H^\omega)_{(\omega_{\vec{1}_d}\to t\omega_{\vec{1}_d})} \Big)\Big|\nonumber\\
    &\qquad\qquad= \Big|\operatorname{Tr}\!\Big( u_{\vec{1}_d}(x)\Big(\,(H^\omega-E)^{-\left( 1 + \lfloor d/2 \rfloor \right)}\,(P_n(H^\omega)-\tilde{f}(H^\omega))\Big)_{(\omega_{\vec{1}_d}\to t\omega_{\vec{1}_d})}\Big|\nonumber\\
    &\qquad\qquad\leq C\|u\|_\infty|S|~\|P_n-\tilde{f}\|_\infty,~~S=\operatorname{supp}(u).
\end{align*}
Since $\|P_n-\tilde{f}\|_\infty\to 0$ as $n\to\infty$, the above implies
\[
\lim_{n\to\infty}\operatorname{Tr}\!\Big( u_{\vec{1}_d}(x)\,
    Q_n'(H^\omega)_{(\omega_{\vec{1}_d}\to t\omega_{\vec{1}_d})} \Big)
    = \operatorname{Tr}\!\Big( u_{\vec{1}_d}(x)\,
    f'(H^\omega)_{(\omega_{\vec{1}_d}\to t\omega_{\vec{1}_d})} \Big), 
    ~ \text{a.e. }\omega.
\]
It now follows from (\ref{lmv-lor}) that
\begin{align}
\label{lmof-lor}
\lim_{n\to\infty} \sigma^2_{Q_n,D} 
&= \mathbb{E}\!\Bigg[\omega_{\vec{1}_d}\,
  \mathbb{E}\!\bigg(
    \int_0^1 
    \operatorname{Tr}\!\Big( u_{\vec{1}_d}(x)\,
    f'(H^\omega)_{(\omega_{\vec{1}_d}\to t\omega_{\vec{1}_d})} \Big)\,dt
  \;\Big|\; \mathcal{F}^1_{\vec{1}_d}  \bigg) \nonumber\\
&\quad
- \mathbb{E}\!\bigg(
    \int_0^1 \omega_{\vec{1}_d}\,
    \operatorname{Tr}\!\Big( u_{\vec{1}_d}(x)\,
    f'(H^\omega)_{(\omega_{\vec{1}_d}\to t\omega_{\vec{1}_d})} \Big)\,dt
    \;\Big|\; \mathcal{F}^1_{\vec{1}_{d-1,0}}
  \bigg)
\Bigg]^2.
\end{align}
Using Minkowski's inequality, we obtain
\begin{align}
\label{minkw1}
&\frac{1}{|\Lambda_L|^{1/2}}\bigg(\mathbb{E}\!\big[|Y_{Q_n,D,L}|^2 \big]\bigg)^{1/2}\nonumber\\
&\qquad\qquad\leq \frac{1}{|\Lambda_L|^{1/2}}\bigg(\mathbb{E}\!\big[|Y_{Q_n,D,L}-Y_{f,D,L}|^2 \big]\bigg)^{1/2}
+ \frac{1}{|\Lambda_L|^{1/2}}\bigg(\mathbb{E}\!\big[|Y_{f,D,L}|^2 \big]\bigg)^{1/2}\nonumber\\
&\qquad\qquad=\frac{1}{|\Lambda_L|^{1/2}}\bigg(\mathbb{E}\!\big[|Y_{(Q_n-f),D,L}|^2 \big]\bigg)^{1/2}
+ \frac{1}{|\Lambda_L|^{1/2}}\bigg(\mathbb{E}\!\big[|Y_{f,D,L}|^2 \big]\bigg)^{1/2}. 
\end{align}
Swapping the roles of $Y_{Q_n,D,L}$ and $Y_{f,D,L}$ gives
\begin{align}
\label{minkw2}
&\frac{1}{|\Lambda_L|^{1/2}}\bigg(\mathbb{E}\!\big[|Y_{f,D,L}|^2 \big]\bigg)^{1/2}\nonumber\\
&\qquad\qquad\leq \frac{1}{|\Lambda_L|^{1/2}}\bigg(\mathbb{E}\!\big[|Y_{Q_n,D,L}-Y_{f,D,L}|^2 \big]\bigg)^{1/2}
+ \frac{1}{|\Lambda_L|^{1/2}}\bigg(\mathbb{E}\!\big[|Y_{Q_n,D,L}|^2 \big]\bigg)^{1/2}\nonumber\\
&\qquad\qquad=\frac{1}{|\Lambda_L|^{1/2}}\bigg(\mathbb{E}\!\big[|Y_{(Q_n-f),D,L}|^2 \big]\bigg)^{1/2}
+ \frac{1}{|\Lambda_L|^{1/2}}\bigg(\mathbb{E}\!\big[|Y_{Q_n,D,L}|^2 \big]\bigg)^{1/2}. 
\end{align}
Let us define 
\[
\tilde{g}_n(x)=(x-E)^{\,1 + \lfloor d/2 \rfloor}\,(Q_n'(x)-f'(x))=P_n(x)-\tilde{f}(x).
\]
Taking $\limsup$ (as $L\to\infty$) on both sides of (\ref{minkw1}) and (\ref{minkw2}), together with (\ref{lmv-lor}), we obtain
\begin{align}
\label{lsvf}
\Big|\sigma_{Q_n,D}-\limsup_{L\to\infty}\tfrac{1}{|\Lambda_L|^{1/2}}\big(\mathbb{E}\!\big[|Y_{f,D,L}|^2 \big]\big)^{1/2}\Big|
&\leq \limsup_{L\to\infty}\tfrac{1}{|\Lambda_L|^{1/2}}\mathbb{E}\!\big[|Y_{(Q_n-f),D,L}|^2 \big]^{1/2}\nonumber\\
&\leq C\|\tilde{g}_n \|_\infty=C\|P_n-\tilde{f}  \|_\infty.
\end{align}
Similarly, taking $\liminf$ (as $L\to\infty$) in (\ref{minkw1}) and (\ref{minkw2}) yields
\begin{align}
\label{livf}
\Big|\sigma_{Q_n,D}-\liminf_{L\to\infty}\tfrac{1}{|\Lambda_L|^{1/2}}\big(\mathbb{E}\!\big[|Y_{f,D,L}|^2 \big]\big)^{1/2}\Big|
&\leq C\|P_n-\tilde{f}  \|_\infty.
\end{align}
Taking the limit $n\to\infty$ in (\ref{livf}) and (\ref{lsvf}), we get
\[
\liminf_{L\to\infty}\frac{1}{|\Lambda_L|^{1/2}}\big(\mathbb{E}\!\big[|Y_{f,D,L}|^2 \big]\big)^{1/2}
=\lim_{n\to\infty}\sigma_{Q_n,D}
=\limsup_{L\to\infty}\frac{1}{|\Lambda_L|^{1/2}}\big(\mathbb{E}\!\big[|Y_{f,D,L}|^2 \big]\big)^{1/2}.
\]
Since all the relevant sequences are bounded, we conclude
\[
\lim_{L\to\infty}\frac{1}{|\Lambda_L|}\mathbb{E}\!\big(|Y_{f,D,L}|^2 \big)
=\lim_{n\to\infty}\sigma^2_{Q_n,D}.
\]
Finally, (\ref{lm-f-lv}) follows from (\ref{lmof-lor}).
\end{proof}
\noindent We now determine the weak limit of the random variable $\frac{1}{|\Lambda_L|^{\frac{1}{2}}}Y_{f,D,L}$ as $L$ becomes large.  For this, we make use of Theorem \ref{app-con-thm}.
\begin{lem}
\label{nr-appr}
Let $f \in C^1_{d,0}\big[-\|V\|_\infty,\infty\big)$, and let $Y_{f,D,L}$ be as in \eqref{rv}.  
Then under Hypothesis~\ref{hypo}, we have
\begin{align}
\label{onlynr}
\frac{1}{|\Lambda_L|^{\frac{1}{2}}}Y_{f,D,L}\xrightarrow[L\to\infty]{\text{in distribution}} \mathcal{N}\big(0,\sigma^2_{f,D}\big).
\end{align}
\end{lem}
\begin{proof}
Let $Q_n$ be as given in Remark~\ref{appr-fun-pol}. By Proposition~\ref{wk-eq-cn},
we obtain, for each fixed $n$,
\begin{align}
\label{con1}
\frac{1}{|\Lambda_L|^{\frac{1}{2}}}Y_{Q_n,D,L}\xrightarrow[L\to\infty]{\text{in distribution}}
\mathcal{N}\big(0,\sigma^2_{Q_n,D}\big).
\end{align}
By Lemma~\ref{lmvrfrf}, it follows that
\begin{align}
\label{con2}
\mathcal{N}\big(0,\sigma^2_{Q_n,D}\big)\xrightarrow[n\to\infty]{\text{in distribution}}
\mathcal{N}\big(0,\sigma^2_{f,D}\big).
\end{align}
Applying Markov's inequality, for each $\delta>0$ we obtain
\begin{align}
\label{con3}
\mathbb{P}\bigg(\frac{1}{|\Lambda_L|^{\frac{1}{2}}}\big|Y_{f,D,L}-Y_{Q_n,D,L}\big| \ge \delta \bigg)
&= \mathbb{P}\bigg(\frac{1}{|\Lambda_L|^{\frac{1}{2}}}\big|Y_{(f-Q_n),D,L}\big|\ge\delta\bigg)\nonumber \\
&\leq \frac{1}{\delta^2|\Lambda_L|}\,\mathbb{E}\big(\big|Y_{(f-Q_n),D,L}\big|^2\big).
\end{align}
By Proposition~\ref{bdd-var-f}, this implies
\begin{align}
\label{con31}
\limsup_{L\to\infty}\mathbb{P}\bigg(\frac{1}{|\Lambda_L|^{\frac{1}{2}}}\big|Y_{f,D,L}-Y_{Q_n,D,L}\big| \ge \delta \bigg) 
\leq \frac{C}{\delta^2} \|\tilde{g}_n\|^2_\infty,
\end{align}
where $\tilde{g}_n(x) = (x-E)^{\,1 + \lfloor d/2 \rfloor}\,\big(f'(x)-Q_n'(x)\big)$.
From the construction of $Q_n$ in Remark~\ref{appr-fun-pol}, we have $\tilde{g}_n(x) = \tilde{f}(x)-P_n(x)$, and moreover $\|\tilde{g}_n\|_\infty \to 0$ as $n \to \infty$.  
Hence, \eqref{con31} yields
\begin{align}
\label{con3t}
\lim_{n\to\infty}\limsup_{L\to\infty}\mathbb{P}\bigg(\frac{1}{|\Lambda_L|^{\frac{1}{2}}}\big|Y_{f,D,L}-Y_{Q_n,D,L}\big| \ge \delta \bigg) = 0, \quad \delta>0.
\end{align}
Combining \eqref{con1}, \eqref{con2}, and \eqref{con3t}, we verify all assumptions of Theorem~\ref{app-con-thm}. Therefore,
\begin{equation}
\label{frvnor}
\frac{1}{|\Lambda_L|^{\frac{1}{2}}}Y_{f,D,L}\xrightarrow[L\to\infty]{\text{in distribution}}
\mathcal{N}\big(0,\sigma^2_{f,D}\big).
\end{equation}
This completes the proof.
\end{proof}
\noindent We now characterize the situation in which the limiting variance $\sigma^2_{f,D}$ is strictly positive; namely, if the test function $f$ is strictly monotone, then the limiting variance is strictly positive.
\begin{lem}
\label{ps-vr}
Let \( f \in C^1_{d,0}\big[-\|V\|_\infty, \infty\big) \), and let \( \sigma^2_{f,D} \) be as in Lemma~\ref{lmvrfrf}.  
Suppose that \( f'(x)>0 \) or \( f'(x)<0 \) on entire domain \( \big[-\|V\|_\infty, \infty\big) \), and let \( u \ge 0 \) or \( u \leq 0 \) with \( \|u\|_2 \neq 0 \).  
Then, under Hypothesis~\ref{hypo}, we have \( \sigma^2_{f,D} > 0 \).
\end{lem}

\begin{proof}
First, observe that for any \( b \in \mathbb{R} \), the random operator \( H^\omega \) from~\eqref{model} can be written as
\begin{align}
\label{mdmdl}
H^\omega 
&= (i\nabla + A)^2 
+ b \sum_{n \in \mathbb{Z}^d} u(x - n) 
+ \sum_{n \in \mathbb{Z}^d} (\omega_n - b) u(x - n).
\end{align}
Since \( b \sum_{n \in \mathbb{Z}^d} u(x - n) \) is a bounded operator,  
all preceding results remain valid if we replace  
\( (i\nabla + A)^2 \) by \( (i\nabla + A)^2 + b \sum_{n \in \mathbb{Z}^d} u(x - n) \)  
and \( \omega_n \) by \( \omega_n - b \).  
Because \( \omega_0 \) is a bounded real random variable, it follows from~\eqref{mdmdl} that, without loss of generality, we may assume that \( \omega_0 \) takes both positive and negative values with nonzero probability.
\vspace{0.5em}
Now define a filtration of $\sigma$-algebras $\{\mathcal{K}_j\}_{j=1}^\infty$ by
\begin{equation}
\label{sgmfl}
\begin{split}
\mathcal{K}_0 &= \{\emptyset, \Omega\}, 
\quad \mathcal{K}_1 = \sigma\{\omega_{\vec{1}_d}\}, 
\quad \vec{1}_d = (1,1,\dots,1) \in \mathbb{Z}^d,\\
\mathcal{K}_j &= \sigma\big\{\omega_k : k \in \mathbb{Z}^d,\, \|k\|_\infty \le j \big\}, 
\quad j \ge 2.
\end{split}
\end{equation}
Let $\mathcal{K}_\infty$ denote the $\sigma$-algebra generated by $\cup_{j=0}^\infty \mathcal{K}_j$.  
Then $\mathbb{E}\big(Z \mid \mathcal{K}_\infty \big) = Z$.
\vspace{0.5em}
Define the random variable \( Z \) (with mean zero) by
\begin{align}
\label{cnrv}
Z &= \omega_{\vec{1}_d}\,
  \mathbb{E}\!\bigg(
    \int_0^1 
    \operatorname{Tr}\!\Big( u_{\vec{1}_d}(x)\,
    f'(H^\omega)_{(\omega_{\vec{1}_d}\to t\omega_{\vec{1}_d})} \Big)\,dt
  \;\Big|\; \mathcal{F}^1_{\vec{1}_d}  \bigg) \nonumber\\
&\quad
- \mathbb{E}\!\bigg(
    \int_0^1 \omega_{\vec{1}_d}\,
    \operatorname{Tr}\!\Big( u_{\vec{1}_d}(x)\,
    f'(H^\omega)_{(\omega_{\vec{1}_d}\to t\omega_{\vec{1}_d})} \Big)\,dt
    \;\Big|\; \mathcal{F}^1_{\vec{1}_{d-1,0}}
  \bigg).
\end{align}
Here, the $\sigma$-algebras $\mathcal{F}^1_{\vec{1}_d}$ and $\mathcal{F}^1_{\vec{1}_{d-1,0}}$ are defined as in (\ref{smag}).\\
Clearly, $\mathbb{E}\big(Z \mid \mathcal{K}_\infty\big) = Z$.  
By the martingale convergence theorem \cite[Theorem 35.6]{PB}, we have
\begin{equation}
\label{mrtcng}
\mathbb{E}\big(Z \mid \mathcal{K}_n\big) 
\xrightarrow[n\to\infty]{\text{a.e. } \omega} 
\mathbb{E}\big(Z \mid \mathcal{K}_\infty\big) = Z.
\end{equation}
From Corollary~\ref{tr-est-fn}, we obtain the uniform bound
\begin{equation}
\label{cnd-bd}
|Z| \le 2C\, \|\omega_0\|_\infty \|u\|_\infty |S| \|\tilde{f}\|_\infty\,,~~\tilde{f}(x)=(x-E)^{\,1 + \lfloor d/2 \rfloor}\ f'(x)~~\text{a.e. } \omega,
\end{equation}
which allows us to apply the dominated convergence theorem to conclude that
\begin{equation}
\label{recng}
\mathbb{E}\big[\big(\mathbb{E}(Z \mid \mathcal{K}_n)\big)^2\big]
\xrightarrow{n\to\infty}
\mathbb{E}\big[Z^2\big].
\end{equation}
Using the martingale difference technique, we can rewrite~Lemma \ref{lmvrfrf} as
\begin{align}
\label{stps}
\sigma^2_{f,D}
&= \mathbb{E}\big[Z^2\big]\nonumber\\
&= \lim_{n\to\infty} 
   \mathbb{E}\big[\big(\mathbb{E}(Z \mid \mathcal{K}_n)\big)^2\big],
   \quad \text{with } \mathbb{E}(Z \mid \mathcal{K}_0) = \mathbb{E}[Z] = 0 \nonumber\\
&= \lim_{n\to\infty} 
   \sum_{j=1}^n 
   \mathbb{E}\!\bigg[\Big(\mathbb{E}(Z \mid \mathcal{K}_j) 
   - \mathbb{E}(Z \mid \mathcal{K}_{j-1})\Big)^2\bigg]\nonumber\\
&\ge 
\mathbb{E}\!\bigg[\Big(\mathbb{E}(Z \mid \mathcal{K}_1) 
- \mathbb{E}(Z \mid \mathcal{K}_0)\Big)^2\bigg]
= \mathbb{E}\big[\big(\mathbb{E}(Z \mid \mathcal{K}_1)\big)^2\big].
\end{align}
We next show that the random variable $\mathbb{E}(Z \mid \mathcal{K}_1)$ is not identically zero.  
From~\eqref{cnrv} we have
\begin{align}
\label{excnz}
\mathbb{E}(Z \mid \mathcal{K}_1)
&= \omega_{\vec{1}_d}\,
  \mathbb{E}\!\bigg(
    \int_0^1 
    \operatorname{Tr}\!\Big( u_{\vec{1}_d}(x)\,
    f'(H^\omega)_{(\omega_{\vec{1}_d}\to t\omega_{\vec{1}_d})} \Big)\,dt
  \;\Big|\; \omega_{\vec{1}_d}  \bigg) \nonumber\\
&\quad
- \mathbb{E}\!\bigg(
    \int_0^1 \omega_{\vec{1}_d}\,
    \operatorname{Tr}\!\Big( u_{\vec{1}_d}(x)\,
    f'(H^\omega)_{(\omega_{\vec{1}_d}\to t\omega_{\vec{1}_d})} \Big)\,dt
  \bigg).
\end{align}
Define a measure $\nu_{\omega_{\vec{1}_d}}(\cdot)$ on $\big[-\|V\|_\infty, \infty\big)$ by
\begin{equation*}
\nu_{\omega_{\vec{1}_d}}(\cdot)
= \mathbb{E}\!\bigg(
    \int_0^1 
    \operatorname{Tr}\!\Big( u_{\vec{1}_d}(x)\,
    E_{H^\omega}(\cdot)_{(\omega_{\vec{1}_d}\to t\omega_{\vec{1}_d})} \Big)\,dt
  \;\Big|\; \omega_{\vec{1}_d} \bigg).
\end{equation*}
This defines a non-trivial measure on $\big[-\|V\|_\infty, \infty\big)$, since by Proposition~\ref{trnf} we have  
$\nu_{\omega_{\vec{1}_d}}\big(\big[-\|V\|_\infty, \infty\big)\big) = \infty$ for a.e. $\omega_{\vec{1}_d}$.
Without loss of generality, assume that \( f'(x) > 0 \) on \( \big[-\|V\|_\infty, \infty\big) \).  
Then for a.e. $\omega_{\vec{1}_d}$,
\begin{equation*}
\label{pstvty}
\mathbb{E}\!\bigg( 
\int_0^1 
\operatorname{Tr}\!\Big( u_{\vec{1}_d}(x)\,
f'(H^\omega)_{(\omega_{\vec{1}_d}\to t\omega_{\vec{1}_d})} \Big)\,dt
\;\Big|\; \omega_{\vec{1}_d} \bigg)
= \int f'(x)\,d\nu_{\omega_{\vec{1}_d}}(x) > 0.
\end{equation*}
Since $\omega_{\vec{1}_d}$ takes both positive and negative values with non-zero probability,  
and the expression above is strictly positive, the random variable
\[
\omega_{\vec{1}_d}\,
  \mathbb{E}\!\bigg(
    \int_0^1 
    \operatorname{Tr}\!\Big( u_{\vec{1}_d}(x)\,
    f'(H^\omega)_{(\omega_{\vec{1}_d}\to t\omega_{\vec{1}_d})} \Big)\,dt
  \;\Big|\; \omega_{\vec{1}_d} \bigg)
\]
takes both positive and negative values with nonzero probability.  
However, 
\[
\mathbb{E}\!\bigg[
    \int_0^1 \omega_{\vec{1}_d}\,
    \operatorname{Tr}\!\Big( u_{\vec{1}_d}(x)\,
    f'(H^\omega)_{(\omega_{\vec{1}_d}\to t\omega_{\vec{1}_d})} \Big)\,dt
  \bigg]
\]
is a constant real number.  
Therefore, by~\eqref{excnz}, $\mathbb{E}(Z \mid \mathcal{K}_1)$ is a nonzero bounded random variable (see~\eqref{cnd-bd}). In particular, since a nonzero square-integrable random variable has strictly positive second moment, we obtain
$\mathbb{E}\!\left[ \big( \mathbb{E}(Z \mid \mathcal{K}_1) \big)^2 \right] > 0$ .
From~\eqref{stps}, we conclude that \( \sigma^2_{f,D} > 0 \).\\~\\
For $u \leq 0$, we can rewrite (\ref{stps}) as
\[
\sigma_{f,D}^2 = \mathbb{E}[(-Z)^2] \geq \mathbb{E}\!\big[\left(\mathbb{E}[-Z \mid \mathcal{K}_1]\right)^2\big].
\]
Since $-u_{\vec{1}_d} \geq 0$, the same argument as in the case $u \geq 0$ implies that $\sigma_{f,D}^2 > 0$ whenever $u \leq 0$.
\end{proof}
\noindent We now proceed to show that the weak limits of the two sequences of random variables
$\big\{|\Lambda_L|^{-1/2} Y_{f,D,L}\big\}_L$ and $\big\{ |\Lambda_L|^{-1/2} Y_{f,N,L}\big\}_L$
coincide. First, we prove this for \(f \in \mathcal{A}_{E,d}\big[-\|V\|_\infty,\infty\big)\). Then, the denseness of \(\mathcal{A}_{E,d}\big[-\|V\|_\infty,\infty\big)\) in \(C^1_{d,0}\big[-\|V\|_\infty,\infty\big)\) allows us to extend the result to a general \(f\).\\
Let $P \in \mathcal{A}_{E,d}\big[-\|V\|_\infty,\infty\big)$ be as in Definition~\ref{sbcl}, and define
\begin{align}
\label{trnrs}
Y_{P,X,L}(\omega) := \operatorname{Tr}\Big(P(H^\omega_{\Lambda_L,X})\Big) - \mathbb{E}\Big[\operatorname{Tr}\big(P(H^\omega_{\Lambda_L,X})\big)\Big],~~X=D,N
\end{align}
where $P(x) = \frac{1}{(x-E)^m} \sum_{k=0}^p \frac{a_k} {(x-E)^{k}}$, $ m > d+1,  a_k \in \mathbb{R},  p \in \mathbb{N}\cup\{0\},  x \in \big[-\|V\|_\infty,\infty\big)$,  and $E < -\|V\|_\infty$. 
\begin{lem}
\label{trn=d}
Consider the random variables $Y_{P,N,L}$ and $Y_{P,D,L}$ as in (\ref{trnrs}) and (\ref{rv-pol}). Then, under Hypothesis \ref{hypo}, we have
\begin{equation}
\label{n=d}
\frac{1}{|\Lambda_L|}\mathbb{E}\bigg[\bigg(Y_{P,N,L}-Y_{P,D,L}  \bigg)^2\bigg]\xrightarrow{L\to\infty}0.
\end{equation}
\end{lem}
\begin{proof}
For $m > d+1$, we define the centered random variable $\mathcal{K}_L$ as
\begin{align}
\label{n-d}
\mathcal{K}_L(\omega) 
&= \operatorname{Tr}\!\big((H^{\omega}_{\Lambda_L,D} - E)^{-m}\big)
- \mathbb{E}\big[\operatorname{Tr}\!\big((H^{\omega}_{\Lambda_L, D} - E)^{-m}\big)\big] \nonumber\\
&\qquad\qquad \qquad- \operatorname{Tr}\!\big((H^{\omega}_{\Lambda_L, N} - E)^{-m}\big)
+ \mathbb{E}\big[\operatorname{Tr}\!\big((H^{\omega}_{\Lambda_L, N} - E)^{-m}\big)\big].
\end{align}
In view of Minkowski's inequality, it suffices to prove
\begin{equation}
\label{rsdf-nd}
\frac{1}{|\Lambda_L|}\,\mathbb{E}\Big[\big(\mathcal{K}_L(\omega)\big)^2\Big] \xrightarrow{L \to \infty} 0.
\end{equation}
Let $B_L$ be as defined in (\ref{indx-set}),  and enumerate its elements by $\{n_j\}_{j=1}^{\# B_L}$. Define a filtration of $\sigma$-algebras $\{\mathcal{D}_k\}_{k=0}^{\# B_L}$ by $\mathcal{D}_k = \sigma\{\omega_{n_j} : j \le k\}$ with $\mathcal{D}_0 = \{\emptyset, \Omega\}$.  
Using martingale difference techniques, we can write
\begin{align}
\label{n-dmrt}
\mathbb{E}\Big[\big(\mathcal{K}_L(\omega)\big)^2\Big]
= \sum_{k=0}^{\#B_L} \mathbb{E}\Big[\mathbb{E}(\mathcal{K}_L \mid \mathcal{D}_k) - \mathbb{E}(\mathcal{K}_L \mid \mathcal{D}_{k-1})\Big]^2.
\end{align}
Next, we estimate each martingale difference:
\begin{align}
 \label{est-df} 
 &\mathbb{E}\big(\mathcal{K}_L\big|\mathcal{D}_k\big)- \mathbb{E}\big(\mathcal{K}_L\big|\mathcal{D}_{k-1}\big)\\ 
 &\quad=\mathbb{E}\big(\mathcal{K}_L(\omega)\big|\mathcal{D}_k\big) - \mathbb{E}\big(\mathcal{K}_L(\omega:\omega_{n_k}=0)\big|\mathcal{D}_k\big)\nonumber\\ &\qquad\qquad-\mathbb{E}\big(\mathcal{K}_L(\omega)\big|\mathcal{D}_{k-1}\big) +\mathbb{E}\big(\mathcal{K}_L(\omega:\omega_{n_k}=0)\big|\mathcal{D}_{k-1}\big)\nonumber\\ 
 &\quad=\mathbb{E}\bigg(\int_0^1\frac{d}{dt}\mathcal{K}_L\big(\omega:\omega_{n_k}\to t\omega_{n_k}\big)dt\big| \mathcal{D}_k  \bigg)-\mathbb{E}\bigg(\int_0^1\frac{d}{dt}\mathcal{K}_L\big(\omega:\omega_{n_k}\to t\omega_{n_k}\big)dt\big| \mathcal{D}_{k-1} \bigg)\nonumber\\ 
 &\quad=-m\mathbb{E}\bigg(\omega_{n_k}\int_0^1\operatorname{Tr}\!\bigg(u_{n_k}\big(H^{\omega}_{\Lambda_L,D} - E\big)^{-m-1}_{(\omega_{n_k}\to t\omega_{n_k})}\bigg)dt\bigg| \mathcal{D}_k \bigg)\nonumber\\ 
 &\qquad\qquad +m\mathbb{E}\bigg(\omega_{n_k}\int_0^1\operatorname{Tr}\!\bigg(u_{n_k}\big(H^{\omega}_{\Lambda_L,D} - E\big)^{-m-1}_{(\omega_{n_k}\to t\omega_{n_k})}\bigg)dt\bigg| \mathcal{D}_{k-1} \bigg)\nonumber\\ 
&\qquad\qquad+m\mathbb{E}\bigg(\omega_{n_k}\int_0^1\operatorname{Tr}\!\bigg(u_{n_k}\big(H^{\omega}_{\Lambda_L,N} - E\big)^{-m-1}_{(\omega_{n_k}\to t\omega_{n_k})}\bigg)dt\bigg| \mathcal{D}_{k} \bigg)\nonumber\\ 
 &\qquad\qquad -m\mathbb{E}\bigg(\omega_{n_k}\int_0^1\operatorname{Tr}\!\bigg(u_{n_k}\big(H^{\omega}_{\Lambda_L, N} - E\big)^{-m-1}_{(\omega_{n_k}\to t\omega_{n_k})}\bigg)dt\bigg| \mathcal{D}_{k-1} \bigg)\nonumber\\ 
 & \quad= -m\mathbb{E}\bigg[\omega_{n_k}\int_0^1\bigg(\operatorname{Tr}\!\bigg(u_{n_k}\big(H^{\omega}_{\Lambda_L,D} - E\big)^{-m-1}_{(\omega_{n_k}\to t\omega_{n_k})}\bigg)\nonumber\\ 
 &\qquad\qquad\qquad\qquad\qquad\qquad-\bigg(u_{n_k}\big(H^{\omega}_{\Lambda_L,N} - E\big)^{-m-1}_{(\omega_{n_k}\to t\omega_{n_k})}\bigg)\bigg)dt \bigg| \mathcal{D}_k \bigg]\nonumber\\ & \qquad\qquad+m\mathbb{E}\bigg[\omega_{n_k}\int_0^1\bigg(\operatorname{Tr}\!\bigg(u_{n_k}\big(H^{\omega}_{\Lambda_L,D} - E\big)^{-m-1}_{(\omega_{n_k}\to t\omega_{n_k})}\bigg)\nonumber\\ &\qquad\qquad\qquad\qquad\qquad\qquad-\bigg(u_{n_k}\big(H^{\omega}_{\Lambda_L,N} - E\big)^{-m-1}_{(\omega_{n_k}\to t\omega_{n_k})}\bigg)\bigg) dt\bigg| \mathcal{D}_{k-1} \bigg]\nonumber\\ 
 & \quad= -m\mathbb{E}\bigg[\omega_{n_k}\int_0^1\bigg(\operatorname{Tr}\!\bigg(u_{n_k}\chi_{S_{n_k}}\big(H^{\omega}_{\Lambda_L,D} - E\big)^{-m-1}_{(\omega_{n_k}\to t\omega_{n_k})}\chi_{S_{n_k}}^*\bigg)\nonumber\\ 
 &\qquad\qquad\qquad\qquad\qquad\qquad-\bigg(u_{n_k}\chi_{S_{n_k}}\big(H^{\omega}_{\Lambda_L,N} - E\big)^{-m-1}_{(\omega_{n_k}\to t\omega_{n_k})}\chi_{S_{n_k}}^*\bigg)\bigg)dt \bigg| \mathcal{D}_k \bigg]\nonumber\\ 
 & \qquad\qquad+m\mathbb{E}\bigg[\omega_{n_k}\int_0^1\bigg(\operatorname{Tr}\!\bigg(u_{n_k}\chi_{S_{n_k}}\big(H^{\omega}_{\Lambda_L,D} - E\big)^{-m-1}_{(\omega_{n_k}\to t\omega_{n_k})}\chi_{S_{n_k}}^*\bigg)\nonumber\\ 
 &\qquad\qquad\qquad\qquad\qquad\qquad-\bigg(u_{n_k}\chi_{S_{n_k}}\big(H^{\omega}_{\Lambda_L,N} - E\big)^{-m-1}_{(\omega_{n_k}\to t\omega_{n_k})}\chi_{S_{n_k}}^*\bigg)\bigg) dt\bigg| \mathcal{D}_{k-1} \bigg]\nonumber. 
 \end{align}
In the last equality, we used Corollary~\ref{cr-tr-rst} together with the notation $S_{n_k} = \operatorname{supp}(u_{n_k})$, $u_{n_k} = u(x-n_k)$, and $S = \operatorname{supp}(u)$.  
Using (\ref{int-bl}), define 
\[
B_L^\circ = \{ n_k \in B_L : S_{n_k} \subset (\Lambda_L)_{\tilde{\ell}_L}^\circ \}.
\]  
Applying Proposition~\ref{df-n-d-in} to (\ref{est-df}), we obtain
\begin{equation}
\label{ex-dry}
\big|\mathbb{E}(\mathcal{K}_L \mid \mathcal{D}_k) - \mathbb{E}(\mathcal{K}_L \mid \mathcal{D}_{k-1})\big| \le 2 \|\omega_0\|_\infty m \|u\|_\infty e^{-\beta \tilde{\ell}_L / 2},~\text{for } k \in \{0,1,\ldots,\#B_L\}.
\end{equation}
We split the sum in (\ref{n-dmrt}) into two parts:
\begin{align}
\label{nd-dmrt}
\mathbb{E}\Big[\big(\mathcal{K}_L(\omega)\big)^2\Big]
&= \sum_{n_k \in B_L^\circ} \mathbb{E}\Big[\mathbb{E}(\mathcal{K}_L \mid \mathcal{D}_k) - \mathbb{E}(\mathcal{K}_L \mid \mathcal{D}_{k-1})\Big]^2 \nonumber\\
&\quad + \sum_{n_k \in B_L \setminus B_L^\circ} \mathbb{E}\Big[\mathbb{E}(\mathcal{K}_L \mid \mathcal{D}_k) - \mathbb{E}(\mathcal{K}_L \mid \mathcal{D}_{k-1})\Big]^2 \nonumber\\
&\le (\# B_L^\circ) \, 4 \|\omega_0\|_\infty^2 m^2 \|u\|_\infty^2 e^{-\beta \tilde{\ell}_L} + (\# B_L \setminus B_L^\circ) \, \big(2 \|\omega_0\|_\infty m \|u\|_\infty C_1 |S|\big)^2.
\end{align}
In the above, we used (\ref{ex-dry}) to estimate the first sum and applied Proposition~\ref{est-tr-chf} to estimate the second sum. Since we have $\displaystyle\lim_{L \to \infty} \frac{\# B_L^\circ}{|\Lambda_L|} = 1$ and $\displaystyle\lim_{L \to \infty} \frac{\# B_L \setminus B_L^\circ}{|\Lambda_L|} = 0$,  therefore (\ref{rsdf-nd}) follows from (\ref{nd-dmrt}). This completes the proof.
\end{proof}
\noindent Now, using the denseness of 
\(\mathcal{A}_{E,d}\big[-\|V\|_\infty,\infty\big)\) in 
\(C^1_{d,0}\big[-\|V\|_\infty,\infty\big)\), together with (\ref{vr-est}) and (\ref{vr-estn}), 
we can extend the convergence in (\ref{n=d}) to any test function 
\(f \in C^1_{d,0}\big[-\|V\|_\infty,\infty\big)\).
\begin{cor} 
\label{neu=dir}
Let $f \in C^1_{d,0}\big[-\|V\|_\infty,\infty\big)$, and consider the random variables 
$Y_{f,N,L}$ and $Y_{f,D,L}$ as in (\ref{rv}). 
Then, under Hypothesis \ref{hypo}, we have
\begin{equation}
\label{n=fd}
\frac{1}{|\Lambda_L|} \, \mathbb{E}\Bigg[\bigg(Y_{f,N,L} - Y_{f,D,L}\bigg)^2\Bigg] 
\xrightarrow{L \to \infty} 0.
\end{equation}
\end{cor}
\begin{proof}
Let \( Q_n \in \mathcal{A}_{E,d}\big[-\|V\|_\infty, \infty\big) \) be as in Remark~\ref{appr-fun-pol}. We can write the difference as
\begin{align}
Y_{f,N,L} - Y_{f,D,L}
&= Y_{f,N,L} - Y_{Q_n,N,L} + Y_{Q_n,D,L} - Y_{f,D,L} + Y_{Q_n,N,L} - Y_{Q_n,D,L} \nonumber\\
&= Y_{(f-Q_n),N,L} + Y_{(Q_n-f),D,L} + \big(Y_{Q_n,N,L} - Y_{Q_n,D,L}\big). \nonumber
\end{align}
Applying Minkowski's inequality, we obtain
\begin{align}
\label{n-d=0}
&\frac{1}{|\Lambda_L|^{\frac{1}{2}}}
\bigg(\mathbb{E}\Big[\big(Y_{f,N,L} - Y_{f,D,L}\big)^2\Big]\bigg)^{\frac{1}{2}} \nonumber\\
&\quad \leq
\frac{1}{|\Lambda_L|^{\frac{1}{2}}}
\bigg(\mathbb{E}\Big[\big(Y_{(f-Q_n),N,L}\big)^2\Big]\bigg)^{\frac{1}{2}}
+ \frac{1}{|\Lambda_L|^{\frac{1}{2}}}
\bigg(\mathbb{E}\Big[\big(Y_{(Q_n-f),D,L}\big)^2\Big]\bigg)^{\frac{1}{2}} \nonumber\\
&\qquad\qquad
+ \frac{1}{|\Lambda_L|^{\frac{1}{2}}}
\bigg(\mathbb{E}\Big[\big(Y_{Q_n,N,L} - Y_{Q_n,D,L}\big)^2\Big]\bigg)^{\frac{1}{2}}.
\end{align}
Define $\tilde{g}_n(x) = (x-E)^{1 + \lfloor d/2 \rfloor}\,\big(f'(x) - Q_n'(x)\big)
= \tilde{f}(x) - P_n(x)$.
Then, by Remark~\ref{appr-fun-pol}, we have \(\|\tilde{g}_n\|_\infty \to 0\) as \(n \to \infty\).
Now, taking \(\limsup\) as \(L \to \infty\) on both sides of~\eqref{n-d=0}, and using~\eqref{vr-est}, \eqref{vr-estn}, and~\eqref{n=d}, we obtain
\begin{align}
\limsup_{L\to\infty}
\frac{1}{|\Lambda_L|^{\frac{1}{2}}}
\bigg(\mathbb{E}\Big[\big(Y_{f,N,L} - Y_{f,D,L}\big)^2\Big]\bigg)^{\frac{1}{2}}
\leq 2\sqrt{C}\,\|\tilde{g}_n\|_\infty.
\end{align}
Since \(\|\tilde{g}_n\|_\infty \to 0\) as \(n \to \infty\) and the left-hand side of the above inequality is independent of \(n\), it follows that
\[
\limsup_{L\to\infty}
\frac{1}{|\Lambda_L|^{\frac{1}{2}}}
\bigg(\mathbb{E}\Big[\big(Y_{f,N,L} - Y_{f,D,L}\big)^2\Big]\bigg)^{\frac{1}{2}} = 0.
\]
Hence, \eqref{n=fd} holds.
\end{proof}
\noindent We are now ready to summarize the proof of our main result. \\~\\
\textbf{Proof of Theorem~\ref{mnthm}:}  
The convergence in~\eqref{clt-cn} was proved for the Dirichlet boundary condition \((X = D)\) in Lemma~\ref{nr-appr}.  
The corresponding convergence for the Neumann boundary condition \((X = N)\) follow from Corollary~\ref{neu=dir} and part~(i) of Proposition~\ref{cn-vr-rv}. \\
The finiteness and equality of the two variances,
$\sigma_f^2 := \sigma_{f,N}^2 = \sigma_{f,D}^2 < \infty$,
follows from Corollary~\ref{neu=dir} and part~(ii) of Proposition~\ref{cn-vr-rv}, together with~\eqref{vr-est}. \\
The expression for \(\sigma_{f,D}^2\) is given in Lemma~\ref{lmvrfrf}.  \\
The statement on the positivity of the limiting variance $\sigma^2_f$ is established in Lemma~\ref{ps-vr}.
Hence, the theorem is proved. \qed

\medskip
\medskip
\medskip

\noindent\textbf{Conflict of Interest}
The authors declare that they have no conflict of interest and no relevant financial or non-financial interests to disclose.

\medskip

\noindent\textbf{Data Availability}
No datasets were generated or analysed during the current study.

\appendix
\section{Appendix}

In this section, we prove several results concerning the calculus of the free 
magnetic Laplacian $(i\nabla + A)^2$ and its finite-volume restrictions. We also 
establish identities describing how the operators $H^{X}_{A,O} + V_O$, defined on 
$L^2(O)$, act on products $\varphi \psi$, provided the product belongs to the corresponding operator domain, for $O \subseteq \mathbb{R}^d$ and $X \in \{D, N\}$. In addition, we derive 
a formula for the derivative of
$\operatorname{Tr}\!\left( H^{X}_{A,O} + U + \lambda V - E \right)^{-m}$
with respect to $\lambda$. In the final part of the appendix, we 
collect several auxiliary results from probability theory.\\
All of the results results presented here are used in the proof of our central limit 
theorem. Although some of these results may be known (in other forms) and 
scattered throughout the literature, we present them here in the precise form 
required for our purposes.
\\~\\
Let's start by defining some notation. For any open set $O\subseteq \mathbb{R}^d$ we denote $C_c^\infty(O)$ to be the set of all infinitely differentiable compactly supported functions whose support is strictly contained in $O$.
Let $A:O\to\mathbb{R}^d$ be a real-valued vector potential such that $A\in \big(L_{\mathrm{loc}}^2(O)\big)^d$, in other words if $A(x)=\big(A_1(x), A_2(x),\ldots,A_d(x)\big)$ then $A_j\in L_{\mathrm{loc}}^2(O)~\forall~j$ and they are real-valued. For an open set $O\subseteq\mathbb{R}^d$, the inner product $\big\langle f, g\big\rangle_O$ on $\big(L^2(O)\big)^d$, $d\ge 1$, is defined by
$$\langle f, g\rangle_O:=\sum_{j=1}^d\int_O \bar{f}_j g_j dx,~f=(f_j)_{j=1}^d~\&~g=(g_j)_{j=1}^d\in\big(L^2(O)\big)^d.$$ 
On the finite-dimensional Hilbert space $\mathbb{C}^d$ the inner product is given by
$$z \cdot w:= \sum_{j=1}^d \bar{z}_jw_j,~z=(z_1,z_2,\ldots,z_d)~\&~w=(w_1,w_2,\ldots, w_d)\in \mathbb{C}^d.$$
\begin{defin}
The support of a function $f:O\subseteq\mathbb{R}^d\longrightarrow \mathbb{C}$ is defined by $\operatorname{supp}f:=\overline{\{x\in O: f(x)\neq 0  \}}$.
\end{defin}
\begin{defin}
\label{mg-sobo} 
Associated with the vector potential $A\in \big(L_{\mathrm{loc}}^2(O)\big)^d$, we define the magnetic Sobolev space $W_A^{1,2}(O)$ as
\begin{align}
\label{sobo}
W_A^{1,2}(O):=\big\{f\in L^2(O):(i\nabla+A)f\in \big(L^2(O)\big)^d  \big\},~~i=\sqrt{-1}.
\end{align}
Here, we denote $i\nabla+A$ to be the gauge-covariant gradient (with respect
to $A$) in the sense of distribution defined on the class of test functions $C_c^\infty(O)$. 
\end{defin}
\begin{rem}
\label{eq-sobo}
Let $A\in \big(L_{\mathrm{loc}}^2(O)\big)^d$; then using the definition of distributional (weak) derivative, we can rewrite magnetic Sobolev space $W_A^{1,2}(O)$ as
\begin{align}
\label{sobo1}
W_A^{1,2}(O):=&\bigcap_{k=1}^d\bigg\{f\in L^2(O):\text{there exists}~ f_k\in L^2(O) ~\text{such that}\nonumber \\
 &\qquad\qquad\big\langle f,i\partial_k \varphi+A_k\varphi\big\rangle_O=\big\langle f_k,\varphi \big\rangle_O\nonumber\\
 &\qquad\text{or}~\int_O\bar{f}(i\partial_k\varphi+A_k\varphi)dx=\int_O\bar{f_k}\varphi dx~\forall~\varphi\in C_c^\infty(O) \bigg\}.
\end{align}
Here $\partial_k \varphi$ denotes the partial derivative $\frac{\partial \varphi(x)}{\partial x_k},~x=(x_1,x_2,\ldots,x_d)$.
\end{rem}
\begin{nota}
\label{ex-g-gra}
For each $f\in W_A^{1,2}(O)$,  $O\subseteq\mathbb{R}^d$, the description in (\ref{sobo1}) of the distributional gauge-covariant gradient of $f$ gives 
$$(i\nabla+A)f=\big((i\nabla+A)_1f, (i\nabla+A)_2f,\ldots, (i\nabla+A)_df \big)\in\big(L^2(O)\big)^d.$$
In the above (\ref{sobo1}) we denote $f_k:=(i\nabla+A)_kf$ for $k=1,2,\ldots, d$. For $A\equiv0$,
one has $i\partial_k f=(i\nabla)_kf$ and $W_A^{1,2}(O)=W^{1,2}(O)$, the usual Sobolev space.
\end{nota}
\begin{nota}
For $\varphi\in C_c^\infty(O)$ and real-valued vector potential $A\in \big( L^2(O)\big)^d$, 
\big(here $O\subseteq \mathbb{R}^d$\big) we have
$(i\nabla+A)_k\varphi=i\partial_k\varphi+A_k\varphi$.
\end{nota}
\begin{nota}
\label{gaug-div}
Consider $W_A^{1,2}(O)$, $O\subseteq\mathbb{R}^d$ as in the Remark \ref{eq-sobo} and
a vector $F=(F_1,F_2,\ldots,F_d)\in \big( W_A^{1,2}(O) \big)^d $. Then, the gauge-covariant
divergence of $F$ is defined in the distributional sense by
$$(i\nabla+A)\cdot F=(i\nabla+A)_1F_1+ (i\nabla+A)_2F_2+\cdots+(i\nabla+A)_dF_d \in L^2(O). $$
Using Remark \ref{eq-sobo}, for each $k$ we have
\begin{equation}
 \label{in-by-p}
\big\langle F_k, (i\nabla+A)_k\varphi\big\rangle_O = \big\langle (i\nabla+A)_k F_k, \varphi\big\rangle_O ~~\forall~\varphi\in C_c^\infty(O).
\end{equation}
Subsequently, adding these identities over $k$, yields
\begin{equation}
\label{in-b-p-div}
\big \langle F,(i\nabla+A)\varphi \big\rangle_O=\big \langle(i\nabla+A)\cdot F,\varphi \big\rangle_O  ~~\forall~\varphi\in C_c^\infty(O).
\end{equation}
In the above,  the l.h.s. denotes the inner product on $\big(L^2(O) \big)^d$ and the r.h.s. inner product is on $L^2(O)$.
\end{nota}
\begin{nota}
\label{sc-vc}
For an open set $O\subseteq \mathbb{R}^d$, consider $f\in L^2(O)$ and  a real-valued vector potential $A=\big(A_1,A_2,\ldots,A_d  \big):O\to\mathbb{R}^d$ such that $A_kf\in L^2(O)$ for each $k$, then we define 
$$ Af:=\big(A_1 f, A_2f,\ldots, A_d f \big)\in \big(L^2(O) \big)^d.$$
\end{nota}
\noindent Let $O\subseteq \mathbb{R}^d$ be an open set. We define the sesquilinear form $q_{A,O}(\cdot,\cdot)$  on the magnetic Sobolev space $W_A^{1,2}(O)$ as
\begin{align}
\label{ses}
q_{A,O}(\varphi, \psi)=\big\langle (i\nabla+A)\varphi, (i\nabla+A)\psi\big\rangle_O,~~\varphi, \psi\in W_A^{1,2}(O).
\end{align}
Then, $W_A^{1,2}(O)$ is a Hilbert space equipped with the norm 
\begin{equation}
\label{sbnrm}
\|\psi\|^2_{A,O}=q_{A,O}(\psi,\psi)+\|\psi\|^2_O,~~\|\psi\|^2_O=\langle \psi, \psi\rangle_O,~~\psi\in W_A^{1,2}(O) .
\end{equation} 
More details can be found in \cite[Sections 7.3 and 7.20]{LL} and \cite[Lemma A.1]{TLMS}.
\begin{defin}
\label{cc-cl}
Denote $W^{1,2}_{A, cl}(O)$ as the closure of $C_c^\infty(O)$ with respect to the norm $\|\cdot\|_{A,O}$ described above. Here, $O\subseteq \mathbb{R}^d$ is an open set.
\end{defin}
\noindent Therefore, we have the inclusion
\begin{equation}
\label{sb-inq}
C_c^\infty(O)\subset W_{A, cl}^{1,2}(O)\subseteq W_A^{1,2}(O)\subset L^2(O),~~~O\subseteq \mathbb{R}^d.
\end{equation}
\begin{rem}
\label{sb-full}
It is well known that when $O=\mathbb{R}^d$ the equality $W_{A, cl}^{1,2}(\mathbb{R}^d)= W_A^{1,2}(\mathbb{R}^d)$ holds, we refer to \cite{{bary}, {KT}} for more details. For a proper open subset $O\subset \mathbb{R}^d$, the equality generally does not hold,  see \cite{RS4}.
\end{rem}
\noindent We next define two sesquilinear (or quadratic) forms associated with the spaces $W_A^{1,2}(O)$ and $W_{A, cl}^{1,2}(O)$.
\begin{defin}
\label{neu}
Let $O\subseteq \mathbb{R}^d$ be an open set and $A\in \big(L_{\mathrm{loc}}^2(O) \big)^d$.
The sesquilinear form $h^N_{A, O}(\cdot,\cdot)$ whose form core is $Q(h^N_{A, O})=W_A^{1,2}(O)$,
is defined by
\begin{equation}
\label{ses-neu}
h^N_{A, O}(\varphi, \psi)=\big\langle (i\nabla+A)\varphi, (i\nabla+A)\psi\big\rangle_{O},~~\varphi, \psi\in W_A^{1,2}(O)
\end{equation}
\end{defin}
\begin{defin}
\label{dir}
Let $O\subseteq \mathbb{R}^d$ be an open set and $A\in \big(L_{\mathrm{loc}}^2(O) \big)^d$.
The sesquilinear form $h^D_{A, O}(\cdot,\cdot)$ whose form core is
$Q(h^D_{A, O})=W_{A, cl}^{1,2}(O)$, is defined by
\begin{equation}
\label{ses-dir}
h^D_{A, O}(\varphi, \psi)=\big\langle (i\nabla+A)\varphi, (i\nabla+A)\psi\big\rangle_{O},~~\varphi, \psi\in W_{A, cl}^{1,2}(O).
\end{equation}
\end{defin}
\begin{rem}
Since $C_c^\infty(O)$ is dense in $L^2(O)$ with respect to the usual $\|\cdot\|_2$ norm on $L^2(O)$. It follows from $(\ref{sb-inq})$, that the sesquilinear forms $h^X_{A, O}(\cdot, \cdot),~X=D, N$ are densely defined on $L^2(O)$.
\end{rem}
\noindent The Neumann and Dirichlet magnetic Laplacians are defined via their closed sesquilinear
forms.
\begin{lem}
\label{dir-neu-lp}
Let $O\subseteq \mathbb{R}^d$ be open and the forms $h^N_{A, O}$ and $h^D_{A, O}$ as in (\ref{ses-neu}) and (\ref{ses-dir}).
\begin{enumerate}
\item[(i)] Both sesquilinear forms are positive, symmetric, and closed. Hence they define unique positive self-adjoint operators 
$H^N_{A,O}$ and $H^D_{A,O}$.
\item[(ii)] The restricted Neumann magnetic Laplacian $H^N_{A, O}$ is given by
 \begin{equation}
\label{neu-lapla}
h^N_{A,O}(\psi, u)=\big\langle \psi, H^N_{A, O}u\big\rangle,~\forall~\psi\in W^{1,2}_A(O),~u\in 
D(H^N_{A,O}).
\end{equation}
with domain
\[
D(H^N_{A,O})
=
\{u\in W^{1,2}_A(O):
\exists\,\tilde u\in L^2(O)
\text{ such that }
h^N_{A,O}(\psi,u)
=
\langle \psi,\tilde u\rangle
\ \forall \psi\in W^{1,2}_A(O)\},
\]
and for $u\in D(H^N_{A, O})$, we write $H^N_{A, O}u=\tilde{u}$.
\item[(iii)] The restricted Dirichlet magnetic Laplacian $H^D_{A, O}$ is given by
 \begin{equation}
\label{dir-lapla}
h^D_{A,O}(\psi, u)=\big\langle \psi, H^D_{A, O}u\big\rangle,~\forall~\psi\in W^{1,2}_{A,cl}(O),~u\in 
D(H^D_{A,O}).
\end{equation}
with domain
\[
D(H^D_{A,O})
=
\{u\in W^{1,2}_{A,cl}(O):
\exists\,\tilde u\in L^2(O)
\text{ such that }
h^D_{A,O}(\psi,u)
=
\langle \psi,\tilde u\rangle
\ \forall \psi \in W^{1,2}_{A,cl}(O)\}.
\]
and for $u\in D(H^D_{A, O})$, we write $H^D_{A, O}u=\tilde{u}$.
\end{enumerate}
\end{lem}
\begin{proof}
The proof can be found in \cite[Theorem VIII.15]{RS} (see also \cite{kato}).
\end{proof}

\noindent We now discuss the full magnetic Laplacian operator on $L^2(\mathbb{R}^d)$.
\begin{thm}
\label{fll-eq}
 Consider $W^{1,2}_{A,cl}(\mathbb{R}^d)$ and $W^{1,2}_{A}(\mathbb{R}^d)$ as in Definitions \ref{mg-sobo} and \ref{cc-cl}, then
 \begin{itemize}
 \item[(i)] We have $W^{1,2}_{A, cl}(\mathbb{R}^d)=
 W^{1,2}_{A}(\mathbb{R}^d)$.
 \item[(ii)] The forms $h^D_{A,\mathbb{R}^d}(\cdot, \cdot)$ and $h^N_{A,\mathbb{R}^d}(\cdot, \cdot)$, now defined
on the same domain, coincide.
 \item[(iii)] Consequently, $H^D_{A,\mathbb{R}^d}=H^N_{A,\mathbb{R}^d}$.
 \end{itemize}
 \end{thm}
\begin{proof}
The proof is given in \cite{bary} (see also \cite{KT}).
\end{proof}
\begin{rem}
\label{fl-d-lp}
We denote $H_A:=H^D_{A,\mathbb R^d}=H^N_{A,\mathbb R^d}=(i\nabla+A)^2
$, and its domain is
\begin{align}
D(H_A)&=\bigg\{\psi\in W^{1,2}_{A, cl}(\mathbb{R}^d)=
 W^{1,2}_{A}(\mathbb{R}^d): \text{there~exists}~\tilde{\psi}\in L^2(\mathbb{R}^d)~\text{such that}\nonumber\\
&\qquad \qquad\qquad \qquad h^D_{A,\mathbb{R}^d}(\varphi,\psi)=h^N_{A,\mathbb{R}^d}(\varphi,\psi)=\langle \varphi, \tilde{\psi}\rangle~
\nonumber\\
&\qquad \qquad\qquad \qquad \qquad \qquad \forall~\varphi\in W^{1,2}_{A, cl}(\mathbb{R}^d)=W^{1,2}_A(\mathbb{R}^d)  \bigg\},
\end{align}
and for $\psi\in D(H_A)$ we have $H_A\psi=\tilde{\psi}$.
\end{rem}
\begin{rem}
\label{fnv-dir-neu}
For a bounded open $O\subset \mathbb{R}^d$, the operators $H^N_{A,O}$, and $H^D_{A,O}$ are the finite-volume Neumann and Dirichlet magnetic Laplacian on $L^2(O)$, respectively. The operator $H_A$ is the free magnetic Laplacian on $L^2(\mathbb{R}^d)$.
\end{rem}

\noindent We consider the perturbations of magnetic Laplacians by a bounded operator $V$. Let $V\in L^\infty(\mathbb R^d)$ be real-valued and define the multiplication operator
 on $L^2(O)$
\begin{align}
\label{mult-v}
(Vf)(x)=V(x)f(x),~f\in L^2(O)~\text{and}~O\subseteq\mathbb{R}^d.
\end{align}
\begin{prop} 
\label{pertur}
Let $O\subseteq \mathbb{R}^d$ be an open set.  Consider $H^X_{A,O}$ \big($X=N, D$\big) and $H_A$ as defined in Lemma \ref{dir-neu-lp} and Remark \ref{fl-d-lp}. Define the perturbation operators as
\[
H^{X,V}_{A,O}=H^X_{A,O}+V
\quad\text{on }L^2(O),
\qquad
H_A^V=H_A+V
\quad\text{on }L^2(\mathbb R^d).
\]
Then $H^{X, V}_{A,O}$ and $H^V_A$ are self-adjoint with $D\big(H^{X, V}_{A,O} \big)=D\big(H^{X}_{A,O} \big)$ and 
$D\big(H^{V}_A\big)=D\big(H_A \big)$.
\end{prop}
\begin{proof}
Since $V$ is real-valued bounded function, it defines a bounded self-adjoint multiplication operator.
Thus $H^{X, V}_{A,O}$ and $H^V_A$ are self-adjoint perturbations of $H^{X}_{A,O}$
and $H_A$, respectively. Now, the proposition is immediate.
\end{proof}
\begin{rem}
\label{dn-pertur}
In the case $O=\mathbb R^d$, we have $H^{D, V}_{A,\mathbb{R}^d}=H^{N, V}_{A,\mathbb{R}^d}=H^V_A$ on $L^2(\mathbb{R}^d)$.
\end{rem}

\noindent The magnetic Schrödinger operators are semi-bounded, and their spectra lie in a half-line.
\begin{cor}
\label{spec}
Let $O\subset \mathbb{R}^d$ be open and $V\in L^\infty(\mathbb R^d)$ real--valued on $\mathbb{R}^d$, i.e., $\|V\|_\infty<\infty$.  $H^{X,V}_{A,O}$ ($X=D,N$) on $L^2(O)$ and $H_A^V$ on $L^2(\mathbb R^d)$ (as in Proposition \ref{pertur}) are semi--bounded  and self--adjoint.  In particular, we have
$$ \big\langle \varphi, H^{X,V}_{A,O}\varphi\big\rangle\ge -\|V\|_\infty \| \varphi \|_2^2~~~\forall~\varphi\in D\big(H^{X,V}_{A,O}  \big)\subset L^2(O) $$
and, 
$$ \big\langle \varphi, H^{V}_A\varphi\big\rangle\ge -\|V\|_\infty \| \varphi \|_2^2~~~\forall~\varphi\in D\big(H^V_A  \big)\subset L^2(\mathbb{R}^d) .$$
As a consequence, we obtain $\sigma\big(H^{X,V}_{A,O}  \big)\subseteq\big[-\|V\|_\infty, \infty\big)$ and $\sigma\big(H^{V}_A  \big)\subseteq\big[-\|V\|_\infty, \infty\big)$.
\end{cor}
\begin{proof}
As noted in \big(\ref{mult-v} \big), the multiplication operator $V$ satisfies
\begin{equation}
\label{semibd-V}
\big\langle \varphi, V\varphi \big\rangle\ge -\|V\|_\infty \| \varphi \|_2^2~~~\forall~~\varphi\in L^2(O).
\end{equation} 
Here $H^{X}_{A,O}$ and $H_A$ are unbounded positive self-adjoint operators. The claim follows because $H^{X,V}_{A,O}$ and $H_A^V$ are self-adjoint perturbation of 
$H^{X}_{A,O}$, $X=D, N$ and $H_A$ by bounded $V$.
\end{proof}

\begin{defin}\label{cn-ext-rst}
Let $O\subset \tilde O\subset\mathbb R^d$ be open.  
Define the restriction operator 
$\chi_O:L^2(\tilde O)\to L^2(O)$ by $(\chi_O f)(x)=f(x), x\in O.
$ Its adjoint $\chi_O^*:L^2(O)\to L^2(\tilde O)$ is the extension-by-zero operator,
\[
(\chi_O^* f)(x)=
\begin{cases}
f(x), & x\in O,\\
0, & x\in \tilde O\setminus O.
\end{cases}
\]
\end{defin}

\begin{nota}\label{cmp-spt}
For $u\in L^2(O)$, 
$\operatorname{supp}u\subset\subset O$ denotes $\operatorname{supp}u \subset O$ and compact.
\end{nota}

\noindent The gauge--covariant gradient preserves supports.
\begin{prop}
\label{sp-wd}
Let $A\in \big(L_{\mathrm{loc}}^2(O)  \big)^d$ be a real-valued vector potential, and $u\in W^{1,2}_A(O)$, where $O\subseteq \mathbb{R}^d$ is open.  Then, for each $k$, we have 
$\operatorname{supp}(i\nabla+A)_ku\subseteq \operatorname{supp}u$.
\end{prop}
\begin{proof}
Let $S := \operatorname{supp} u$.  
By definition of the distributional gauge-covariant derivative, for every
$\varphi \in C_c^\infty(O)$, $\int_O \overline{(i\nabla + A)_k u}\,\varphi \, dx
= \int_O \overline{u}\,(i\partial_k \varphi + A_k \varphi)\, dx.$
If $\varphi \in C_c^\infty(O \setminus S)$, then
$\operatorname{supp}(i\partial_k \varphi + A_k \varphi) \subset O \setminus S$,
and since $u$ vanishes on $O \setminus S$.
Hence,
\[
\int_{O \setminus S} \overline{(i\nabla + A)_k u}\,\varphi \, dx
= \int_{O \setminus S} \overline{u}\,(i\partial_k \varphi + A_k \varphi)\, dx=0
\quad \text{for all } \varphi \in C_c^\infty(O \setminus S).
\]
This implies $(i\nabla + A)_k u = 0$ a.e. on $O \setminus S$.
Therefore, $\operatorname{supp}(i\nabla + A)_k u \subseteq \operatorname{supp} u.
$
\end{proof}

\noindent It follows from the above that $\operatorname{supp}H^X_{A,O}u\subseteq \operatorname{supp}u$, for $X=D,N$.
\begin{cor}
\label{sp-neudir}
Let $H^X_{A,O},~X=D,N$ and $H_A$ be as in Lemma \ref{dir-neu-lp} and Remark \ref{fl-d-lp}. Then, under the same assumption as in Proposition \ref{sp-wd}, we have
$\operatorname{supp}H^X_{A,O}u\subseteq \operatorname{supp}u$ and $\operatorname{supp}H_A u\subseteq \operatorname{supp}u$
for all $u\in D\big(H^X_{A,O}  \big)$ and $u\in D\big(H_A  \big)$, respectively.
\end{cor}
\begin {proof}
For $u\in D\big(H^X_{A,O}  \big)$ we have $H^X_{A,O} u\in L^2(O)$,  where $O\subseteq \mathbb{R}^d$ .  Denote $S=\operatorname{supp}u$.  Now we have $C_c^\infty(O\setminus S)\subseteq C_c^\infty (O)\subset W^{1,2}_{A,cl}(O)\subseteq W^{1,2}_{A}(O)$. Therefore, for each $\varphi\in C^\infty_c\big( O\setminus S  \big)\subseteq C_c^\infty(O)$ 
$$\big\langle \varphi,  H^X_{A,O}u\big\rangle=\big\langle (i\nabla+A)\varphi, (i\nabla+A)u\big\rangle_O=0$$
In the above, we have used $\operatorname{supp}(i\nabla+A)u\subseteq \operatorname{supp}u$.  Since $C^\infty_c\big( O \setminus S \big)$ is dense in $L^2\big( O \setminus S \big)$ with respect to the usual norm $\|\cdot\|_2$, it follows that
$\operatorname{supp} H^X_{A,O} u \subseteq S = \operatorname{supp} u$.
In the same way, it is also true that
$\operatorname{supp}H_A u\subseteq \operatorname{supp}u$.
\end{proof}

\noindent The gauge–covariant gradient satisfies a Leibniz rule in the distributional sense. 
\begin{prop}
\label{prd-rul}
Let $O\subseteq \mathbb{R}^d$ be an open set, $A\in \big(L^2_{\mathrm{loc}}(O)\big)^d$ be a real-valued vector potential, $u\in W^{1,2}_A(O)$ and $\psi\in C^\infty(O)\cap L^\infty(O)$ with $\nabla \psi\in \big(  L^\infty(O)\big)^d$.
Then
\begin{equation}
\label{rul-gra}
(i\nabla+A)(\psi u)=\psi(i\nabla+A)u+(i\nabla \psi)u\in \big(L^2(O)\big)^d.
\end{equation}
\end{prop}
\begin{proof}
Let $\varphi \in C_c^\infty(O)$. For each $k=1,\dots,d$, using the definition of the gauge-covariant gradient in the sense of distributions and the product rule, we compute
\begin{align}
\label{leib}
\int_O \overline{\psi u}\,(i\partial_k\varphi + A_k\varphi)\,dx
&= \int_O \overline{u}\,\big(i\partial_k(\overline{\psi}\varphi) + A_k\overline{\psi}\varphi\big)\,dx
- \int_O \overline{u}\,(i\partial_k\overline{\psi})\,\varphi\,dx\nonumber\\
&=\int_O\overline{\psi((i\nabla+A)_ku)}~\varphi dx+\int_O\overline{(i\partial_k\psi)u}~\varphi dx.
\end{align}
Now using (\ref{leib}) and the Notation \ref{ex-g-gra}, the claim follows.
\end{proof}

\begin{rem}
\label{prd-cmp}
Under the assumptions of Proposition \ref{prd-rul}, for each $k=1,\dots,d$,
\begin{equation}
\label{cmp-prd}
(i\nabla+A)_k(\psi u)= \psi(i\nabla+A)_k u+(i\partial_k \psi)u\in L^2(O).
\end{equation}
\end{rem}

\begin{cor}
\label{prd-div}
Let for each $\psi_k$ and $u$ satisfy the assumptions of Proposition \ref{prd-rul}  and define $\Psi=\big( \psi_1,\psi_2,\ldots,\psi_d \big)$. Then
\begin{equation}
\label{rul-dv}
(i\nabla+A)\cdot (\Psi u)=\overline{\Psi}\cdot (i\nabla+A)u+(i\nabla\cdot \Psi)u\in L^2(O).
\end{equation}
\end{cor}
\begin{proof}
Apply Remark~\ref{prd-cmp} to each component $\psi_k u$ and sum over $k=1,\dots,d$.
\end{proof}

\noindent The gauge–covariant gradient is preserved under restriction.
\begin{prop}
\label{rst-gra}
Let $\tilde O\subset\mathbb R^d$ be open and 
$A\in(L^2_{\mathrm{loc}}(\tilde O))^d$ real-valued.  
If $u\in W^{1,2}_A(\tilde O)$ and $O\subset\tilde O$ is open, then we have $\chi_{O} u\in W^{1,2}_A(O)$ and $(i\nabla+A)(\chi_O u)=\chi_O(i\nabla+A)u$.
\end{prop} 
\begin{proof}
Let $\varphi\in C_c^\infty(O)$. 
Since $O\subset\tilde O$, we regard 
$\varphi\in C_c^\infty(\tilde O)$, and 
$\chi_O u=u$ on $\operatorname{supp}\varphi$.
For each $k=1,\dots,d$, we compute
\begin{align}
\label{drv}
\int_O \overline{\chi_O u}
       (i\partial_k\varphi + A_k\varphi)\,dx\nonumber
&= \int_{\tilde O}
      \overline{u}
      (i\partial_k\varphi + A_k\varphi)\,dx
 = \int_{\tilde O}
      \overline{((i\nabla+A)_k u)}\,\varphi\,dx \\
&= \int_O
      \overline{((i\nabla+A)_k u)}\,\varphi\,dx
 = \int_O
      \overline{\chi_O((i\nabla+A)_k u)}\,\varphi\,dx.
\end{align}
Thus $(i\nabla+A)(\chi_O u)
=\chi_O(i\nabla+A)u$ in distribution, and hence
$\chi_O u\in W^{1,2}_A(O)$.
\end{proof}

\noindent The above result also holds for $W^{1,2}_{A,cl}(\tilde{O})$, provided $\operatorname{supp}u\subset \subset O\subset \tilde{O}$.
\begin{prop}
\label{rst-qrd}
Let $\tilde{O}\subseteq \mathbb{R}^d$ be an open set and $A\in \big(L_{\mathrm{loc}}^2(\tilde{O}) \big)^d$ be a real-valued vector potential. If $u\in W^{1,2}_{A,cl}(\tilde{O})$, then for any open set $O\subset \tilde{O}$ with $\operatorname{supp}u\subset\subset O\subset \tilde{O}$ we have
$\chi_O u\in W^{1,2}_{A,cl}(O)$ and $(i\nabla+A)(\chi_O u)=\chi_O(i\nabla+A)u$.
\end{prop}
\begin{proof}
Since $\operatorname{supp}u\subset\subset O$, choose an open set 
$O_u$ such that $\operatorname{supp}u\subset O_u\subset O$,
and take $\psi\in C_c^\infty(O)$ with $\psi=1$ on $O_u$.  
Because $O\subset\tilde O$, we regard 
$\psi\in C_c^\infty(\tilde O)$ with 
$\operatorname{supp}\psi\subset O \subset\tilde O$.
Since $u\in W^{1,2}_{A,cl}(\tilde O)$, there exists 
$\{u_n\}_n\subset C_c^\infty(\tilde O)$ such that $\|u_n-u\|_{A,\tilde O}\to 0.
$
Since $u\in W^{1,2}_{A,cl}(\tilde{O})\subset  W^{1,2}_{A}(\tilde{O})$ and $\psi\in C_c^\infty(\tilde{O})$, by Proposition \ref{prd-rul}, we have $\psi u\in W^{1,2}_A(\tilde{O})$ and $\psi u_n\in C_c^\infty(\tilde{O})\subset W^{1,2}_{A}(\tilde{O})$. Because $\operatorname{supp}\psi\subset O$, $\chi_O(\psi u_n)\in C_c^{\infty}(O)\subset W^{1,2}_{A,cl}(O)\subset W^{1,2}_{A}(O)$, by Proposition \ref{rst-gra},
\[
(i\nabla+A)(\chi_O(\psi u_n))
=\chi_O(i\nabla+A)(\psi u_n)
=\chi_O\big[\psi(i\nabla+A)u_n+(i\nabla \psi)u_n  \big]
\]
Similarly,
\[
(i\nabla+A)(\chi_O(\psi u))
=\chi_O(i\nabla+A)(\psi u)
=\chi_O\big[\psi(i\nabla+A)u+(i\nabla \psi)u \big]
\]
Since $\|u_n-u\|_{A,\tilde O}\to 0$, it follows that $\|\chi_O(\psi u_n)-\chi_O(\psi u)\|_{A,O}\to 0.
$
Thus $\chi_O(\psi u)\in W^{1,2}_{A,cl}(O)$.
Since $\psi\equiv1$ on $O_u\supset\operatorname{supp}u$, we have $\chi_O(\psi u)=\chi_O u,
$
hence $\chi_O u\in W^{1,2}_{A,cl}(O) \subset W^{1,2}_A(O)$.
Finally, applying Proposition \ref{rst-gra}, gives $(i\nabla+A)(\chi_O u)
=\chi_O(i\nabla+A)u.$
\end{proof}

\noindent If $u$ belongs to the domain of a Dirichlet magnetic Laplacian, then for any smooth function $\psi$, the product $u\psi$ also lies in its domain.
\begin{lem}
\label{mult-dom}
Let $O\subseteq \mathbb{R}^d$ be an open set and $A\in \big(L^2_{\mathrm{loc}}(O)  \big)^d$ be a real-valued vector potential. Assume $u\in W^{1,2}_{A,cl}(O)$ and $\psi\in C^\infty(O)$ with $\psi,~\Delta\psi\in L^\infty(O)$ and $\nabla\psi\in \big(L^\infty(O)  \big)^d$.  
\begin{enumerate}
\item [(i)] Then the product $\psi u$ is also in the $W^{1,2}_{A,cl}(O)$.
\item[(ii)] If $u\in D \big( H^D_{A,O} \big)$, then
$\psi u \in D\big( H^D_{A,O} \big)$, and
\begin{equation}
\label{formula-prd}
H^D_{A,O}(\psi u)=\psi \big(H^D_{A,O}u  \big)+2\big(-i\nabla\bar{\psi}\big)\cdot (i\nabla+A)u- (\Delta \psi)u.
\end{equation}
\end{enumerate}
\end{lem}
\begin{proof}
(i) Since $u\in W^{1,2}_{A,cl}(O)$, there exists 
$\{u_n\}\subset C_c^\infty(O)$ such that 
$\|u_n-u\|_{A,O}\to 0$.  
Then $\psi u_n\in C_c^\infty(O)$, we have $(i\nabla+A)(\psi u_n)
=\psi(i\nabla+A)u_n+(i\nabla\psi)u_n,$
and similarly $(i\nabla+A)(\psi u)
=\psi(i\nabla+A)u+(i\nabla\psi)u.$
Because $\psi\in L^\infty(O)$ and 
$\nabla\psi\in (L^\infty(O))^d$, hence $\|\psi(u_n-u)\|_{A,O}\to 0$, so 
$\psi u\in W^{1,2}_{A,cl}(O)$.\\
(ii) Let $v\in C_c^\infty(O)$, then we write
\begin{align}
\label{mult-in-dom}
h^D_{A,O}(v,\psi u )&=\big\langle(i\nabla+A)v,  (i\nabla+A)(\psi u)\big\rangle_O\nonumber\\
&=\big\langle(i\nabla+A)v,  \psi(i\nabla+A) u\big\rangle_O+
\big\langle(i\nabla+A)v,  (i\nabla \psi)u\big\rangle_O\nonumber\\
&=\big\langle\overline{\psi}(i\nabla+A)v,  (i\nabla+A) u\big\rangle_O+\big\langle(i\nabla+A)v,  (i\nabla \psi)u\big\rangle_O\nonumber\\
&=\big\langle(i\nabla+A)(\overline{\psi}v),  (i\nabla+A) u\big\rangle_O
-\big\langle(i\nabla \overline{\psi})v,  (i\nabla+A) u\big\rangle_O\nonumber\\
&\qquad\qquad\qquad\qquad
\qquad\qquad~+\big\langle(i\nabla+A)v,  (i\nabla \psi)u\big\rangle_O\nonumber\\
&=h^D_{A,O}(\overline{\psi}v,  u)+\big\langle v,  (-i\nabla\overline{\psi})\cdot (i\nabla+A)u\big\rangle_O +\big\langle v,  (i\nabla+A)\cdot(i\nabla \psi)u\big\rangle_O\nonumber\\
&=\big\langle\overline{\psi}v, H^D_{A,O}u\big\rangle_O +\big\langle v,  (-i\nabla\overline{\psi})\cdot (i\nabla+A)u\big\rangle_O\nonumber\\
&\qquad\qquad +\big\langle v,  (-i\nabla\overline{\psi})\cdot (i\nabla+A)u\big\rangle_O
+\big\langle v,  (i\nabla\cdot i\nabla \psi)u\big\rangle_O\nonumber\\
&=\big\langle v,  \psi H^D_{A,O}u\big\rangle_O
+2\big\langle v,  (-i\nabla\overline{\psi})\cdot (i\nabla+A)u\big\rangle_O-\big\langle v, (\Delta \psi)u\big\rangle_O.
\end{align}
For the first part in the fifth line of above, we used the fact $\overline{\psi}v\in C_c^\infty(O)$ and $u\in D \big( H^D_{A,O} \big)$, together with identity (\ref{in-b-p-div}) for the third part.   In the sixth line, we have applied (\ref{rul-dv}).  Since $C_c^\infty(O)$ is dense in $W^{1,2}_{A,cl}(O)$ in the norm $\|\cdot\|_{A,O}$, therefore for each $v\in W^{1,2}_{A,cl}(O)$  using the above we get 
\begin{equation}
\label{qd-id}
h^D_{A,O}(v,\psi u )=\big\langle v,  \psi H^D_{A,O}u\big\rangle_O
+2\big\langle v,  (-i\nabla\overline{\psi})\cdot (i\nabla+A)u\big\rangle_O
-\big\langle v, (\Delta \psi)u\big\rangle_O.
\end{equation}
For fixed $\psi$ and $u$ satisfying the above assumption, it follows from (\ref{qd-id}) that $h^D_{A,O}(v,\psi u )$ defines a bounded linear functional of $v$ on $L^2(O)$, because $C_c^\infty(O)\subset W^{1,2}_{A, cl}(O)\subset L^2(O)$. So we get $\psi u\in D\big(H^D_{A,O}  \big)$ and $h^D_{A,O}(v,\psi u )=\big\langle v, H^D_{A,O}(\psi u)\big\rangle~\forall~v\in W^{1,2}_{A, cl}(O)$ and now (\ref{formula-prd}) is immediate from (\ref{qd-id}). Since $C_c^\infty(O)$ is dense in  $L^2(O)$ in the usual norm $\|\cdot\|_2:=\|\cdot\|_O$.  Hence the lemma.
\end{proof}

\begin{cor}
\label{cr-d-prt}
Let $H^{D,V}_{A,O}$ and $H^V_A$ be as in Proposition \ref{pertur} and assume $u\in D\big(H^{D,V}_{A,O} \big)$.  Then, under the same assumptions of Lemma \ref{mult-dom}, we obtain
\begin{equation}
\label{formula-prd-ptn}
H^{D,V}_{A,O}(u\psi)=\psi \big(H^{D,V}_{A,O}u  \big)+2\big(-i\nabla\bar{\psi}\big)\cdot (i\nabla+A)u- (\Delta \psi)u.
\end{equation}
When $O=\mathbb{R}^d$, we denote $H^{D,V}_{A,\mathbb{R}^d}:=H^V_A=H_A+V,~H_A=H^{D}_{A,\mathbb{R}^d}$.
\end{cor}
\begin{proof}
Since $H^{D,V}_{A,O}=H^D_{A,O}+V$ on $L^2(O)$ for an open set $O\subseteq\mathbb{R}^d$ and we have 
$D\big(H^{D,V}_{A,O} \big)=D\big(H^{D}_{A,O} \big)$, therefore (\ref{formula-prd-ptn}) follows from (\ref{formula-prd}).
\end{proof}

\begin{rem}\label{sp-prd}
Identity~\eqref{formula-prd-ptn} implies that $\operatorname{supp}\big(H^{D,V}_{A,O}(u\psi)\big)
\subset \operatorname{supp}\psi$
\end{rem}

\begin{cor}
\label{inv-mupt-op}
Let $O\subseteq\mathbb{R}^d$ be an open set and let $V\in L^\infty(O)$ be real-valued. Then, for each $\lambda<-\|V\|_\infty$, under the same assumptions of Lemma \ref{mult-dom} we have 
\begin{equation}
\label{inv-mult}
\big\|\big(H^{D,V}_{A,O}-\lambda\big)\psi\big(H^{D,V}_{A,O}-\lambda\big)^{-1}g\|_O\leq C^V_{\lambda,\psi}\|g\|_O~~\forall~g\in L^2(O).
\end{equation}
The positive constant $C^V_{\lambda,\psi}$ depends on $\lambda$, $V$, and $\psi$ but independent of $A$, $O$, and $g$.
\end{cor}
\begin{proof}
Let $u=\big(H^{D,V}_{A,O}-\lambda\big)^{-1}g$. Then $u\in D\big(H^{D,V}_{A,O} \big)$ and $H^{D,V}_{A,O}u=g+\lambda u$ and using (\ref{formula-prd-ptn}), we get
\begin{align}
\label{inv-bd}
\big(H^{D,V}_{A,O}-\lambda\big)\psi\big(H^{D,V}_{A,O}-\lambda\big)^{-1}g&=\psi g
+2\big(-i\nabla\bar{\psi}\big)\cdot (i\nabla+A)\bigg(\big(H^{D,V}_{A,O}-\lambda\big)^{-1}g\bigg)\nonumber\\
&\qquad\qquad\qquad- (\Delta \psi)\big(H^{D,V}_{A,O}-\lambda\big)^{-1}g.
\end{align}
Since $\| \psi g\|_O\leq \|\psi \|_\infty\|g\|_O$,  $\| u\|_O=\big\| \big(H^{D,V}_{A,O}-\lambda\big)^{-1}g \big\|_O\leq C_{V,\lambda}\|g\|_O$ and the derivatives (up to second order) of $\psi$ are bounded, it suffices to prove that:
\begin{equation}
\label{gra-inv-est}
\bigg\| (i\nabla+A)\bigg(\big(H^{D,V}_{A,O}-\lambda\big)^{-1}g\bigg) \bigg\|_{\big(L^2(O)\big)^d}\leq C^V_\lambda \|g\|_{O}.
\end{equation}
Since $H^{D,V}_{A,O}u=g+\lambda u$ i.e. $H^{D}_{A,O}u=g+(\lambda-V)u$, the definition of the sesquilinear form $h^D_{A,O}(\cdot, \cdot)$ yields
\begin{align}
\big\|(i\nabla+A)u  \big\|_O^2&=h^D_{A,O}(u, u)=\big\langle u,  H^{D}_{A,O}u\big\rangle_O =\langle u, g\rangle_O +\big\langle u, (\lambda-V)u\big\rangle_O.
\end{align}
Since $\| u\|_O\leq C_{V,\lambda}\|g\|_O$ and $\| V\|_\infty<\infty$, inequality (\ref{gra-inv-est}) is follows from the above. Hence, (\ref{inv-mult}) follows.
\end{proof}

\noindent An analogue of Lemma \ref{mult-dom} also holds for the Neumann magnetic Laplacian, and to prove it, we require an integration-by-parts formula for the magnetic Laplacian for weak derivative.
\begin{prop}
\label{bpt-wk}
Let $\tilde{O}\subseteq \mathbb{R}^d$ be an open set and let $A\in \big(L^2_{\mathrm{loc}}(\tilde{O})\big)^d$ be a real-valued vector potential. Assume $v\in W^{1,2}_A(\tilde{O})$ with $Av\in\big(L^2_{\mathrm{loc}}(\tilde{O})\big)^d$.
If $\psi=\big(\psi_1,\psi_2,\ldots,\psi_d\big)\in \big(W^{1,2}_A(\tilde{O})\big)^d$ satisfies $\operatorname{supp}\psi_k\subset \subset \tilde{O}$ for each $k$, and $A\psi\in\big(L^2_{\mathrm{loc}}(\tilde{O})  \big)^d$, then
\begin{equation}
\label{i-bpts}
\big\langle (i\nabla+A)v, \psi\big\rangle_{\tilde{O}}=\big\langle v,  (i\nabla+A)\cdot \psi\big\rangle_{\tilde{O}}.
\end{equation}
\end{prop}
\begin{proof}
Since $\operatorname{supp}\psi_k\subset\subset\tilde O$, 
choose a bounded open set $O\subset\tilde O$ with smooth boundary such that 
$\operatorname{supp}\psi_k\subset O$ for all $k$.
Using Propositions \ref{rst-gra} and \ref{wk-dr-nn-mg}, 
\begin{align*}
\langle (i\nabla+A)v,\psi\rangle_{\tilde O}
&=\langle \chi_O(i\nabla+A)v,\chi_O\psi\rangle_O =\langle (i\nabla+A)(\chi_O v),\chi_O\psi\rangle_O \\
&=\langle i\nabla(\chi_O v),\chi_O\psi\rangle_O
 +\langle A(\chi_O v),\chi_O\psi\rangle_O \\
&=\langle \chi_O v, i\nabla\!\cdot(\chi_O\psi)\rangle_O
 +\langle \chi_O v, A\!\cdot(\chi_O\psi)\rangle_O \\
&=\langle \chi_O v,(i\nabla+A)\!\cdot(\chi_O\psi)\rangle_O =\langle \chi_O v,\chi_O((i\nabla+A)\!\cdot\psi)\rangle_O \\
&=\langle v,(i\nabla+A)\!\cdot\psi\rangle_{\tilde O}.
\end{align*}
In the third line, we have used the Gauss identity \cite[Exercise II.4.3]{gal}, \big(an application of Gauss divergence theorem on $W^{1,2}(O)$\big) together with the fact that $\psi=0$ on $\partial O$. In the last line we used $\operatorname{supp}\psi_k\subset O\subset\tilde O$.
\end{proof}

\noindent Using (\ref{i-bpts}), we obtain a product rule for the Neumann magnetic Laplacian.
\begin{lem}
\label{mult-dom-neu}
Let $O\subseteq \mathbb{R}^d$ be an open set and let $A\in \big(L^\infty_{\mathrm{loc}}(O)  \big)^d$ be a real-valued vector potential.  Assume $u\in W^{1,2}_{A}(O)$ and $\psi\in C^\infty(O)\cap L^\infty(O)$ with $\nabla\psi\in \big(C_c(O)\big)^d$.
\begin{enumerate}
\item [(i)] Then the product $\psi u$ is also in the $W^{1,2}_{A}(O)$.
\item[(ii)] If $u\in D \big( H^N_{A,O} \big)$, then $\psi u \in D\big( H^N_{A,O} \big)$ and
\begin{equation}
\label{formula-prd-neu}
H^N_{A,O}(\psi u)=\psi \big(H^N_{A,O}u  \big)+2\big(-i\nabla\bar{\psi}\big)\cdot (i\nabla+A)u- (\Delta \psi)u.
\end{equation}
\end{enumerate}
When $O=\mathbb{R}^d$, we denote $H^N_{A,\mathbb{R}^d}=H_A$.
\end{lem}
\begin{proof}
(i) This follows from Proposition~\ref{prd-rul}.\\
\medskip
(ii) Let $v\in W^{1,2}_A(O)$. Arguing as in Lemma~\ref{mult-dom}(ii), we compute
\begin{align}
h^N_{A,O}(v,\psi u )&=\big\langle(i\nabla+A)(\overline{\psi}v), (i\nabla+A) u\big\rangle_O -\big\langle(i\nabla \overline{\psi})v, (i\nabla+A) u\big\rangle_O\nonumber\\ &\qquad\qquad\qquad\qquad \qquad\qquad~+\big\langle(i\nabla+A)v, (i\nabla \psi)u\big\rangle_O\nonumber\\
&=h^N_{A,O}(\overline{\psi}v, u)+\big\langle v, (-i\nabla\overline{\psi})\cdot (i\nabla+A)u\big\rangle_O+\big\langle v, (i\nabla+A)\cdot(i\nabla \psi)u\big\rangle_O\nonumber\\ 
&=\big\langle v, \psi H^N_{A,O}u\big\rangle_O +2\big\langle v, (-i\nabla\overline{\psi})\cdot (i\nabla+A)u\big\rangle_O-\big\langle v, (\Delta \psi)u\big\rangle_O. \nonumber
\end{align}
In the second line we used that $\overline{\psi}v\in W^{1,2}_A(O)$ and $u\in D(H^N_{A,O})$, together with Proposition~\ref{bpt-wk}, noting that each component of $(i\nabla\psi)u$ has compact support contained in $O$. The remaining steps follow exactly as in Lemma~\ref{mult-dom}(ii). Hence for each $v\in W^{1,2}_{A}(O)$, we  get 
\begin{equation}
\label{qd-id-neu}
h^N_{A,O}(v,\psi u )=\big\langle v,  \psi H^N_{A,O}u\big\rangle_O
+2\big\langle v,  (-i\nabla\overline{\psi})\cdot (i\nabla+A)u\big\rangle_O
-\big\langle v, (\Delta \psi)u\big\rangle_O.
\end{equation}
For fixed $\psi$ and $u$ satisfying the above assumptions, it follows from (\ref{qd-id-neu}) that $h^N_{A,O}(v,\psi u )$ extends to a bounded linear functional (of $v$) on $L^2(O)$ because $C_c^\infty(O)\subset W^{1,2}_{A}(O)\subset L^2(O)$. So we get $\psi u\in D\big(H^N_{A,O}  \big)$ and $h^N_{A,O}(v,\psi u )=\big\langle v, H^N_{A,O}(\psi u)\big\rangle~\forall~v\in W^{1,2}_{A}(O)$ and now (\ref{formula-prd-neu}) is immediate from (\ref{qd-id-neu}) as $C_c^\infty(O)$ is dense in the $L^2(O)$ in the usual norm $\|\cdot\|_2:=\|\cdot\|_O$.  Hence the lemma.
\end{proof}

\begin{cor}
\label{cr-d-prt-neu}
Consider $H^{N,V}_{A,O}$ and $H^V_A$ be as in Proposition \ref{pertur}, and assume $u\in D\big(H^{N,V}_{A,O} \big)$.  Then, under the assumptions of Lemma \ref{mult-dom-neu}, we obtain
\begin{equation}
\label{formula-prd-ptn-nu}
H^{N,V}_{A,O}(u\psi)=\psi \big(H^{N,V}_{A,O}u  \big)+2\big(-i\nabla\bar{\psi}\big)\cdot (i\nabla+A)u- (\Delta \psi)u.
\end{equation}
In this case, when $O=\mathbb{R}^d$, we denote $H^{N,V}_{A,\mathbb{R}^d}:=H^V_A=H_A+V,~H_A=H^{N}_{A,\mathbb{R}^d}$.
\end{cor}
\begin{proof}
The proof is similar to that of Corollary \ref{cr-d-prt}.
\end{proof}
\begin{rem}
\label{sp-prd-neu}
Identity (\ref{formula-prd-ptn-nu}) gives, $\operatorname{supp}(H^{N,V}_{A,O}(u\psi))\subset \operatorname{supp}\psi$.
\end{rem}
\begin{cor}
\label{inv-mupt-op-neu}
Let $O\subseteq\mathbb{R}^d$ be an open set and let $V\in L^\infty(O)$ be real-valued function . Then, for each $\lambda<-\|V\|_\infty$ under the same assumptions of Lemma \ref{mult-dom-neu} we  have 
\begin{equation}
\label{inv-mult-neu}
\big\|\big(H^{N,V}_{A,O}-\lambda\big)\psi\big(H^{N,V}_{A,O}-\lambda\big)^{-1}g\|_O\leq C^V_{\lambda,\psi}\|g\|_O~~\forall~g\in L^2(O).
\end{equation}
The positive constant $C^V_{\lambda,\psi}$ depends on $\lambda$, $V$, and $\psi$, but independent of $A$, $O$, and $g$.
\end{cor}
\begin{proof}
The proof is similar to that of Corollary \ref{inv-mupt-op}.
\end{proof}

\noindent Let $O\subset \tilde{O}\subseteq \mathbb{R}^d$ are open sets. Then we can canonically embed 
$W^{1,2}_{A,cl}(O)$ into $W^{1,2}_{A,cl}(\tilde{O})$.
\begin{prop}
\label{cn-embd}
Let $\tilde{O}\subseteq\mathbb{R}^d$ be open and $A\in \big(L_{\mathrm{loc}}^2(\tilde{O})\big)^d $ real-valued. If $O\subset\tilde O$ is open and $u\in W^{1,2}_{A,cl}(O)$, then $\chi_O^* u\in W^{1,2}_{A,cl}(\tilde{O})$ and $(i\nabla+A)(\chi_O^* u)=\chi_O^*(i\nabla+A)u\in \big( L^2(\tilde{O}) \big)^d$,  where $\chi_O^*:L^2(O)\to L^2(\tilde{O})$ as in the Definition \ref{cn-ext-rst}.
\end{prop}
\begin{proof}
Since $u\in W^{1,2}_{A,cl}(O)$, there exists 
$\{u_n\}_n\subset C_c^\infty(O)$ such that $\|u_n-u\|_{A,O}\to 0$. Because $O\subset\tilde O$, we have $\chi_O^*u_n\in C_c^\infty(\tilde{O})$ and $(i\nabla+A)(\chi_O^*u_n)=\chi_O^*(i\nabla+A)u_n$. Moreover, $\|\chi_O^*u_n-\chi_O^*u_m\|_{A,\tilde O}
=\|u_n-u_m\|_{A,O}$,
so $\{\chi_O^*u_n\}_n$ is Cauchy in $W^{1,2}_A(\tilde O)$.
Since $(i\nabla+A)(\chi_O^*u_n)=\chi_O^*(i\nabla+A)u_n\to \chi_O^*(i\nabla+A)u$ and $\chi_O^*u_n\to \chi_O^*u$ as $n\to\infty$ in $L^2(\tilde O)$, it follows that $\chi_O^* u\in W^{1,2}_{A}(\tilde{O})$ and $(i\nabla+A)(\chi_O^* u)=\chi_O^*(i\nabla+A)u$. Because $\chi_O^*u_n\in C_c^\infty(\tilde O)$ and 
$\chi_O^*u_n\to\chi_O^*u$ in $\|\cdot\|_{A,\tilde O}$,
we conclude $\chi_O^*u\in W^{1,2}_{A,cl}(\tilde O)$.
\end{proof}

\noindent Now we can derive a relation, in a suitable sense of restriction, between the two Dirichlet magnetic Laplacians $H^D_{A,O}$ and $H^D_{A,\tilde{O}}$, where $O\subset \tilde{O}\subseteq \mathbb{R}^d$ are open sets. 
\begin{lem}
\label{rest-prt-dir}
Let $u\in D \big( H^D_{A,\tilde{O}}\big)$ with $\operatorname{supp}u\subset\subset O$ where $O\subset\tilde{O}\subseteq\mathbb{R}^d$ are open. Assume $A\in \big(L_{\mathrm{loc}}^2(\tilde{O})\big)^d $ is a real-valued vector potential, then $\chi_O u\in D \big( H^D_{A,O}\big)$ and $H^D_{A,O}(\chi_O u)=\chi_O\big(H^D_{A,\tilde{O}} u\big)$,  where $\chi_O:L^2(\tilde{O})\to L^2(O)$ as in the Definition \ref{cn-ext-rst}.
\end{lem}
\begin{proof}
Since $u\in D(H^D_{A,\tilde O})\subset W^{1,2}_{A,cl}(\tilde O)$ and 
$\operatorname{supp}u\subset\subset O\subset\tilde O$, 
Proposition \ref{rst-qrd} yields $\chi_O u\in W^{1,2}_{A,cl}(O)$ and $(i\nabla+A)(\chi_Ou)=\chi_O(i\nabla+A)u$. Let $v\in W^{1,2}_{A,cl}(O)$. 
By Proposition \ref{cn-embd}, 
$\chi_O^*v\in W^{1,2}_{A,cl}(\tilde O)$.
We compute
\begin{align}
\label{to-dom}
h^D_{A,O}(v,\chi_O u)
&=\langle (i\nabla+A)v,
          (i\nabla+A)(\chi_O u)\rangle_O =\langle (i\nabla+A)v,
          \chi_O(i\nabla+A)u\rangle_O \nonumber\\
&=\langle \chi_O^*(i\nabla+A)v,
          (i\nabla+A)u\rangle_{\tilde O} =\langle (i\nabla+A)(\chi_O^*v),
          (i\nabla+A)u\rangle_{\tilde O}\nonumber \\
&=h^D_{A,\tilde O}(\chi_O^*v,u) =\langle \chi_O^*v,
          H^D_{A,\tilde O}u\rangle_{\tilde O} =\langle v,
          \chi_O H^D_{A,\tilde O}u\rangle_O.
\end{align}
Since $C_c^\infty(O)\subset W^{1,2}_{A,cl}(O)\subset L^2(O)$ and $C_c^\infty(O)$ is dense in $L^2(O)$ in $\|\cdot \|_2=\| \cdot\|_O$, identity (\ref{to-dom}) shows that
$h^D_{A,O}(v,  \chi_Ou)$ defines a bounded linear functional of $v$ on $L^2(O)$ and $\chi_O u \in D\big(H^D_{A, O}  \big)$ with 
$h^D_{A, O}(v,  \chi_O  u)=\big\langle v,  H^D_{A, O}(\chi_O u)\big\rangle_O$. Finally, (\ref{to-dom}) implies 
$$\big\langle v,  H^D_{A, O}(\chi_O u)\big\rangle_O=\big\langle v,  \chi_O H^D_{A,\tilde{O}} u\big\rangle_{O}~~\forall~v\in W^{1,2}_{A,cl}(O).$$
Since $C_c^\infty(O)\subset W^{1,2}_{A,cl}(O)$, we conclude $H^D_{A, O}(\chi_O  u)=\chi_O\big( H^D_{A,\tilde{O}}u\big)$ for $u\in D\big(H^D_{A, \tilde{O}} \big)$.
\end{proof}

\noindent For $u\in W^{1,2}_A(O)$ with $\operatorname{supp}u\subset \subset O\subset\tilde{O} \subset\mathbb{R}^d$ can be canonically embed $u$ into $W^{1,2}_A(\tilde{O})$.
\begin{lem}
\label{ext-0}
Let $\tilde{O}\subseteq \mathbb{R}^d$ be open and $A\in \big(L_{\mathrm{loc}}^2(\tilde{O})\big)^d$ real-valued. If $u\in W^{1,2}_A(O)$ with $\operatorname{supp} u\subset \subset O \subset \tilde{O}$ for an open set $O$, then we have $\chi_{O}^*u\in W^{1,2}_A(\tilde{O})$ and $(i\nabla+A)\big(\chi_O^* u\big)=\chi_O^*\big((i \nabla+A) u\big)\in \big(L^2(\tilde{O})\big)^d$. 
\end{lem}
\begin{proof}
Since $\operatorname{supp}u\subset\subset O$, choose an open set 
$O_u$ such that $\operatorname{supp}u\subset O_u\subset O$
and take $\psi\in C_c^\infty(O)$ with $\psi=1$ on $O_u$.
Define
$\widetilde{u\psi}\in L^2(\tilde{O})$ as
\[
\widetilde{u\psi}(x)=
\begin{cases}
u(x)\psi(x), & x\in O,\\
0, & x\in \tilde O\setminus O.
\end{cases}
\]
Let $\varphi\in C_c^\infty(\tilde O)$. For each $k=1,\dots,d$,
\begin{align*}
\big\langle \widetilde{u\psi}, (i\partial_k+A_k)\varphi \big\rangle_{\tilde O}
&= \int_{\tilde O}\overline{\widetilde{u\psi}}
   (i\partial_k\varphi+A_k\varphi)\,dx = \int_O \overline{u\psi}
   (i\partial_k\varphi+A_k\varphi)\,dx \\
&= \int_{(\operatorname{supp}u)} \overline{u}
   \big(i\partial_k(\bar{\psi}\varphi)
   - i(\partial_k\bar{\psi})\varphi
   + A_k(\bar{\psi}\varphi)\big)\,dx.
\end{align*}
Since $\partial_k\psi=0$ on $\operatorname{supp}u \subset \subset O$, the middle term vanishes, hence
\begin{align*}
&= \int_O \overline{u}
   \big(i\partial_k(\bar{\psi}\varphi)
   + A_k(\bar{\psi}\varphi)\big)\,dx.
\end{align*}
Because $\overline{\psi}\varphi\in C_c^\infty(O)$ and 
$u\in W^{1,2}_A(O)$, using that 
$\operatorname{supp}((i\nabla+A)_k u)
\subset \operatorname{supp}u$, we obtain
\[
= \int_O
\overline{((i\nabla+A)_k u)}\,\bar{\psi}\varphi\,dx= \int_{(\operatorname{supp}u)}
\overline{((i\nabla+A)_k u)}\,\bar{\psi}\varphi\,dx
\]
Using that $\psi=1$ on $\operatorname{supp}u$, $\operatorname{supp}((i\nabla+A)_k u)
\subset \operatorname{supp}u\subset O_u\subset O\subset \tilde{O}$, we get
\[
= \int_O
\overline{((i\nabla+A)_k u)}\,\varphi\,dx
= \int_{\tilde O}
\overline{\chi_O^*((i\nabla+A)_k u)}\,\varphi\,dx.
\]
Since $\widetilde{u\psi}=\chi_O^*u$ in $L^2(\tilde O)$, we obtain $(i\nabla+A)(\chi_O^*u)
=
\chi_O^*((i\nabla+A)u)$ in the distributional sense,
we conclude that $\chi_O^*u\in W^{1,2}_A(\tilde O)$.
\end{proof}

\noindent The Neumann magnetic Laplacian satisfies an analogous restriction property.
\begin{lem}
\label{rest-prt-dir-neu}
Let $\tilde{O}\subseteq\mathbb{R}^d$ be open and $A\in (L_{\mathrm{loc}}^2(\tilde{O})\big)^d $ real-valued. If $u\in D \big( H^N_{A,\tilde{O}}\big)$ with
$\operatorname{supp}u\subset\subset O$ for an open set $O\subset \tilde{O}$, then $\chi_O u\in D \big( H^N_{A,O}\big)$ and $H^N_{A,O}(\chi_O u)=\chi_O\big(H^N_{A,\tilde{O}} u\big)$, where $\chi_O:L^2(\tilde{O})\to L^2(O)$ as in the Definition \ref{cn-ext-rst}.
\end{lem}
\begin{proof}
Since $u\in D(H^N_{A,\tilde O})\subset W_A^{1,2}(\tilde O)$, 
Proposition \ref{rst-gra} yields, $\chi_O u \in W_A^{1,2}(O)$ and $(i\nabla+A)(\chi_O u)=\chi_O(i\nabla+A)u$. Because $\operatorname{supp}u\subset\subset O$, choose an open set 
$O_u$ such that $\operatorname{supp}u\subset O_u\subset O$ and take $\psi\in C_c^\infty(O)$ with $\psi=1$ on $O_u$.
Let $v\in W_A^{1,2}(O)$. Then
\begin{align*}
h^N_{A,O}(v,\chi_O u)
&=\langle (i\nabla+A)v, (i\nabla+A)(\chi_O u)\rangle_O=\langle (i\nabla+A)v, \chi_O(i\nabla+A)u\rangle_O
\end{align*}
since $\operatorname{supp}(i\nabla+A)u
\subseteq \operatorname{supp}u\subset O_u \subset O$, and $\psi\equiv 1$ on $O_u$, we get
\begin{align*}
&=\langle \chi_{O_u}(i\nabla+A)v,
          \chi_{O_u}(i\nabla+A)u\rangle_{O_u}=\langle \chi_{O_u}\psi(i\nabla+A)v,
          \chi_{O_u}(i\nabla+A)u\rangle_{O_u}.
\end{align*}
Using the product rule (Proposition \ref{prd-rul}),
\begin{align*}
&=\langle \chi_{O_u}(i\nabla+A)(\psi v),
          \chi_{O_u}(i\nabla+A)u\rangle_{O_u}
   - \langle \chi_{O_u}(i\nabla\psi)v,
          \chi_{O_u}(i\nabla+A)u\rangle_{O_u}.
\end{align*}
Since $\nabla\psi=0$ on $O_u$, the second term vanishes, and again using $\operatorname{supp}(i\nabla+A)u
\subseteq \operatorname{supp}u\subset O_u \subset O$, we get
\begin{align*}
&=\langle (i\nabla+A)(\psi v),
          \chi_O(i\nabla+A)u\rangle_O =\langle \chi_O^*(i\nabla+A)(\psi v),
          (i\nabla+A)u\rangle_{\tilde O}.
\end{align*}
Because $\operatorname{supp}(\psi v)\subset\subset O$,
Proposition \ref{ext-0} gives
\begin{align*}
&=\langle (i\nabla+A)(\chi_O^*(\psi v)),
          (i\nabla+A)u\rangle_{\tilde O} = h^N_{A,\tilde O}(\chi_O^*(\psi v),u) \\
&= \langle \chi_O^*(\psi v),
          H^N_{A,\tilde O}u\rangle_{\tilde O} = \langle \psi v,
          \chi_O H^N_{A,\tilde O}u\rangle_O.
\end{align*}
Finally, since 
$\operatorname{supp}(H^N_{A,\tilde O}u)
\subseteq \operatorname{supp}u\subset O_u$
and $\psi\equiv 1$ on $O_u$, we obtain
\begin{equation}
\label{cnt-lf}
h^N_{A,O}(v,\chi_O u)
=
\langle v,
        \chi_O H^N_{A,\tilde O}u\rangle_O.
\end{equation}
By the same density argument used at the end of 
Lemma \ref{mult-dom-neu}, we conclude that $\chi_O u \in D(H^N_{A,O})$ and $H^N_{A,O}(\chi_O u)
=
\chi_O(H^N_{A,\tilde O}u).$
\end{proof}

\begin{cor}
\label{res-pertur-neu-dir}
Let $O\subset \tilde{O}\subseteq \mathbb{R}^d$ be open sets and consider $H^{X,V}_{A,\tilde{O}}$,  $X=D,N$ as in Proposition \ref{pertur}. If  $u\in D\big(H^{X,V}_{A,\tilde{O}} \big)$ with $\operatorname{supp}u\subset \subset O \subset \tilde{O}$.  Then, under the same assumptions of Lemmas \ref{rest-prt-dir} and \ref{rest-prt-dir-neu}, we get $\chi_O u\in  D\big(H^{X,V}_{A,O} \big)$ and $H^{X,V}_{A,O}(\chi_O u)=\chi_O\big(H^{X,V}_{A,\tilde{O}} u\big)$.
\end{cor}
\begin{proof}
Since $H^{X,V}_{A,O}=H^X_{A,O}+V$ with $D(H^{X,V}_{A,O})=D(H^X_{A,O})$ (as $V$ is bounded), the claim follows directly from Lemmas~\ref{rest-prt-dir} and~\ref{rest-prt-dir-neu}.
\end{proof}

\noindent Compactly supported Dirichlet domain vectors extend canonically to larger Dirichlet domains.
\begin{lem}
\label{ext-dir}
Let $O\subset \tilde{O}\subseteq \mathbb{R}^d$ be open sets and $A\in\big( L_{\mathrm{loc}}^2(\tilde{O}) \big)^d$ real-valued. Consider $H^{D}_{A,O}$ and $H^{D}_{A,\tilde{O}}$ as in
Lemma \ref{dir-neu-lp}. If $u\in D\big(H^{D}_{A,O}\big)$ with $\operatorname{supp}u\subset\subset O$, then
$$\chi^*_O u\in D\big( H^{D}_{A,\tilde{O}} \big)~~\text{and}~~H^{D}_{A,\tilde{O}}\big(\chi^*_O u\big)=\chi^*_O \big(H^{D}_{A,O}u\big)\in L^2(\tilde{O}),$$
\end{lem}
\begin{proof}
Let $u\in D(H^{D}_{A,O})\subset W^{1,2}_{A,cl}(O)$ with 
$\operatorname{supp}u\subset\subset O$. 
By Proposition \ref{cn-embd}, we have $\chi_O^*u\in W^{1,2}_{A, cl} (\tilde{O})$ and
$(i\nabla+A)(\chi_O^* u)=\chi_O^*(i\nabla+A)u$. Since $\operatorname{supp}u\subset\subset O\subset\tilde O$, 
choose an open set $O_u$  with $\operatorname{supp}u\subset O_u\subset O,
$ and take $\psi\in C_c^\infty(\tilde O)$ with 
$\psi=1$ on $O_u$. For $v\in W^{1,2}_{A,cl}(\tilde O)$. 
arguing as in \eqref{cnt-lf}, we obtain
\begin{align*}
h^D_{A,\tilde O}(v,\chi_O^*u)
&=\langle (i\nabla+A)(\chi_O\psi v),
        (i\nabla+A)u\rangle_O.
\end{align*}
Since $\chi_O\psi v\in W^{1,2}_{A,cl}(O)$ and $u\in D(H^D_{A,O})$, this equals
\begin{align*}
&= h^D_{A,O}(\chi_O\psi v,u) = \langle \chi_O\psi v,
          H^D_{A,O}u\rangle_O
\end{align*}
As 
$\operatorname{supp}(H^D_{A,O}u)
\subset \operatorname{supp}u
\subset O_u \subset O \subset \tilde{O}$
and $\psi\equiv 1$ on $O_u$, we obtain
\begin{align*}
&= \langle \chi_O v,
          H^D_{A,O}u\rangle_O = \langle v,
          \chi_O^* H^D_{A,O}u\rangle_{\tilde O}.
\end{align*}
Thus $h^D_{A,\tilde O}(v,\chi_O^*u)
=
\langle v,
        \chi_O^* H^D_{A,O}u\rangle_{\tilde O}.$
Applying the same density argument as in 
Lemma \ref{rest-prt-dir}, we conclude that $\chi_O^*u\in D(H^D_{A,\tilde O})$ and $H^D_{A,\tilde O}(\chi_O^*u)=\chi_O^*(H^D_{A,O}u)\in L^2(\tilde O)$.
\end{proof}

\noindent The canonical extension of compactly supported domain vectors also holds for the magnetic Neumann Laplacian.
\begin{lem}
\label{ext-neu}
Let $O\subset \tilde{O}\subseteq \mathbb{R}^d$ be open sets and $A\in\big( L_{\operatorname{loc}}^2(\tilde{O}) \big)^d$ real-valued. Consider $H^{N}_{A,O}$ and $H^{N}_{A,\tilde{O}}$ as in Lemma \ref{dir-neu-lp}. If $u\in D\big(H^{N}_{A,O}\big)$ with $\operatorname{supp}u\subset\subset O$, then
$$\chi^*_O u\in D\big( H^{N}_{A,\tilde{O}} \big)~~\text{and}~~H^{N}_{A,\tilde{O}}\big(\chi^*_O u\big)=\chi^*_O \big(H^{N}_{A,O}u\big)\in L^2(\tilde{O}).$$
\end{lem}
\begin{proof}
Since $u\in  D\big(H^{N}_{A,O}\big)\subset W^{1,2}_{A}(O)$ with $\operatorname{supp}u\subset \subset O$, Lemma \ref{ext-0} implies that $\chi_O^*u\in W^{1,2}_{A} (\tilde{O})$ and
$(i\nabla+A)(\chi_O^* u)=\chi_O^*(i\nabla+A)u$. Now for any $v\in W^{1,2}_A(\tilde{O})$, we have
\begin{align}
\label{extsn-neu}
h^N_{A,\tilde{O}}\big(v, \chi^*_Ou \big)&=\big\langle (i\nabla+A)v,  (i\nabla+A)(\chi_O^*u)\big\rangle_{\tilde{O}}=\big\langle (i\nabla+A)v,  \chi_O^*(i\nabla+A)u\big\rangle_{\tilde{O}}\nonumber\\
&=\big\langle \chi_O (i\nabla+A)v,  (i\nabla+A)u\big\rangle_{O}=\big\langle  (i\nabla+A)(\chi_Ov),  (i\nabla+A)u\big\rangle_{O}\nonumber\\
&=h^N_{A,O}\big(\chi_Ov, u \big)=\big\langle \chi_Ov, H^N_{A,O}u \rangle_O=\big\langle v, \chi_O^*H^N_{A,O}u \rangle_{\tilde{O}}.
\end{align}
Now applying the same density argument used at the end of Lemma \ref{mult-dom-neu}, we conclude that $H^N_{A,\tilde{O}} (\chi_O^*u)=\chi_O^*\big(H^N_{A,O}u\big)$ for all $u\in D\big( H^N_{A,O} \big)$ with $\operatorname{supp}u\subset \subset O$.
\end{proof}

\begin{cor}
\label{ext-bth}
Let $\tilde O\subset\mathbb R^d$ be open and $V\in L^\infty(\tilde O)$. Consider $H^{X,V}_{A,\tilde O}$, $X=D,N$, as in Proposition~\ref{pertur}.  
If $u\in D(H^{X,V}_{A,O})$ with $\operatorname{supp}u\subset\subset O\subset\tilde O$. Then, under the assumptions of Lemmas \ref{ext-dir} and \ref{ext-neu}, we get $\chi_O^* u\in  D\big(H^{X,V}_{A,\tilde{O}} \big)$ and
$H^{X,V}_{A,\tilde{O}}(\chi_O^* u)=\chi_O^*\big(H^{X,V}_{A,O} u\big)$.
\end{cor}
\begin{proof}
Since $H^{X,V}_{A,O}=H^X_{A,O}+V$ with $D(H^{X,V}_{A,O})=D(H^X_{A,O})$ (as $V$ is bounded), the claim follows directly from Lemmas~\ref{ext-dir} and~\ref{ext-neu}.
\end{proof}

\noindent If a domain vector is compactly supported in the interior of $O$, the Dirichlet and Neumann magnetic Laplacian coincide on that vector.
\begin{prop}
\label{dr-nu-agr}
Let $O\subset\mathbb R^d$ be open and consider $H^{X,V}_{A,O}$, $X=D,N$, as in Proposition~\ref{pertur}. If $u\in D\big(H^{D,V}_{A,O}  \big)\cup D\big(H^{N,V}_{A,O}  \big)$ with $\operatorname{supp}u\subset \subset O$. Then $u\in D\big(H^{D,V}_{A,O}  \big)\cap D\big(H^{N,V}_{A,O}  \big)$ and $H^{D,V}_{A,O} u=H^{N,V}_{A,O} u$.
\end{prop}
\begin{proof}
If $O=\mathbb R^d$, then $H^{D,V}_{A,\mathbb R^d}=H^{N,V}_{A,\mathbb R^d}$ by Remark~\ref{dn-pertur}, so the claim is immediate.  
Assume $O\subset \mathbb{R}^d$ and let 
$u\in D\big(H^{D,V}_{A,O}\big)$ with $\operatorname{supp}u\subset\subset O$.
By Corollary \ref{ext-bth} with $\tilde O=\mathbb{R}^d$,
\[
\chi_O^*u\in D\big(H^{D,V}_{A,\mathbb{R}^d}\big), 
\qquad 
H^{D,V}_{A,\mathbb{R}^d}(\chi_O^*u)
=\chi_O^*\big(H^{D,V}_{A,O}u\big).
\]
Since $H^{D,V}_{A,\mathbb{R}^d}=H^{N,V}_{A,\mathbb{R}^d}$, it follows that
$\chi_O^*u\in D\big(H^{N,V}_{A,\mathbb{R}^d}\big)$, and because
$\operatorname{supp}\chi_O^*u\subset\subset O$, Corollary \ref{res-pertur-neu-dir}
yields $u=\chi_O(\chi_O^*u)\in D\big(H^{N,V}_{A,O}\big)$ and
\[
H^{D,V}_{A,O}u
=\chi_O\!\left(H^{D,V}_{A,\mathbb{R}^d}(\chi_O^*u)\right)
=\chi_O\!\left(H^{N,V}_{A,\mathbb{R}^d}(\chi_O^*u)\right)
=H^{N,V}_{A,\mathbb{R}^d}\!\left(\chi_O(\chi_O^*u)\right)
=H^{N,V}_{A,O}u.
\]
Thus $D(H^{D,V}_{A,O})\subset D(H^{N,V}_{A,O})$ on such vectors and the operators agree. The converse inclusion is obtained by the same argument
with the Dirichlet and Neumann operators interchanged.
\end{proof}

\begin{prop}
\label{wk-dr-nn-mg}
Let $A\in\big( L_{\mathrm{loc}}^2(O) \big)^d$ be real-valued and $O\subseteq \mathbb{R}^d$ open. If $f\in W^{1,2}_A(O)$ and $Af\in \big(L^2(O)\big)^d$, then the distributional gradient (without magnetic field) exists and $i\nabla f=(i\nabla+A)f-Af\in \big(L^2(O)\big)^d.$
Conversely, if $f\in L^2(O)$, $i\nabla f\in\big(L^2(O)\big)^d$, and $Af\in \big(L^2(O)  \big)^d$. 
Then $f\in W^{1,2}_A(O)$ and $(i\nabla+A)f= i\nabla f+Af\in \big(L^2(O)\big)^d.$
\end{prop}
\begin{proof}
Since $f\in W^{1,2}_A(O)$ and $Af\in (L^2(O))^d$, identity \eqref{sobo1} implies that for every $\varphi\in C_c^\infty(O)$ and each $k$, $\langle f,i\partial_k\varphi\rangle_O
=
\langle (i\nabla+A)_k f - A_k f,\varphi\rangle_O$.
Hence the distributional gradient (without magnetic field) $i\nabla f$ exists and $i\nabla f=(i\nabla+A)f-Af\in (L^2(O))^d$.
Conversely, suppose $f\in L^2(O)$, $i\nabla f\in (L^2(O))^d$, and $Af\in (L^2(O))^d$. For $\varphi\in C_c^\infty(O)$,
\begin{equation}
\label{con-cal}
\langle f,i\partial_k\varphi+A_k\varphi\rangle_O
=
\langle (i\nabla)_k f,\varphi\rangle_O
+
\langle A_k f,\varphi\rangle_O
=
\langle (i\nabla)_k f+A_k f,\varphi\rangle_O .
\end{equation}
Thus $(i\nabla+A)f=i\nabla f+Af\in (L^2(O))^d
$, and therefore $f\in W^{1,2}_A(O)$.
\end{proof}

\noindent It is well known that for bounded open sets $O\subset \tilde{O}\subseteq\mathbb{R}^d$, the operators $\chi_{O}(H^X_{A,\tilde{O}}+V-E)^{-m}\chi_{O}^*$ and $\chi_{O}(H_{A}+V-E)^{-m}\chi_{O}^*$ :$L^2(O)$ on $L^2(O)$, are trace class whenever $m>\frac{d}{2}$ and $E<-\|V\|_\infty$. Moreover, their trace norms admit explicit upper bounds.
\begin{prop}
\label{est-tr-chf}
let $\tilde{O}\subset \mathbb{R}^d$ be  bounded open set and consider $H^X_{A,\tilde{O}}~(X=D,N)$ and $H_A$ as in Lemma \ref{dir-neu-lp} and Remark \ref{fl-d-lp}. Let \( O \subset \tilde{O} \) be open and \( V \) be a bounded function.
Then, for $E<-\|V \|_\infty$ and $m>\frac{d}{2}$, the operators $\chi_{O}(H^X_{A,\tilde{O}}+V-E)^{-m}\chi_{O}^*$ and $\chi_{O}(H_A+V-E)^{-m}\chi_{O}^*$ are trace-class, and satisfy
\begin{equation}
\label{tr-est}
\big\| \chi_O (H^X_{A,\tilde O}+V-E)^{-m} \chi_O^* \big\|_1
\le C_1 |O|,
\qquad
\big\| \chi_O (H_A+V-E)^{-m} \chi_O^* \big\|_1
\le C_2 |O|,
\end{equation}
where
$C_1,C_2>0$ depend only on $E$ and $m$, and are independent of
$O$, $\tilde O$, $V$, and $A$.
\end{prop}
\begin{proof}
The second inequality in (\ref{tr-est}) is proved in \cite[(3.20)]{DIM}.  
For $V=0$, an estimate (see \cite[(4.39)]{TLMS1}) for $(H^X_{A,\tilde{O}}-E)^{-n}\chi_{O}^*:L^2(O)\to L^2(\tilde{O})$ is given as 
\begin{equation}
\label{hl-sm-p0}
\big\| (H^X_{A,\tilde{O}}-E)^{-n}\chi^*_O  \big\|^2_2 \leq C_1|O|,~~n>\frac{d}{4}.
\end{equation}
We may decompose $V-E=(V-E_1)-E_2$, where $E=E_1+E_2$ with fixed 
$E_1<-\|V \|_\infty$, and $E_2<0$ such that $V-E_1>0$.
Denote by $H^X_{0,\tilde{O}}$ the operator with $A\equiv 0$. By the diamagnetic inequality (see \cite[Theorem 3.3,  Remarks 3.4, (3.2) \& (3.3)]{HS}), 
\begin{equation*}
\big|(H^X_{A,\tilde{O}}+(V-E_1)-E_2)^{-n}\varphi\big|\leq (H^X_{0,\tilde{O}}-E_2)^{-n}|\varphi|~~\forall~\varphi\in L^2(\tilde{O}).
\end{equation*}
Hence for $\psi \in L^2(O)$,
\begin{equation}
\label{dm-iq}
\big| (H^X_{A,\tilde{O}}+V-E)^{-n}\chi_{O}^*\psi \big |\leq 
(H^X_{0,\tilde{O}}-E_2)^{-n}\chi_{O}^*|\psi|.
\end{equation}
The pointwise domination of $(H^X_{A,\tilde{O}}+V-E)^{-n}\chi_{O}^*$ by $(H^X_{0,\tilde{O}}-E_2)^{-n}\chi_{O}^*$ implies (see, \cite[Proposition 3.1]{DIM} \& \cite[Theorem 2.13]{bstr})
\begin{equation}
\label{in-hlsn}
\big\| (H^X_{A,\tilde{O}}+V-E)^{-n}\chi_{O}^* \big\|_2\leq 
\big\| (H^X_{0,\tilde{O}}-E_2)^{-n}\chi_{O}^* \big\|_2\leq \sqrt{C_1}~|O|^{\frac{1}{2}}.
\end{equation}
where the last step uses (\ref{hl-sm-p0}) with $A=0$. Writing $m=n_1+n_2$ with $n_1,n_2>\frac{d}{4}$,  we obtain
\begin{equation*}
\big\| \chi_{O}(H^X_{A,\tilde{O}}+V-E)^{-m}\chi_{O}^* \big\|_1\leq 
\big\| \chi_{O}(H^X_{A,\tilde{O}}+V-E)^{-n_1} \big\|_2\big\| (H^X_{A,\tilde{O}}+V-E)^{-n_2}\chi_{O}^* \big\|_2.
\end{equation*}
which together with \eqref{in-hlsn} yields the first estimate in \eqref{tr-est}. We use that an operator and its adjoint have the same Hilbert--Schmidt norm.
\end{proof}

\noindent Let \(f\) be a continuous function that decays
sufficiently fast at infinity. 
Then, using the above result, we can estimate the trace of the operator 
$\chi_O f\!\left(H_{A,  \tilde{O}}^X + V\right)$.
\begin{cor}
\label{tr-est-fn}
Let \(O\subset\tilde O\subset\mathbb{R}^d\) be bounded open sets, and let \(H_{A,\tilde{O}}^X\) and \(V\) be as in Proposition~\ref{est-tr-chf}.
Assume \(f \in C_0\big[-\|V\|_\infty, \infty\big)\) satisfies $|f(x)| = \mathcal{O}(|x|^{-m})$ as $ |x| \to \infty$,
for some \(m > 1 + \lfloor d/2 \rfloor\), and let \(E < -\|V\|_\infty\).  
Then 
\begin{equation}
\label{fn-tr-est}
\bigl|\operatorname{Tr}\!\bigl(\chi_O f(H_{A,  \tilde{O}}^X + V)\bigr)\bigr|
\;\leq\;  \|\breve{f}\|_\infty \, C \, |O|,
\quad 
\breve{f}(x) = (x - E)^{\,1 + \lfloor d/2 \rfloor} f(x).
\end{equation}
Here \(C > 0\) is independent of \(O\), \(\tilde{O}\), \(V\), and \(A\).
\end{cor}
\begin{proof}
Let $E < -\|V\|_\infty$ and choose integers $n_1, n_2 >\frac{d}{4}$ with 
$n_1 + n_2 = 1 + \lfloor d/2 \rfloor$.
By Corollary~\ref{cr-tr-rst},
\begin{align}
\operatorname{Tr}\!\bigl(\chi_O f(H_{A,  \tilde{O}}^X + V)\bigr)
&= \operatorname{Tr}\!\bigl(\chi_O f(H_{A,  \tilde{O}}^X + V)\chi_O^*\bigr)\nonumber\\
&= \operatorname{Tr}\!\biggl(\chi_O\,(H_{A,  \tilde{O}}^X + V - E)^{-n_1}
\,(H_{A,\tilde{O}}^X + V - E)^{1 + \lfloor d/2 \rfloor}\nonumber\\
& \qquad\qquad \times f(H_{A,  \tilde{O}}^X + V)\,(H_{A,  \tilde{O}}^X + V - E)^{-n_2}\chi_O^*\biggr).\nonumber
\end{align}
Using $|\operatorname{Tr}(T)| \leq \|T\|_1$ and $\|T_1T_2\|_1 \leq \|T_1\|_2 \|T_2\|_2$, we obtain
\begin{align*}
\bigl|\operatorname{Tr}\!\bigl(\chi_O f(H_{A,  \tilde{O}}^X + V)\bigr)\bigr|
&\leq C \,\bigl\| \chi_O (H_{A,  \tilde{O}}^X + V - E)^{-n_1} \bigr\|_2\,
\bigl\| (H_{A,  \tilde{O}}^X + V - E)^{-n_2}\chi_O^* \bigr\|_2 \\
&\leq  \|\breve{f}\|_\infty C |O|.
\end{align*}
In the last step, we used the Hilbert--Schmidt estimate~\eqref{in-hlsn} and the boundedness of \(\breve f\).
\end{proof}

\noindent A non-negative multiplication operator on $L^2(\mathbb{R}^d)$ cannot be trace class.
\begin{prop}
\label{trnf}
Let $u \in L^2(\mathbb{R}^d)$ with $\|u\|_2 \neq 0$ and $u \ge 0$. 
Define the multiplication operator $u$ on $L^2(\mathbb{R}^d)$ by $(u f)(x) = u(x) f(x), f \in L^2(\mathbb{R}^d)$. Then $\operatorname{Tr}(u) = \infty$.
\end{prop}
\begin{proof}
Since $\|u\|_2\neq 0$ and $u\ge0$, the set $\{x:u(x)>0\}$ has positive Lebesgue measure on $\mathbb{R}^d$. Hence there exists $n_0\in\mathbb{N}$ such that
$E:=\{x:u(x)>n_0^{-1}\}$ has positive measure. Let $\{f_n\}_{n=1}^\infty$
be an orthonormal basis of $L^2(E)$. Then
\[
\operatorname{Tr}(u)\ge \operatorname{Tr}(\chi_E u)
= \sum_{n=1}^\infty \int_E u(x)\,|f_n(x)|^2\,dx
> n_0^{-1}\sum_{n=1}^\infty \int_E |f_n(x)|^2\,dx
= n_0^{-1}\sum_{n=1}^\infty 1=\infty.
\]
Thus $\operatorname{Tr}(u)$ is not finite.
\end{proof}

\noindent The following result gives a differentiation formula for the trace of a function of a parameter-dependent self-adjoint operator, based on the Hellmann--Feynman theorem.
\begin{prop}
\label{tr-drv-frm}
Let \( T_\lambda = T + \lambda K \) be a family of self-adjoint
operators operator on a Hilbert space \( \mathcal{H} \), where \( K \) is bounded and self-adjoint. Assume that \( \sigma(T_\lambda) \) is discrete for each \( \lambda \in (a, b) \). Let \( f \) be differentiable real-valued function on \( \sigma(T_\lambda) \) such that \( \operatorname{Tr}\big(f(T_\lambda)\big) \) and \( \operatorname{Tr}\big( K f'(T_\lambda)\big) \) are finite, and $\operatorname{Tr}\big(|f'(T_\lambda)|\big) \leq \sum_n M_n < \infty ~ \text{for all } \lambda \in (a,b)$,
where \( M_n \) is independent of \( \lambda \). Then the derivative of the trace of \( f(T_\lambda) \) is given by
\begin{equation}
\label{d-tr-frm}
\frac{d}{d\lambda} \operatorname{Tr}\big(f(T_\lambda)\big) = \operatorname{Tr}\big( K f'(T_\lambda) \big).
\end{equation}
\end{prop}
\begin{proof}
Let $\{E_{n,\lambda}\}_n$ be the eigenvalues of $T_\lambda$ with normalized
eigenvectors $\{\psi_n\}_n$. Then $\operatorname{Tr}\big(f(T_\lambda)\big)=\sum_n f(E_{n,\lambda})$. Using the Hellmann-Feynman theorem (see \cite{IM}), we obtain $\frac{d}{d\lambda}E_{n,\lambda}=\langle \psi_n,K\psi_n\rangle=\langle K\psi_n,\psi_n\rangle,
$ and hence
\begin{equation}
\label{drv-tr}
\frac{d}{d\lambda}\operatorname{Tr}\big(f(T_\lambda)\big)
=\sum_n f'(E_{n,\lambda})\langle K\psi_n,\psi_n\rangle
=\sum_n \langle K\psi_n,f'(T_\lambda)\psi_n\rangle
=\operatorname{Tr}\big(Kf'(T_\lambda)\big).
\end{equation}
The differentiation term by term is justified since $\sum_n\big|f'(E_{n,\lambda})\langle \psi_n,K\psi_n\rangle\big|
\le \|K\|\sum_n M_n<\infty,$ uniformly for $\lambda\in(a,b)$.
\end{proof}

\noindent The following corollary applies the trace differentiation formula to resolvent powers of magnetic Schr\"{o}dinger operators with bounded potentials. It provides an explicit expression for the derivative of the trace with respect to a coupling parameter.
\begin{cor}
\label{drv-tr-rsl}
Let \( O \subset \mathbb{R}^d \) be bounded open and consider \( H^X_{A,O} \) for \( X = D, N \) as in Lemma~\ref{dir-neu-lp}. Assume \( U,V \) are bounded functions on \( O \) and \( \lambda \in (-b, b) \). Then, for each \( E < -\|U\|_\infty - b\|V\|_\infty \), and $m > \frac{d}{2}$, we have the following trace derivative formula:
\begin{align}
\label{drv-tr-rs} 
\frac{d}{d\lambda}
\operatorname{Tr}\big( H^X_{A,O}+U+\lambda V-E \big)^{-m}
=
-m\,\operatorname{Tr}
\!\left(
V\big(H^X_{A,O}+U+\lambda V-E\big)^{-m-1}
\right).
\end{align}
\end{cor}
\begin{proof}
Set \( T_\lambda = H^X_{A,O} + U + \lambda V \) and \( f(x) = (x - E)^{-m} \). Let \( \{M_n\} \) be the eigenvalues of the positive operator \( (H^X_{A,O} - \|U\|_\infty - b \|V\|_\infty + E)^{-m-1} \). It follows that 
$\operatorname{Tr}\big( \big|f'(T_\lambda)\big| \big) \leq m \sum_n M_n < \infty.
$
Thus, \eqref{drv-tr-rs} follows from Proposition~\ref{tr-drv-frm}.
\end{proof}
\noindent The following Combes--Thomas estimate provides exponential off-diagonal decay of the resolvent of magnetic Schr\"odinger operators.
\begin{thm}
Let \( O \subseteq \mathbb{R}^d \) be an open and consider \( H^X_{A,O} \) for \( X = D, N \) and \( H_A \) as in Lemma~\ref{dir-neu-lp} and Remark~\ref{fl-d-lp}. Assume \( V \) is bounded function on \( O \), let \( m > \frac{d}{2p} \), and let \( F,G \subset O \) be bounded. Then, for every \( E < -\|V\|_\infty \), 
\begin{equation}
\label{exp-decy}
\big\| \chi_{F} \big(H^X_{A,O} + V - E\big)^{-m} \chi_{G} \big\|_p \leq C e^{-\beta \operatorname{dist}(F,G)} \quad \text{for all}\quad p \in \mathbb{N} \cup \{\infty\},
\end{equation}
where \( C,\beta > 0 \) depend only on \( E \), \( m \), and \( p \), but are independent of \( V \), \( A \), \( O \), and \( X \). 
Here, \( H^X_{A,\mathbb{R}^d} := H_A \), and for \( p = \infty \), the norm refers to the operator norm.
\end{thm}
\begin{proof}
For a proof we refer to \cite{gerkl1,gerkl, CHF, BCH, She}.
\end{proof}

\noindent We decompose $O$ into finitely many disjoint open subsets covering it up to a null set. For trace-class operators, the trace can then be computed locally on each piece and then summed.
\begin{prop}
Let $O \subseteq \mathbb{R}^d$ be open and $\{O_k\}_{k=1}^M$ pairwise disjoint with $O_k \subset O$,
$O = \left( \overline{\bigcup_{k=1}^M O_k} \right)^\circ$, and
$\left| O \setminus \bigcup_{k=1}^M O_k \right| = 0.$
If $T$ is trace-class on $L^2(O)$. Then
\begin{equation}
\label{tr-idt}
\operatorname{Tr}(T) 
= \sum_{k=1}^M \operatorname{Tr} \left( \chi_{O_k} T \chi_{O_k}^* \right).
\end{equation}
\end{prop}
\begin{proof}
The Hilbert spaces $L^2(O)$ and $\bigoplus_{k=1}^M L^2(O_k)$ are unitarily equivalent via unitary operator,   
$U : L^2(O) \to \bigoplus_{k=1}^M L^2(O_k)$, together with
its adjoint, is given by
\[
Uf = \bigoplus_{k=1}^M \big( \chi_{O_k} f \big)
~~\&~~
U^*\!\left( \bigoplus_{k=1}^M f_k \right) 
= \sum_{k=1}^M \chi_{O_k}^* f_k,~~f\in L^2(O)~~f_k\in L^2(O_k).
\]
Let $\{\psi_{k,n}\}_{n=1}^\infty$ be an orthonormal basis for $L^2(O_k)$.  
Then $$\Bigg\{
  0 \oplus \cdots \oplus 0 
  \oplus \underbrace{\psi_{k,n}}_{\text{$k$-th position}} 
  \oplus ~0 \oplus \cdots \oplus 0 
  : n = 1, 2, \ldots
\Bigg\}_{k=1}^M$$
forms an orthonormal basis for $\bigoplus_{k=1}^M L^2(O_k)$.  
Since $\operatorname{Tr}(T) = \operatorname{Tr}(UTU^*)$, writing $\operatorname{Tr}(UTU^*)$ with respect to the above basis yields (\ref{tr-idt}).
\end{proof}
\noindent The following corollary is immediate.
\begin{cor}
\label{cr-tr-rst}
Let $T$ be a trace-class operator on $L^2(\tilde{O})$, where $\tilde{O}\subseteq \mathbb{R}$ is open.  If $O \subset \tilde{O}$, then 
$\operatorname{Tr}(\chi_O T) \;=\; \operatorname{Tr}(\chi_O T \chi_O^*)$.
\end{cor}

\noindent Now we state a criterion for equivalence of convergence in distribution.
\begin{prop}
\label{cn-vr-rv}
Let $\{X_n\}$ and $\{Y_n\}$ be sequences of random variables on $(\Omega, \mathcal{B}_\Omega, \mathbb{P})$ such that $\mathbb{E}\big[X_n - Y_n\big]^2 \to 0$ as $n \to \infty$.
Suppose $Y_n \to Y$ in distribution. Then:\\
(i) $X_n \to Y$ in distribution as $n\to\infty$. \\
(ii) If $\displaystyle\lim_{n\to\infty}\mathbb{E}[Y_n^2]$ exists, then 
$\displaystyle\lim_{n\to\infty}\mathbb{E}[X_n^2] \;=\; \lim_{n\to\infty}\mathbb{E}[Y_n^2]$.
\end{prop}
\begin{proof}
(i) Convergence in quadratic mean, i.e.,  $\mathbb{E}\!\left[(X_n - Y_n)^2\right] \;\to\; 0$ as $n \to \infty$,  implies that $X_n - Y_n \xrightarrow{d} 0.$
Writing $X_n = (X_n - Y_n) + Y_n$,  the claim follows from Slutsky’s theorem.\\~\\
(ii) Using Minkowski's inequality, we can write
\begin{align*}
\big(\mathbb{E}[X_n^2]\big)^{\tfrac{1}{2}}
&\leq \big(\mathbb{E}\big[(X_n - Y_n)^2\big]\big)^{\tfrac{1}{2}}
   + \big(\mathbb{E}[Y_n^2]\big)^{\tfrac{1}{2}}, \\
\big(\mathbb{E}[Y_n^2]\big)^{\tfrac{1}{2}}
&\leq \big(\mathbb{E}\big[(X_n - Y_n)^2\big]\big)^{\tfrac{1}{2}}
   + \big(\mathbb{E}[X_n^2]\big)^{\tfrac{1}{2}}.
\end{align*}
Hence,
$\bigg| \big(\mathbb{E}[X_n^2]\big)^{\tfrac{1}{2}}
      - \big(\mathbb{E}[Y_n^2]\big)^{\tfrac{1}{2}} \bigg|
   \leq \big(\mathbb{E}\big[(X_n - Y_n)^2\big]\big)^{\tfrac{1}{2}}\to0$,
because $X_n - Y_n \to 0$ in quadratic mean as $n \to \infty$, which proves the assertion.
\end{proof}

\begin{prop}
\label{vr-sm-lm}
Let $\{X_n\}$ and $\{Y_n\}$ be sequences of random variables on $(\Omega, \mathcal{B}_\Omega, \mathbb{P})$ such that $\displaystyle\sup_n \mathbb{E}[X_n^2] < \infty$ and $\displaystyle\lim_{n\to\infty} \mathbb{E}[Y_n^2] = 0$. Then
$\displaystyle\lim_{n\to\infty} \left[\mathbb{E}\big[(X_n+Y_n)^2\big] - \mathbb{E}[X_n^2]\right]= 0$. 
\end{prop}
\begin{proof}
Write $\mathbb{E}\big[(X_n+Y_n)^2\big]=\mathbb{E}[X_n^2]+\mathbb{E}[Y_n^2]
+2\mathbb{E}[X_nY_n]$.   An application of the Cauchy-Schwarz inequality then yields the desired result.
\end{proof}

\noindent Next, we state the Burkholder inequality, which is useful for estimating the fourth moments of random variables that will appear in our work.
\begin{thm}
\label{Burkholder}
Let $(M_k)_{k=0}^n$ be a real-valued martingale with $M_0=0$ and quadratic variation
$[M]_n \;=\; \sum_{k=1}^n (M_k - M_{k-1})^2$.
Then there exist constants $0<c_p\le C_p<\infty$ (depending only on $p\ge 1$) such that
\begin{equation}
\label{binq}
c_p\,\mathbb{E}\!\left([M]_n^{\frac{p}{2}}\right)
\;\le\;
\mathbb{E}\!\left( |M_n|^{p}\right)
\;\le\;
C_p\,\mathbb{E}\!\left([M]_n^{\frac{p}{2}}\right).
\end{equation}
\end{thm}
\begin{proof}
A proof can be found in \cite{Bur}.
\end{proof}

\noindent Using the ergodicity of the random operator \( H^\omega \), we show that the distributions of the traces \(\operatorname{Tr}(\omega_n u_n f(H^\omega))\) are identical.
\begin{prop}\label{trids}
Consider \( H^\omega \) as in \eqref{model} under Hypothesis~\ref{hypo}.  
Let \( f \) be a continuous function such that \((x - E)^{1 + \lfloor d/2 \rfloor} f(x)\) is bounded on \( \big[-\|V\|_\infty, \infty\big) \),  
where \( E < -\|V\|_\infty \).  
Then the random variable \( \omega_n \operatorname{Tr}\!\big(u_n f(H^\omega)\big) \) has the same distribution for all \( n \in \mathbb{Z}^d \),  
where \( u_n(x) = u(x - n) \).
\end{prop}
\begin{proof}
Let \( T_m \) and \( U_m \)  be as in \eqref{mpt} and \eqref{unop}.  
By \cite[Lemma~4.5]{Kir}, we have
\[
U_m H^\omega U_m^* = H^{T_m \omega}, 
\quad \text{equivalently,} \quad
U_m f(H^\omega) U_m^* = f(H^{T_m \omega}).
\]
Define $Y_n(\omega) = \omega_n \operatorname{Tr}\big(u_n f(H^\omega)\big)$. Let \( S = \operatorname{supp}(u),~S_n = \operatorname{supp}(u_n) = S + n \) and \( \chi_n = \chi_{S_n}\).  
For each \( m \in \mathbb{Z}^d \), define the unitary operator  
\(\tilde{U}_m : L^2(S_n) \to L^2(S_{n-m})\) by
\[
(\tilde{U}_m \varphi)(x) = e^{-i\Psi_m(x+m)} \varphi(x + m),
\quad \varphi \in L^2(S_n), \quad x \in S_{n-m},
\]
where \(\Psi_m\) is is given in \eqref{unop}.  
By Corollary \ref{cr-tr-rst}, 
\begin{align*}
Y_n(\omega)
= \omega_n \operatorname{Tr}\big(u_n f(H^\omega)\big) 
= \omega_n \operatorname{Tr}\big(u_n \chi_n f(H^\omega) \chi_n^*\big).
\end{align*}
Then
\begin{align}
\label{trmpt}
Y_n(T_m \omega)
&= (T_m \omega)_n \operatorname{Tr}\big(u_n \chi_n f(H^{T_m \omega}) \chi_n^*\big) \nonumber\\
&= (T_m \omega)_n \operatorname{Tr}\big(\tilde{U}_m u_n \chi_n U_m f(H^\omega) U_m^* \chi_n^* \tilde{U}_m^*\big).
\end{align}
One verifies that
$\tilde{U}_m u_n \chi_n U_m = \chi_{n-m} u_{n-m}$ and $U_m^* \chi_n^* \tilde{U}_m^* = \chi_{n-m}^*$.
By \eqref{trmpt},
\begin{align}
\label{smdis}
Y_n(T_m \omega)
&= \omega_{n-m} \operatorname{Tr}\big(u_{n-m} \chi_{n-m} f(H^{ \omega}) \chi_{n-m}^*\big) \nonumber\\
&= \omega_{n-m} \operatorname{Tr}\big(u_{n-m} f(H^{ \omega})\big) = Y_{n-m}(\omega).
\end{align}
Since \(T_m : \Omega \to \Omega\) is a measure-preserving transformation,  
it follows that the family \(\{Y_n\}_{n \in \mathbb{Z}^d}\) consists of identically distributed  
(but not necessarily independent) random variables.
\end{proof}

\begin{cor}
\label{smdscn}
Let \( g(H^\omega)_{(\omega_k \to t\omega_k)} \) be defined as in~\eqref{mdfmdl}, and let 
\( \mathcal{F}^d_{k_1, k_2, \ldots, k_{d-1}, k_d} \) be the $\sigma$-algebra introduced in~\eqref{dsgal}. 
Then, under Hypothesis~\ref{hypo}, the conditional expectations.
\begin{equation*}
\mathbb{E}\!\left(
  g(H^\omega)_{\big(\omega_{(k_1,k_2,\ldots,k_{d-1},k_d)} \to t\omega_{(k_1,k_2,\ldots,k_{d-1},k_d)}\big)}
  \,\middle|\, \mathcal{F}^d_{(k_1,k_2,\ldots,k_{d-1},k_d)}
\right)
\end{equation*}
have identical distributions for all \( (k_1,k_2,\ldots,k_{d-1},k_d) \in \mathbb{Z}^d \).
The same conclusion holds when conditioning on
\( \mathcal{F}^d_{(k_1,k_2,\ldots,k_{d-1},k_d-1)}\).
\end{cor}
\begin{proof}
The result follows from the following two observations. 
First, the collection $\{\omega_n\}_{n \in \mathbb{Z}^d}$ consists of independent and identically distributed random variables. 
Second, 
$\mathcal{F}^d_{(k_1, k_2, \ldots, k_{d-1}, k_d)} 
= \sigma\{\omega_n : n \in A^d_{(k_1, k_2, \ldots, k_{d-1}, k_d)}\}$.
By (\ref{dsgal}), the sets 
$A^d_{(k_1, k_2, \ldots, k_{d-1}, k_d)}$ and 
$A^d_{(m_1, m_2, \ldots, m_{d-1}, m_d)}$ 
are isomorphic. 
The same property also holds for the $\sigma$-algebra 
$$\mathcal{F}^d_{(k_1, k_2, \ldots, k_{d-1}, k_d - 1)} 
= \sigma\{\omega_n : n \in A^d_{(k_1, k_2, \ldots, k_{d-1}, k_d - 1)}\}.$$
\end{proof}


\begin{thebibliography}{99} 

\bibitem{AMSW} Aizenman, M., Warzel, S.: \textsl{Random Operators: Disorder Effects on Quantum Spectra and Dynamics}, Graduate Studies in Mathematics, vol. 168. Amer. Math. Soc., Providence, Rhode Island, 2015.

\bibitem{AW06}
Aizenman, M., Warzel, S.:
\emph{The canopy graph and level statistics for random operators on trees},
Math.\ Phys.\ Anal.\ Geom.\ \textbf{9(4)}, 291--333, 2006.

\bibitem{apw}
Anderson, P. W.: \textsl{Absence of diffusion in certain random lattices}, Phys. Rev. {\bf 109}, 1492--1505, 1958.

\bibitem{BCH} Barbaroux, J.M., Combes, J.M., Hislop, P.D.: \textsl{Localization near band edges for random Schr\"{o}dinger operators}, Helv. Phys. Acta. {\bf 70}, 16--43, 1997.

\bibitem{PB} Billingsley, P.: \textsl{Probability and Measure}, Wiley Series in Probability and Mathematical Statistics, 3rd edn, John Wiley \& Sons, Inc, New York, 1995. 

\bibitem{BGW} Breuer, J., Grinshpon, Y., White, M.J.: \textsl{Spectral Fluctuations for Schr\"{o}dinger Operators with a Random Decaying Potential}, Ann. Henri Poincar\`e, {\bf 22}, 3763--3794, 2021.

\bibitem{BHL} Broderix, K., Hundertmark, D., Leschke, H.: \textsl{Self-averaging, decomposition and asymptotic properties of the density of states for random Schr\"{o}dinger operators with constant magnetic field}, World Scientific Publishing Co., Inc., River Edge, NJ, 98--107, 1993.

\bibitem{Bur} Burkholder, D.L.: \textsl{Martingale transforms}. Ann. Math. Stat. {\bf 37}, 1494--1504, 1966.

\bibitem{CL} Carmona, R., Lacroix, J.: \textsl{Spectral Theory of Random Schr\"{o}dinger Operators}, Boston, Birkhauser, 1990.

\bibitem{CLZ}
Cannizzaro, G., Labb\'e, C., Zuijlen, W.:
\emph{Top of the spectrum of discrete Anderson Hamiltonians with correlated Gaussian potentials},
arXiv:2505.06051.

\bibitem{CHF} Combes, J.M., Hislop, P.D., Klopp, F.: \textsl{An optimal Wegner estimate and its application to the global continuity of the integrated density of states for random Schr\"{o}dinger operators}, Duke Math. J. {\bf 140(3)}, 469--498, 2007.

\bibitem{conwy} Conway, J. B.: \textsl{A Course in Functional Analysis}, (Graduate Texts in Mathematics, 96), New York, Springer, 1990.

\bibitem{CFKS} Cycon, H.L., Froese, R.G., Kirsch, W., Simon, B.: \textsl{Schr\"{o}dinger operators: with application to quantum mechanics and global geometry}, Texts and Monographs in Physics, Springer-Verlag, Berlin, 1987.

\bibitem{DE21}
Dietlein, A., Elgart, A,:
\emph{Level spacing and Poisson statistics for continuum random Schr\"odinger operators},
J.\ Eur.\ Math.\ Soc. \textbf{23(4)}, 1257--1293, 2021.

\bibitem {DCLT} Dolai, D.: \textsl{Central limit theorem for the random variables associated with the IDS of the Anderson model on lattice}; arXiv:2309.07529.

\bibitem{DIM} Doi, S., Iwatsuka, A., Mine, T.: \textsl{The uniqueness of the integrated density of states for the Schr\"{o}dinger operators with magnetic fields}, Math. Z. {\bf 237(2)}, 335--371, 2001.

\bibitem{DL20}
Dumaz, L., Labb\'e, C.:
\emph{Localization of the continuous Anderson Hamiltonian in 1-D},
Probab.\ Theory Relat.\ Fields.\ \textbf{176}, 353--419, 2020.

\bibitem{DL24}
Dumaz, L., Labb\'e, C.
\emph{Localization crossover for the continuous Anderson Hamiltonian in 1-d},
Invent.\ Math.\ \textbf{235(2)}, 345--440, 2024.

\bibitem{fos} Forrester, P.J.: \textsl{A review of exact results for fluctuation formulas in random matrix theory}, Probability Surveys. {\bf 20}, 170--225, 2023.

\bibitem{gal} Galdi, G.: \textsl{An introduction to the mathematical theory of the Navier-Stokes equations: Steady-state problems}, New York, Dordrecht, Heidelberg, London: Springer, 2011.

\bibitem{gerkl} Germinet, F., Klein, A.: \textsl{Explicit finite-volume criteria for localization in continuous random media and applications}, Geom. Funct. Anal. {\bf 13(6)}, 1201--1238, 2003.

\bibitem{gerkl1} Germinet, F., Klein, A.: \textsl{Operator kernel estimates for functions of generalized Schr\"{o}dinger operators}. Proc. Amer. Math. Soc. {\bf 131(3)}, 911--920, 2003.

\bibitem{GK14}
Germinet, F., Klopp, F.:
\emph{Spectral statistics for random Schr\"odinger operators in the localized regime},
J.\ Eur.\ Math.\ Soc.\ \textbf{16(9)}, 1967--2031, 2014.

\bibitem{GYWJ} Grinshpon, Y., White, M.J.: \textsl{Spectral fluctuations for the multi-dimensional Anderson model}, J. Spectr. Theory {\bf 12(2)}, 591--615, 2022.

\bibitem{His} Hislop, P.D.: \textsl{Lectures on random Schr\"{o}dinger operators}, Contemp. Math., Amer. Math. Soc., Providence, RI, {\bf 476}, 41--131, 2008.

\bibitem{HK15}
Hislop, P. D., Krishna, M.:
\emph{Eigenvalue statistics for random Schr\"odinger operators with non rank one perturbations},
Comm.\ Math.\ Phys.\ \textbf{340(1)}, 125--143, 2015.

\bibitem{HKK20}
Hislop, P. D., Kirsch, W., Krishna, M.:
\emph{Eigenvalue statistics for Schr\"odinger operators with random point interactions on $\mathbb{R}^d$, $d=1,2,3$},
J.\ Math.\ Phys.\ \textbf{61(9)}, 092103, 24~pp, 2020.

\bibitem{HS} Hundertmark, D., Simon, B.: \textsl{A diamagnetic inequality for semigroup differences}, J. Reine Angew. Math. {\bf 571}, 107--130, 2004.

\bibitem{TLMS} Hupfer, T., Leschke, H., M\"{u}ller, P., Warzel, S.: \textsl{The absolute continuity of the integrated density of states for magnetic Schr\"{o}dinger operators with certain unbounded random potentials}, Comm. Math. Phys. {\bf 221(2)}, 229--254, 2001.

\bibitem{TLMS1} Hupfer, T., Leschke, H., M\"{u}ller, P., Warzel, S.: \textsl{Existence and uniqueness of the integrated density of states for Schr\"{o}dinger operators with magnetic fields and unbounded random potentials}, Rev. Math. Phys. {\bf 13(12)}, 1547--1581, 2001.

\bibitem{IM} Ismail, M.E.H., Muldoon, M.E.: \textsl{On the variation with respect to a parameter of zeros of Bessel and $q$-Bessel functions}, J. Math. Anal. Appl. {\bf 135(1)}, 187--207, 1988.

\bibitem{kato} Kato, T.: \textsl{Perturbation theory for linear operators}, Berlin, Heidelberg, New York: Springer, 1980.

\bibitem{KT} Kato, T.: \textsl{Remarks on Schr\"{o}dinger operators with vector potentials}, Integral Equations Operator Theory. {\bf 1}, 103--113, 1978.

\bibitem{KB} Kirsch, W., Metzger, B.: \textsl{The integrated density of states for random Schr\"{o}dinger operator}, Proc. Sympos. Pure Math., 76, Part 2, Amer. Math. Soc., Providence, RI, 649--696, 2007.

\bibitem{Kir} Kirsch, W.: \textsl{An invitation to random Schr\"{o}dinger operators}, (With an appendix by Fr{\'e}d{\'e}ric Klopp) Panor. Synth{\'e}ses, {\bf 25}, Random Schr\"{o}dinger operators, {\bf 1}, Soc. Math. France, Paris, 1--119, 2008.

\bibitem{Wernr} Kirsch, W.: \textsl{Random Schr\"{o}dinger operators a course}, Lecture Notes in Phys. {\bf 345}, Springer-Verlag, Berlin, 1989.

\bibitem{KM} Kirsch, W., Martinelli, F.: \textsl{On the density of states of Schr\"{o}dinger operators with a random potential}, J. Phys. A. {\bf 15}, 2139--2156, 1982.

\bibitem{KMcre} Kirsch, W., Martinelli, F.: \textsl{On the ergodic properties of the spectrum of general random operators}, J. Reine Angew. Math. {\bf 334}, 141--156, 1983.

\bibitem{WF} Kirsch, W., Martinelli, F.: \textsl{On the essential selfadjointness of stochastic Schr\"{o}dinger operators}, Duke Math. J. {\bf 50(4)}, 1255--1260, 1983.

\bibitem{KP} Kirsch, W., Pastur, L.A.: \textsl{Analogues of Szeg{\H{o}}'s  theorem for ergodic operators}, Mat. Sb. {\bf 206(1)}, 103--130, 2015.

\bibitem{Klo13}
Klopp, F.:
\emph{Asymptotic ergodicity of the eigenvalues of random operators in the localized phase},
Probab.\ Theory Relat.\ Fields \textbf{155(3--4)}, 867--909, 2013.

\bibitem{KN}
Kotani, S., and Nakano, F.:
\emph{Poisson statistics for 1d Schr\"odinger operators with random decaying potentials},
Electron.\ J.\ Probab.\ \textbf{22}, 1--31, 2017.

\bibitem{KVV}
Kritchevski, E., Valk\'o, B., and Vir\'ag, B.:
\emph{The scaling limit of the critical one-dimensional random Schr\"odinger operator},
Commun.\ Math.\ Phys.\ \textbf{314}, 775--806, 2012.

\bibitem{LL} Lieb, E.H., Loss, M.: \textsl{Analysis}, Providence, Rhode Island: Amer. Math. Soc., 1997.

\bibitem{MD2019}
Mallick, A., Dolai, D.: \textsl {Spectral statistics for one-dimensional Anderson model with unbounded but decaying potential}, \emph{Infinite Dimensional Analysis, Quantum Probability and Related Topics}
\textbf{22(2)}, 1950012, 14pp, 2019.

\bibitem{MH} Matsumoto, H.: \textsl{On the integrated density of states for the Schr\"{o}dinger operators with certain random electromagnetic potentials}, J. Math. Soc. Japan. {\bf 45}, 197--214, 1993.

\bibitem{MMMN} Mashiko, T., Marui, Y., Maruyama, N., Nakano, F.: \textsl{Eigenvalue fluctuations of 1--dimensional random Schr\"{o}dinger operators}, J. Math. Phys. {\bf 65(8)}, Paper No. 083301, 2024.

\bibitem{Min96}
Minami, N.:
\emph{Local fluctuation of the spectrum of a multidimensional Anderson tight binding model},
Comm.\ Math.\ Phys.\ \textbf{177(3)}, 709--725, 1996.

\bibitem{Mol81}
Mol\v{c}anov, S.~A.:
\emph{The local structure of the spectrum of the one-dimensional Schr\"odinger operator},
Comm.\ Math.\ Phys.\ \textbf{78}, 429--446, 1981.

\bibitem{Na} Nakamura, S.: \textsl{A remark on the Dirichlet-Neumann decoupling and the integrated density of states}, J. Funct. Anal. {\bf 179(1)}, 136--152, 2001.

\bibitem{Nas} Nakao, S.: \textsl{On the spectral distribution of the Schr\"{o}dinger operator with random potential}, Japan. J. Math. {\bf 3}, 111--139, 1977.

\bibitem{Nak14}
Nakano, F.:
\emph{Level statistics for one-dimensional Schr\"odinger operators and Gaussian beta ensemble},
J.\ Stat.\ Phys.\ \textbf{156(1)}, 66--93, 2014.

\bibitem{nfjsp} Nakano, F.: \textsl{Fluctuation of density of states for 1d Schr\"{o}dinger operators}, J. Stat. Phys. {\bf 166(6)}, 1393--1404, 2017.

\bibitem{pasl} Pastur, L.: \textsl{On the Schr\"{o}dinger equation with a random potential}, Theor. Math. Phys. {\bf 6}, 299--306, 1971.

\bibitem{PF} Pastur, L., Figotin, A.: \textsl{Spectra of random and almost-periodic operators}, Berlin: Springer, 1992.

\bibitem{pas-sh} Pastur, L., Shcherbina, M.: \textsl{ Szeg{\H{o}}'s-type theorems for one-dimensional Schr\"{o}dinger operator with random potential (smooth case)}, J. Math. Phys. Anal. Geom. {\bf 14(3)}, 362--388, 2018.

\bibitem{PLSM} Pastur, L., Shcherbina, M.: \textsl{Eigenvalue distribution of large random matrices.}, Mathematical Surveys and Monographs, 171. Amer. Math. Soc., Providence, RI, 2011.

\bibitem{Rez} Reznikova, A.Y.: \textsl{A central limit theorem for the spectrum of random Jacobi matrices}, Teor. Veroyatnost. i Primenen. {\bf 25(3)}, 513--522, 1980.

\bibitem{R} Reznikova, A.Y.: \textsl{Limit theorems for the spectra of ID random Schr\"{o}dinger operator}, Random Oper. and Stoch. {\bf 12(3)}, 235--254, 2004.
 
\bibitem {R2} Re\u{z}nikova, A.J.: \textsl{The central limit theorem for the spectrum of the random one-dimensional Schr\"{o}dinger operator}, J. Stat. Phys. {\bf25 (2)}, 291--308, 1981.

\bibitem{RS} Reed, M., Simon, B.: \textsl{Methods of modern mathematical physics I: Functional analysis}, Revised and enlarged edition, San Diego: Academic, 1980.

\bibitem{RS4} Reed, M., Simon, B.: \textsl{Methods of Modern Mathematical Physics IV: Analysis of Operators}, New York, Academic, 1978.

\bibitem{She} Shen, Z.: \textsl{An improved Combes--Thomas estimate of magnetic Schr\"{o}dinger operators}, Arkiv f{\"o}r Matematik, {\bf 52(2)}, 383--414, 2014.

\bibitem{bary} Simon, B.: \textsl{Maximal and minimal Schr\"{o}dinger forms}, J. Operator Theory. {\bf 1}, 37--47, 1979.

\bibitem{bstr} Simon, B.: \textsl{Trace Ideals and Their Applications}, Mathematical Surveys and Monographs, vol. 20. Amer. Math. Soc., Providence, 2005.

\bibitem{U} Ueki, N.: \textsl{On spectra of random Schr\"{o}dinger operators with magnetic fields}, Osaka J. Math. {\bf 31}, 177--187, 1994.

\bibitem{ves} Veseli{\'c}, I.: \textsl{Existence and regularity properties of the integrated density of states of random Schr\"{o}dinger operators}, Lecture Notes in Mathematics, {\bf 1917}, Springer, Berlin, 2008.

\end{thebibliography}
\end{document}